\documentclass[sn-mathphys-num]{sn-jnl}

\usepackage{graphicx}
\usepackage{xcolor}
\usepackage[section]{placeins}
\usepackage{subcaption}
\usepackage{multirow}
\usepackage{multicol}
\usepackage{booktabs}
\usepackage{array}
\usepackage{graphbox}

\usepackage{additionals}
\usepackage[frozencache,cachedir=.]{minted2}

\definecolor{bg}{rgb}{.95, .95, .95}
\setminted{autogobble,
    breaklines,
    bgcolor=bg,
    samepage=true,
    python3=true,
}
\setmintedinline{breaklines=false}

\newcommand{\pyth}[1]{\mintinline{python}{#1}}

\newminted{python}{}

\usepackage{amsmath,amssymb}
\usepackage{amsthm}

\usepackage[utf8]{inputenc}
\usepackage[T1]{fontenc}
\usepackage{lmodern}

\usepackage[UKenglish]{babel}
\usepackage[babel=true]{microtype}

\usepackage[noabbrev]{cleveref}

\usepackage{xspace}

\newtheorem{rmk}{Remark}[section]
\newtheorem{note}[rmk]{Note}

\newtheorem{meth}{Method}

\newcommand{\R}{\mathbb{R}}

\newcommand{\domain}{\Omega}
\newcommand{\domainbnd}{\partial\domain{}}
\newcommand{\extdomain}{B}
\newcommand{\extdomainbnd}{\partial\extdomain{}}
\newcommand{\DD}{Diffuse Domain}

\newcommand{\BC}{\zeta(u, \phi, g)}

\newcommand{\DDMO}{\mbox{DDM1}\xspace}
\newcommand{\MixedZ}{\mbox{Mix0DDM}\xspace}
\newcommand{\MixedZDDM}{\mbox{Mix0DDM1}\xspace}
\newcommand{\MixedO}{\mbox{Mix1DDM}\xspace}
\newcommand{\NSDDM}{\mbox{NSDDM}\xspace}

\newcommand{\ddfem}{\mbox{\textsc{ddfem}}\xspace}
\newcommand{\dunefem}{\mbox{\textsc{Dune-Fem}}\xspace}

\newcommand{\dunefemdg}{\mbox{\textsc{Dune-Fem-DG}}\xspace}
\newcommand{\dolfindg}{\mbox{\textsc{DOLFIN\_DG}}\xspace}
\newcommand{\gmsh}{\mbox{Gmsh}\xspace}
\newcommand{\fenics}{\mbox{FEniCS}\xspace}
\newcommand{\UFL}{\mbox{\textsc{UFL}}\xspace}

\begin{document}
\title[DDFEM]{DDFEM: A Python Package for Diffuse Domain Methods}

\author*[1]{\fnm{Luke} \sur{Benfield}}\email{luke.benfield@warwick.ac.uk}
\author[1]{\fnm{Andreas} \sur{Dedner}}\email{a.s.dedner@warwick.ac.uk}

\affil*[1]{
    \orgdiv{Mathematics Institute},
    \orgname{University of Warwick},
    \orgaddress{\city{Coventry}, \postcode{CV4 7AL}, \country{UK}}
}
\abstract{
    Solving partial differential equations (PDEs) on complex domains can present significant computational challenges.
    The Diffuse Domain Method (DDM) is an alternative that
    reformulates the partial differential equations on a larger, simpler domain.
    The original geometry is embedded into the problem by representing it with a phase-field function.

    This paper introduces \ddfem{}, an extensible Python package
    to provide a framework for transforming PDEs into a \DD{} formulation.
    We aim to make the application of a variety of different \DD{} approaches
    more accessible and straightforward to use.

    The \ddfem{} package includes features to intuitively define complex domains
    by combining signed distance functions and provides a number of \DD{}
    transformers for general second-order evolution equations.
    In addition, we present a new approach for combining
    multiple boundary conditions of different types
    on distinct boundary segments.
    This is achieved by applying a normalised weighting,
    derived from the component signed distance functions,
    to combine the additional boundary terms in the \DD{} formulations
    without requiring additional phase fields.
    The domain definition and \DD{} transformation provided by our package
    are designed to be compatible with a wide range of existing finite element solvers
    without requiring code alterations.
    Both the original (non-linear) PDEs provided by the user and
    the resulting transformed PDEs on the extended domain are defined using
    the Unified Form Language (\UFL{}).
    Here we focus on the software implementation and
    demonstrate its effectiveness across a range of problems.
}

\maketitle

\section{Introduction}
\label{section:intro1}

Complex domains are essential for real world problems when solving partial differential equations (PDEs).
Traditional approaches rely on generating a fitted mesh of the geometry.
To capture the required detail,
the discretisation length is smaller than the problem's characteristic length.
Furthermore, polyhedra are unable to precisely match rounded domains.
Generating a sufficiently refined mesh is computationally expensive
(especially in higher dimensions),
and geometry that evolves in time can require repeated mesh generation/movement.

Many approaches have been developed to resolve these problems by approximating complex boundaries.
Typically, these methods instead require extending the problem
to a larger and simpler computational domain.
For example,
the Fictitious Domain Method \citep{Glowinski1996,Ramiere2007,Vos2008,Parussini2009}
which uses Lagrange multipliers to enforce boundary constraints,
and the Immersed Boundary Method \citep{LeVeque1994,Prenter2023,Griffith2020}
which adds Dirac delta functions for forcing terms.
Alternatively, other approaches modify the discretisation near the interface,
such as the Extended Finite Element Method \citep{Fries2010},
cut-cell methods \citep{Ji2006},
and the Ghost Fluid Method \citep{Macklin2008}.

The \DD{} method extends the original
computational domain into a much larger simple domain, making mesh generation trivial.
The original domain is then embedded using a phase-field function that approximates the domain's characteristic function,
creating a small diffuse interface layer.
The method was introduced by \citet{Kockelkoren2003},
and approximations for Dirichlet, Neumann, and Robin boundary conditions with asymptotic analysis have been developed \citep{Li2009, Lervaag2015}.
Furthermore, significant analysis has established convergence properties and error estimates \citep{Franz2012, Burger2015, Schlottbom2016}.
Also, the method has been extended to surface domains \citep{Raetz2006},
with analysis for complex coupled bulk-surface systems \citep{Teigen2009, Abels2015}.

One of the key advantages of the \DD{} method is
the ability to use existing PDE solver frameworks without
complex modifications to incorporate a boundary structure,
as only the PDE form itself changes with standard functions.
This approach has been demonstrated to be a useful method for many applications,
ranging from biological problems in bone deformation and tumour growth \citep{Aland2012, Nguyen2017, Nguyen2017a, Chen2014},
and materials science \citep{Chadwick2018, Raetz2016}.
Furthermore, it has been applied widely in multiphase fluid dynamics,
including Cahn-Hilliard Navier-Stokes systems \citep{Aland2010, Guo2021},
flows involving soluble surfactants \citep{Teigen2011, Teigen2009},
and other complex coupled problems such as flow in porous domains \citep{Bukac2023, Termuhlen2022}
and surface phase-field crystal models \citep{Aland2012a}.

There is a growing list of different approximations,
for example:
higher order approximations
\citep{Yu2020, Lervaag2015}, %
smoothed boundary method \citep{Yu2009, Termuhlen2022},
and derived from Nitsche's method \citep{Nguyen2017}.
Also, importantly, combining different boundary conditions requires careful consideration.
For example, \citet{Monte2022} introduces additional phase fields to
indicate local boundary conditions,
determined by solid angles extending from a centroid.
This is limited for complex boundaries;
for example, an arc shape has an external centroid,
which forces the inside and outside edges to group together.
Our method is an automatic approach that operates directly on
the signed distance function (SDF) formulation.
We construct domains through constructive solid geometry (CSG),
and obtain component SDFs which isolate boundary segments.
Therefore, we use a single phase field and a weighting approach.

However, evaluating existing research shows little development around
providing solver independent implementations that offer multiple \DD{} approaches.
Our package provides an automated translation of these
boundary constraints into a usable implementation,
whilst allowing users to easily swap between different \DD{} methods
or implement their own.
The derivation and convergence analysis of these methods
are studied in \citep{Benfield2025b},
whereas this paper focuses on
the software design and transformation framework.

Our
goal is to create an extensible Python package
that adapts existing PDE models to different \DD{} formulations.
To minimise library dependence,
we will only rely on the Unified Form Language (\UFL{}) module of \fenics{} \citep{Alnaes2014};
this is an established interface to define the PDE forms for
finite element approximation,
such as \dunefem{} \citep{Dedner2020}, \fenics{} \citep{Baratta2023}, and Firedrake \citep{Ham2023}.

The package will take a PDE model provided by the user using \UFL{}
and translate it into a \DD{} model description of the same format.
It is important that the documentation for each step is clear and accessible to the user,
allowing them to add new methods or change any attributes/parameters of a transformation.

In our package, we focus on a general nonlinear advection-diffusion problem in divergence form:
\begin{equation}
    \label{eqn:general}
    \partial_t U
    =
    - \nabla \cdot \left( F_c(U) - F_v(U,\nabla U) \right)
    + S_{i}(U,\nabla U) + S_{e}(U,\nabla U)
    \quad \text{ in } (0,T] \times \domain{}.
\end{equation}
The viscous flux is denoted by $F_v$,
the convective flux by $F_c$,
and two source terms $S_{i}$ and  $S_{e}$.
The interface used in this package follows that
introduced for the \dunefemdg{} module in \citep{Dedner2022}. There
a number of different implicit-explicit (IMEX) Runge-Kutta methods are implemented which
motivates the splitting of the source term $S=S_{i}+S_{e}$.
In \dunefemdg{} the diffusion $F_v$ and the source term $S_i$ are
treated implicitly while $S_e$ and the advection are treated explicitly.
All these terms may also depend on $x$ and time $t$.
Note that initial and boundary conditions are required to uniquely define
the solution.
We focus on the strong form of the PDE since
some \DD{} transformations are easier to carry out based on this form.
The same approach is followed in other existing work,
such as \dunefemdg{} \citep{Dedner2022} and \dolfindg{} \citep{Houston2018},
where the divergence form is preferred for Discontinuous Galerkin methods.
\begin{listing}[H]
    \begin{pythoncode}
class Model:
    def S_e(t, x, U, DU): ...
    def S_i(t, x, U, DU): ...
    def F_c(t, x, U): ...
    def F_v(t, x, U, DU): ...
    boundary = {...}
    \end{pythoncode}
    \caption{Abbreviated model API for \dunefemdg{}.}
    \label{lst:dunefemdg}
\end{listing}

The Python listings throughout will show individual components of \ddfem{} within the context of the library
and are not intended as standalone scripts.
Abbreviated implementations are identified in their captions or the surrounding text,
with \cref{app:code} containing instructions to obtain the complete code.

Our aim is to start with a minimal class describing the problem in \cref{eqn:general}
and transform it into a class with the same API
but now describing a PDE on a larger domain $\extdomain{}\supset\domain$.
As mentioned, we employ the same basic model API used in \dunefemdg{}.

This paper is organised as follows:
\Cref{section:ddm1} introduces
the fundamental ideas of the \DD{} method
through the construction of a simple problem.
\Cref{section:sdf} details how
multiple simple signed distance functions are combined using \pyth{ddfem.geometry}
to create complex domains.
\Cref{section:mixedboundaries} presents the \pyth{ddfem.boundary} module,
which is our approach to implementing multiple boundary conditions by using the construction hierarchy and weightings.
These components are combined in \cref{section:transformers}
to describe functions in \pyth{ddfem.transformers} to transform a fitted problem to a \DD{} approach.
Finally, \cref{section:experiments} demonstrates the complete \ddfem{} package.

\section{Simple \DD{} Method}
\label{section:ddm1}

\subsection{Dirichlet Boundary Condition}
\label{subsection:ddmethod}

We will explain how a \DD{} transformation is defined using the following example.
Consider a simple Dirichlet problem on a ball, $\domain{} = B_1(0)$:
\begin{subequations}
    \label{eq:standardAD}
    \begin{alignat}{2}
        - \nabla \cdot (D \nabla u) + \nabla \cdot (b u)  + c u & = f & \quad & \text{ in } \domain{},    \\
        u                                                       & = g & \quad & \text{ on } \domainbnd{},
    \end{alignat}
\end{subequations}
where the coefficients and boundary data are defined by
\begin{alignat*}{3}
    D(x) & = 0.1 - 0.05 |x|^2,                &
    b(x) & = 3 x^\perp,                       &
    c(x) & \equiv 0.01,                         \\
    f(x) & = 0.1 (1 - |x|^2),                 &
    g(x) & = \sin{(\pi x_0)} \sin{(\pi x_1)}.
\end{alignat*}
Therefore, we derive \cref{lst:ogModel}.
\begin{listing}[H]
    \begin{pythoncode}
import ufl
D = lambda x: 0.1 - 0.05 * ufl.dot(x, x)
b = lambda x: 3 * ufl.as_vector([x[1], -x[0]])
c = lambda x: 0.01
f = lambda x: 0.1 * ufl.as_vector([1 - ufl.dot(x, x)])
g = lambda x: ufl.as_vector(
    [ufl.sin(ufl.pi * x[0]) * ufl.sin(ufl.pi * x[1])]
)

class Model:
    def S_i(t, x, U, DU): return f(x) - c(x) * U
    def F_c(t, x, U): return ufl.outer(U, b(x))
    def F_v(t, x, U, DU): return D(x) * DU
    boundary = {(1,): lambda t, x: g(x)}
    \end{pythoncode}
    \caption{Model of original problem.}
    \label{lst:ogModel}
\end{listing}

\begin{rmk}
    The \ddfem{} API requires the time parameter, \pyth{t},
    as it follows from
    the \pyth{Model} format
    introduced in \cref{lst:dunefemdg}.
    The coefficients and boundary conditions may be time dependent for a given problem,
    as seen in
    \cref{sec:chemR}.
\end{rmk}

We use \pyth{ufl.as_vector} as the \ddfem{} package assumes all values are vector-valued,
avoiding the need for explicit handling of scalars and vectors.
To produce a \UFL{} form we can use the following code in \cref{lst:model2ufl},
based on a vector-valued $\R^1$ \pyth{space}.
\begin{listing}[H]
    \begin{pythoncode}
def model2ufl(Model, space, DirichletBC):
    u, v = ufl.TrialFunction(space), ufl.TestFunction(space)
    x = ufl.SpatialCoordinate(space)
    t = 0

    si = -Model.S_i(t, x, u, ufl.grad(u))
    fc = -Model.F_c(t, x, u)
    fv = Model.F_v(t, x, u, ufl.grad(u))
    form = ufl.inner(fv + fc, ufl.grad(v)) + ufl.inner(si, v)
    dbcs = [
        DirichletBC(space, value(t, x), region)
        for regions, value in Model.boundary.items()
        for region in regions
    ]
    return [form * ufl.dx == 0, *dbcs]
\end{pythoncode}
    \caption{Convert a model class to a \UFL{} form.}
    \label{lst:model2ufl}
\end{listing}
We have implemented the class \pyth{ddfem.model2ufl.DirichletBC} to store the Dirichlet data.
This follows the same constructor as \dunefem{} \citep{Dedner2020}.

To convert this problem into the \DD{} framework on the extended domain, $\extdomain{} \supset \domain{}$,
we need a phase-field function,
$\phi$, which approximates the characteristic function $\chi_\domain{}$:
\begin{equation}
    \phi(x)
    \approx \chi_{\domain{}}
    = \begin{cases}
        1 & x \in \domain{};                      \\
        0 & x \in \extdomain{}\setminus\domain{}.
    \end{cases}
\end{equation}
A common choice is
\begin{equation}
    \label{eq:phi}
    \phi(x) = \frac{1}{2} \left(1 - \tanh\left(\frac{3 r(x)}{\epsilon}\right)\right),
\end{equation}
where $r(x)$ denotes the signed distance function (SDF),
which is negative for $x \in \domain{}$
and positive for $x \in \extdomain{} \setminus \domain{}$,
and with $|\nabla r| = 1$
\citep{Li2009,Aland2010,Yu2020,Lervaag2015}.
Also, $\epsilon$ is a small parameter to determine the width of the interfacial region.
The PDE is then reformulated onto the larger,
regular domain with additional terms that approximate the boundary conditions.\\
From the example problem in \cref{eq:standardAD}, it is simple to define the SDF for
$\domain{} = B_1(0)$ as
\begin{equation}
    \label{eq:sdfBall}
    r(x) = |x| - 1.
\end{equation}
\begin{pythoncode}
def sdf(x):
    return ufl.sqrt(ufl.dot(x, x)) - 1
\end{pythoncode}
Importantly,
data functions like $f$ and $g$ in the problem in \cref{eq:standardAD}
are only defined on the boundary $\domainbnd{}$
or in the interior of $\domain{}$,
so they must be extended to $\extdomain{}$.
To do this, we will smoothly extend each function
such that
it is constant in the normal direction to the original boundary.
Note, we use an overbar notation to denote extended data and
coefficients on $\extdomain{}$.
For example, the extension of $g$ is defined as,
\begin{equation}
    \label{eq:gbar}
    \overline{g}(x) = g\left(x - r(x) \nabla r(x) \right),
\end{equation}
and the extension of $f$ is defined as,
\begin{equation}
    \label{eq:fbar}
    \overline{f}(x) = f\left(x - ( 1 - \chi_{\domain{}}) r(x) \nabla r(x) \right).
\end{equation}
To resolve the diffuse interface sufficiently
\citet{Yu2020} set $h = \epsilon / 4$,
and \citet{Li2009} use approximately 10 grid points.
We typically use 7 elements.
Therefore, let us define the core functions related to the SDF:
\begin{listing}[H]
    \begin{pythoncode}
def chi(x):
    return ufl.conditional(sdf(x) <= 0, 1, 0)

h = 0.01  # Example grid spacing
elements = 7
epsilon = h * elements / 2
def phi(x):
    return 0.5 * (1 - ufl.tanh(3 * sdf(x) / epsilon))

def projection(x):
    return -ufl.grad(sdf(x)) * sdf(x)

def boundary_projection(x):
    return x + projection(x)

def external_projection(x):
    return x + (1 - chi(x)) * projection(x)
    \end{pythoncode}
    \caption{Core SDF functions.}
    \label{fig:coreSDF}
\end{listing}

A very straightforward \DD{} method, originally introduced by \citet{Li2009},
is derived by extending the domain
in the variational formulation using the characteristic function and surface delta function,
\begin{align}
    \int_{\domain{}}{... \; dx}    & \mapsto \int_{\extdomain{}}{\chi_\domain{} ...\; dx}        \\
    \int_{\domainbnd{}}{... \; dS} & \mapsto \int_{\extdomain{}}{\delta_{\domainbnd{}} ...\; dx}
\end{align}
and adding a term to enforce the boundary condition
\begin{equation}
    \label{eq:ugOutside}
    u = \overline{g} \text{ on } \extdomain{} \setminus \domain{}.
\end{equation}

The strong formulation of \cref{eq:standardAD} for the simplest \DD{} method, DDM1,
is given by the following PDE,
\begin{equation}
    \label{eq:standardDDM1}
    - \nabla \cdot (\phi \overline{D} \nabla u)
    + \nabla \cdot ( \phi \overline{b} u)
    + \phi \overline{c} u
    + \BC{}
    = \phi \overline{f}
    \quad \text{ in } \extdomain{},
\end{equation}
with
\begin{equation}
    \label{eq:standardDDM1dirichlet}
    \BC{} = \frac{\overline{D}}{\epsilon^3}(1-\phi)(u-\overline{g}).
\end{equation}

As established in \citep{Li2009,Lervaag2015,Yu2020},
for a sufficiently extended domain,
the newly imposed boundary conditions on $\extdomainbnd{}$
will not impact the solution in $\domain{}$.
Furthermore,
the penalty term in \cref{eq:standardDDM1dirichlet} ensures that the solution
decays rapidly to $\overline{g}$ in $\extdomain{} \setminus \domain{}$,
meaning the choice of boundary condition on $\extdomainbnd{}$ has a negligible effect.
Consequently, we can naturally use $u = \overline{g}$ on $\extdomainbnd{}$.
Therefore, we get the \DD{} model shown in \cref{lst:DDMO}.
\begin{listing}[htbp]
    \begin{pythoncode}
bp = boundary_projection
ep = external_projection

class ModelDDM1:
    def S_i(t, x, U, DU):
        outFactor = D(ep(x))  # optional weighting

        source = phi(x) * (f(ep(x)) - c(ep(x)) * U)
        bc = ( -outFactor * (U - g(bp(x))) * (1 - phi(x)) /
               (epsilon**3) )
        return source + bc

    def F_c(t, x, U):
        return phi(x) * ufl.outer(U, b(ep(x)))

    def F_v(t, x, U, DU):
        return phi(x) * D(ep(x)) * DU

    boundary = {range(1, 5): lambda t, x: g(bp(x))}
\end{pythoncode}
    \caption{Model for \DDMO{} method.}
    \label{lst:DDMO}
\end{listing}

\subsection{Flux Boundary Condition}
\label{sec:fluxbc}

Alternatively, we take the original problem in \cref{eq:standardAD}
with flux boundary conditions.
We assume that $b u \cdot n = g_c$ and $D \nabla u \cdot n = g_v$
both hold on $\domainbnd{}$, so that
\begin{equation}
    (b u - D \nabla u) \cdot n = g_{c} - g_{v} \quad \text{on } \domainbnd{}.
\end{equation}
Then, we can again use the approach by \citet{Li2009}
and further expanded by \citet{Lervaag2015},
resulting in the same equation in \cref{eq:standardDDM1} but with
\begin{equation}
    \label{eq:standardDDM1neu}
    \BC{} = (\overline{g_{c}} - \overline{g_{v}}) |\nabla \phi|.
\end{equation}
This is derived from using the surface delta approximation
that $\delta_{\domainbnd{}} \approx |\nabla \phi|$.
As mentioned in the introduction, we assume the convective flux $g_{c}$ is taken explicitly,
and the diffusive flux $g_{v}$ is taken implicitly.
Therefore, the \DD{} model with flux boundaries would have the source terms
shown in \cref{lst:fluxbnd}.
\begin{listing}[htbp]
    \begin{pythoncode}
def S_i(t, x, U, DU):
    source = phi(x) * (f(ep(x)) - c(ep(x)) * U)
    g_phi = ufl.grad(phi(x))
    bc = g_v(bp(x)) * ufl.sqrt(ufl.dot(g_phi, g_phi))
    return source + bc

def S_e(t, x, U, DU):
    g_phi = ufl.grad(phi(x))
    bc = -g_c(bp(x)) * ufl.sqrt(ufl.dot(g_phi, g_phi))
    return  bc

ModelDDM1.S_i = S_i
ModelDDM1.S_e = S_e
ModelDDM1.boundary = {range(1,5): lambda t, x: ufl.zero(1)}
\end{pythoncode}
    \caption{Source term for model with flux boundary conditions.}
    \label{lst:fluxbnd}
\end{listing}
For this model, we have used zero Dirichlet boundary conditions on the extended domain.
As $\phi \approx 0$ and $|\nabla \phi| \approx 0$ near $\extdomainbnd{}$,
this choice does not affect the solution in $\domain{}$ provided
the domain is sufficiently extended.

This approach is straightforward to implement manually for a simple model;
however, significant challenges arise with increasing complexity of the
model or the \DD{} method used.
In this current model,
$\phi$ is simultaneously used in the PDE
and for controlling boundary conditions,
which is extended into the new domain using the SDF.
Introducing multiple or different types of boundary conditions
requires careful consideration of how to apply each separately to different regions.
Furthermore, defining a single analytical SDF for a complex domain is often practically impossible.

We provide the subpackages \pyth{ddfem.geometry} for composing simpler SDFs,
and \pyth{ddfem.boundary} to manage complex boundary conditions.
These tools provide a set of helpful functions
and classes that simplify the implementation of more complex \DD{}
approaches, which are made available in
\pyth{ddfem.transformers}.
\section{Signed Distance Functions}
\label{section:sdf}

It is not a trivial task for the user to provide the SDF
for the geometry of a complex domain.
We provide a \pyth{ddfem.geometry} subpackage
which relies on constructive solid geometry (CSG);
notably this is the same concept used in \gmsh{} \citep{Geuzaine2009}.

Existing libraries to generate SDFs using CSG
rely on computing values on a mesh or coordinate array;
however, our goal of transforming \UFL{} expressions cannot depend on a mesh.
While \citet{Persson2004,Quilez2020} provide standard formulations
for primitive SDFs and Boolean operations,
we implement these expressions directly in \UFL{}.
The user of the module can easily define their own shapes
by extending the \pyth{SDF} base class;
the only requirement is to implement
the \pyth{sdf} method containing the relevant function.
For example, we defined the previous ball domain SDF as:
\begin{pythoncode}
class Ball(SDF):
    def __init__( self, radius, center, epsilon=None, name=None,
                  *args, **kwargs):
        self.radius = radius
        self.center = center
        super().__init__(epsilon, name, *args, **kwargs)

    def sdf(self, x):
        center_x = x - self.center
        return ufl.sqrt(ufl.dot(center_x, center_x)) - self.radius
\end{pythoncode}
The other functions in \cref{fig:coreSDF}, e.g.,
\begin{multicols}{2}
    \begin{itemize}
        \item \pyth{chi}
        \item \pyth{phi}
        \item \pyth{boundary_projection}
        \item \pyth{external_projection}
    \end{itemize}
\end{multicols}
are provided by the base class, so they can easily be called on any SDF instance.
Note that the method \pyth{phi} requires that the attribute \pyth{SDF.epsilon} is set.

\begin{rmk}
    From here on, the term signed distance function (SDF) refers to a function
    that is a proper SDF only locally around the boundaries,
    relaxing the definition used previously.
\end{rmk}

Furthermore, to combine and modify these primitive functions,
we provide additional classes:
\begin{multicols}{3}
    \begin{itemize}
        \item \pyth{Intersection}
        \item \pyth{Invert}
        \item \pyth{Rotate}
        \item \pyth{Scale}
        \item \pyth{Subtraction}
        \item \pyth{Translate}
        \item \pyth{Union}
        \item \pyth{Xor}
        \item \pyth{Round}
        \item \pyth{Extrusion}
        \item \pyth{Revolution}
        \item[\vspace{\fill}]
    \end{itemize}
\end{multicols}
These are defined in the same way as the primitive shapes
by extending the \pyth{SDF} base class but through the parent class \pyth{BaseOperator(SDF)}.
This new class is used to transfer
the maximum \pyth{epsilon}.
Any operator the user requires can easily be added.
Furthermore, the method \pyth{search(child_name)}
is provided in the \pyth{SDF} base class, allowing easy traversal of the tree of SDFs
created using these operators via the name attribute.
Individual SDF classes can be constructed with separate $\epsilon$ attribute values
to be used when computing their \pyth{phi}.
Also, to simplify the common case of setting all to the same value,
\pyth{SDF.epsilon} is defined using a Python property attribute, so
the value will be propagated to the other composite SDFs.
For example, the \pyth{Union} class
simply stores the argument classes and uses the minimum to define the \pyth{sdf} method:
\begin{pythoncode}
class Union(BaseOperator):
    def __init__(
        self, sdf1, sdf2, epsilon=None, *args, **kwargs
    ):
        super().__init__(
            epsilon=epsilon, children=[sdf1, sdf2], *args, **kwargs
        )

    def sdf(self, x):
        return ufl.min_value(
            self.child_sdf[0].sdf(x), self.child_sdf[1].sdf(x)
        )
\end{pythoncode}

\begin{figure}[htbp]
    \centering
    \begin{subfigure}[t]{.49\textwidth}
        \centering
        \includegraphics{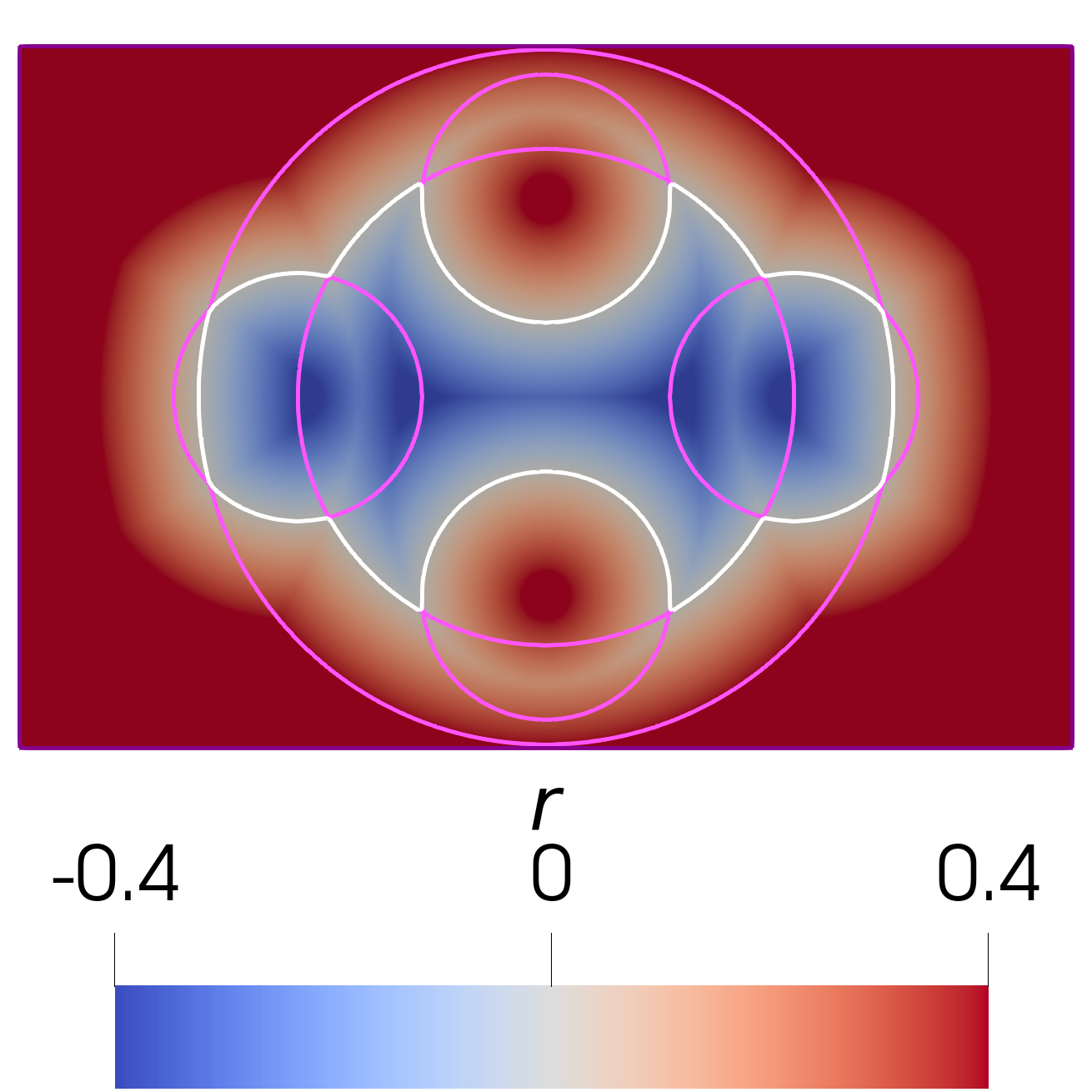}
        \caption{\pyth{omega(x)}}
    \end{subfigure}
    \begin{subfigure}[t]{.49\textwidth}
        \centering
        \includegraphics{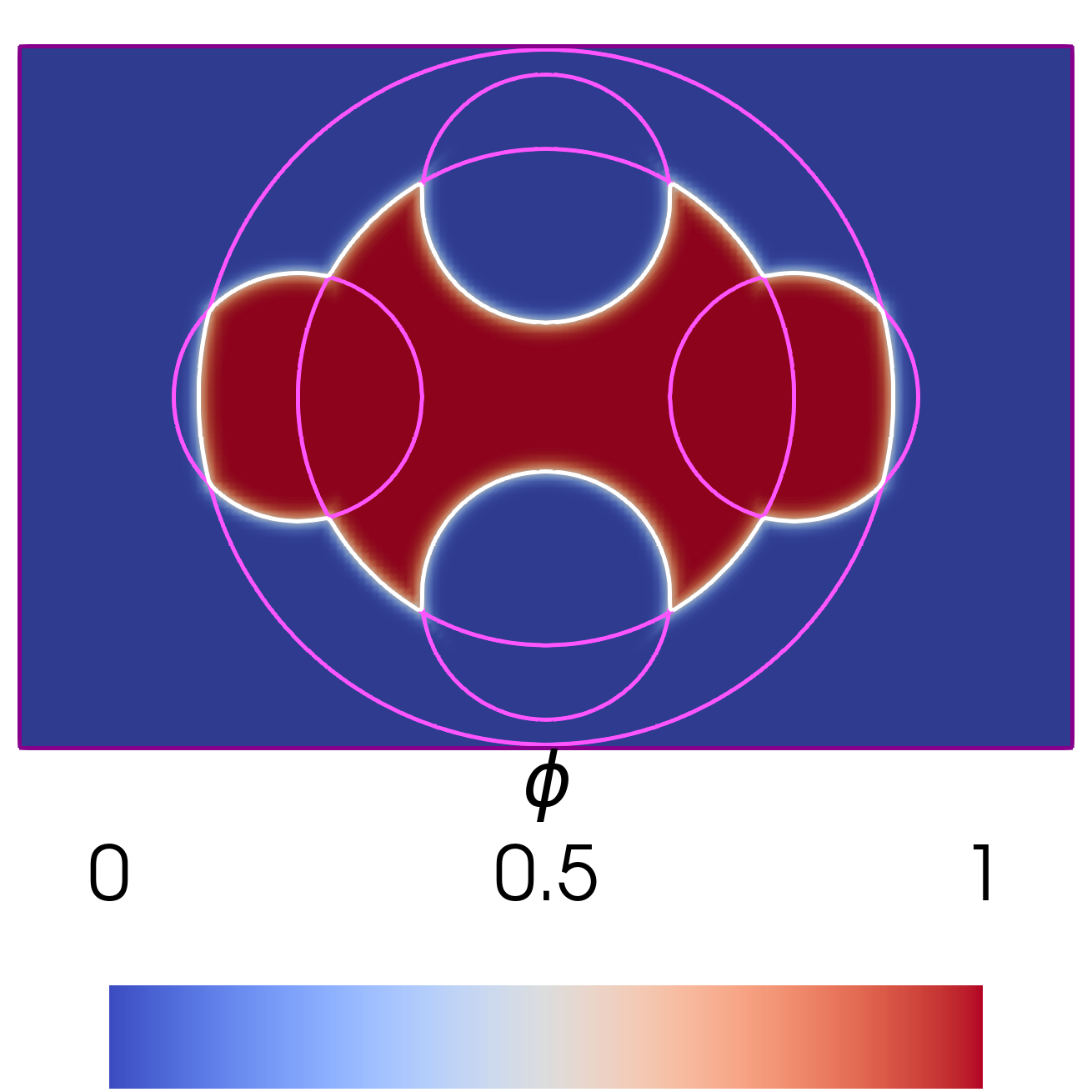}
        \caption{\pyth{omega.phi(x)}}
        \label{fig:circles_phi}
    \end{subfigure}
    \caption{Constructed SDF \pyth{omega} in \cref{lst:sdfsix}.}
    \label{fig:Balls}
\end{figure}

The \pyth{SDF} base class has methods implemented with the same name for each operation.
These several Boolean operations have also been used to overload the built-in Python operators.
For example, the following class shows the union and intersection operators, with other methods omitted:
\begin{pythoncode}
class SDF:
    def __or__(self, other):
        return Union(self, other)

    def __and__(self, other):
        return Intersection(self, other)
\end{pythoncode}
This allows easy manipulation of SDFs, such as in \cref{lst:sdfsix}, creating
the phase-field in \cref{fig:Balls}.
\begin{listing}[htbp]
    \begin{pythoncode}
sdfs = [
    Ball(radius=1, center=(0, 0), name="Center"),
    Ball(radius=0.5, center=(0, 0.8), name="TopCut"),
    Ball(radius=0.5, center=(0, -0.8), name="BotCut"),
    Ball(radius=0.5, center=(1, 0), name="RightAdd"),
    Ball(radius=0.5, center=(-1, 0), name="LeftAdd"),
    Ball(radius=1.4, center=(0, 0), name="Cutoff"),
]
bulk = sdfs[0] - sdfs[1] - sdfs[2] | sdfs[3] | sdfs[4]
bulk.name = "core"
omega = bulk & sdfs[5]
omega.name = "full"
omega.epsilon = 1e-3
    \end{pythoncode}
    \caption{Defining a six ball SDF.}
    \label{lst:sdfsix}
\end{listing}

Note that each \pyth{SDF} has an optional \pyth{name} attribute;
this can be defined during construction or assigned later.
We will see in \cref{section:mixedboundaries}
that this can be used to define boundary conditions
by associating the boundary \pyth{name}
with the boundary segment it defines in the final SDF.

\begin{figure}[htbp]
    \centering
    \includegraphics{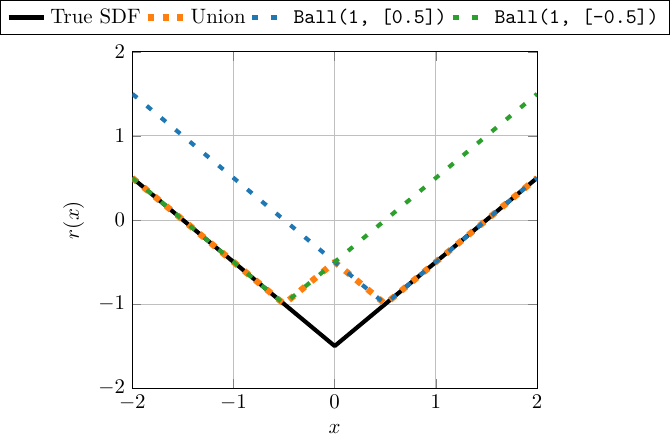}
    \caption{Errors in union of SDFs.}
    \label{fig:union}
\end{figure}

It is well known that Boolean operations using $\min$ and $\max$
do not in general produce a true SDF away from the interface \citep{Breen1998, Persson2004}.
However, they preserve the exact zero level-set defining the boundary $\domainbnd{}$.
This can be seen in \cref{fig:union};
this is a 1D example of two intersecting intervals.
We can see that the SDF for the intersecting region does not match the true SDF,
and the error increases away from the boundary.
The process of applying the \DD{} method requires the SDF for two purposes.
Firstly, we use it to generate a phase-field function, $\phi$.
As $\epsilon$ is small and the width of the interfacial region is approximately $2\epsilon$,
the impact of the imperfect SDF on $\phi$ is negligible and unobservable.
Secondly, the SDF is used in the extension of
the boundary functions and domain functions to the full computational domain $\extdomain{}$.
The \DD{} method does not require
the extension of these functions to be accurate far from the boundary.
\citet{Yu2020} states that the \DD{} solution only requires the extension to be at least $3\epsilon$ from the boundary.
Consequently, the error inherent to approximations from CSG operations
does not affect the \DD{} solution, provided $\epsilon$ is sufficiently small.

The \pyth{ddfem.geometry.Domain} class collects the functions related to the boundaries.
This acts similarly to the SDF class for $r_\domain{}$,
initialising with \pyth{Domain(omega)}.
The key methods implemented for the modification of the integrals are:
\begin{itemize}
    \item \pyth{phi}: $\phi(x)$,
    \item \pyth{scaled_normal}: $ - \nabla(\phi(x))$,
    \item \pyth{surface_delta}: $|\nabla(\phi(x))|$,
    \item \pyth{normal}: $ - \nabla(\phi(x)) / |\nabla(\phi(x))|$,
    \item \pyth{bndProjSDFs(SDFname)}:
          maps from the SDF $r_i$ to the unnormalised weight $w_i(x)$, \cref{eq:unnormweight},
          using its attribute \pyth{SDF.name}.
\end{itemize}
\section{Mixed Boundaries}
\label{section:mixedboundaries}

The preceding sections covered the implementation of boundary conditions
and the construction of SDFs for complex domains.
For general problems, it is critical to be able to combine multiple different boundary conditions.
Our approach achieves this by partitioning the extended domain into sections,
and restricting the boundary terms to the boundary neighbourhood.
However, there is currently a limitation for vector-valued problems.
Each component is restricted to the same type of boundary condition when defined by the same model class.
This requires combining multiple classes to split the system of equations.

We assume the boundary conditions are defined on segments $\Gamma_i$
such that
\begin{equation}
    \domainbnd{} = \bigcup_{i \in I} \Gamma_{i}.
\end{equation}
Furthermore, the different boundaries can be defined using SDFs,
\begin{equation}
    \label{eq:boundaryGamma}
    \Gamma_i = \{ x : r_i(x) = 0 \} \cap \{ x : r_\domain{}(x) = 0 \},
\end{equation}
where the SDFs $r_i$ were combined via
Boolean operators to compose the domain $\domain$.
It is important
that the zero level sets defined by the individual SDFs do not
intersect on $(d-1)$-dimensional sets,
and can only intersect on lower-dimensional sets.
Therefore, the $(d-1)$-dimensional measure of their intersection, denoted by $|\cdot|_{d-1}$, is zero,
\begin{equation}
    \label{eq:boundaryPoints}
    | \{x : r_i(x) = 0 \} \cap \{ x : r_j(x) = 0 \} |_{d-1} = 0
    \quad \forall i,j \in I : i \neq j.
\end{equation}

Given the problem in \cref{eqn:general} with boundary conditions
\begin{subequations}
    \label{eq:boundaryMultiple}
    \begin{alignat}{3}
        u                          & = g_{i}               & \quad & \text{on } \Gamma_{i} & \quad & \forall i \in I_V, \\
        (b u - D \nabla u) \cdot n & = g_{c, i} - g_{v, i} & \quad & \text{on } \Gamma_{i} & \quad & \forall i \in I_F,
    \end{alignat}
\end{subequations}
where the index sets $I_V$ and $I_F$ partition $I$ into the Dirichlet and flux boundary segments.
We need to combine
the additional boundary modification terms
from \cref{eq:standardDDM1dirichlet,eq:standardDDM1neu}
to define a term ${\BC{}}$ combining all boundary conditions.

Based on \cref{eq:standardDDM1dirichlet,eq:standardDDM1neu},
the boundary conditions lead to the following forcing terms in the \DD{} method:
\begin{subequations}
    \label{eq:BCterms}
    \begin{alignat}{2}
        {\BC{}}_i & = \frac{1}{\epsilon^3}(1-\phi)(u-\overline{g_{i}})          & \quad & \forall i \in I_V, \\
        {\BC{}}_i & = (\overline{g_{c, i}} - \overline{g_{v, i}}) |\nabla \phi| & \quad & \forall i \in I_F.
    \end{alignat}
\end{subequations}
To restrict these forcing terms to their specific regions,
we combine them using a weighted sum:
\begin{equation}
    {\BC{}} = \sum_{i \in I_V \cup I_F} \omega_{i} {\BC{}}_i.
\end{equation}
In the following,
we discuss a practical approach for defining the weights $\omega_i$.

Ideally, the weighting should satisfy the condition
\begin{equation}
    \label{eq:idealweight}
    \omega_{i}
    =
    \begin{cases}
        1 & \text{if } x \text{ is closest to } \Gamma_i,           \\
        0 & \text{if } x \text{ is closest to } \Gamma_j, j \neq i.
    \end{cases}
\end{equation}
The constructed boundary terms must remain valid \UFL{} expressions,
so that they can be integrated into the weak form.
We approximate the weighting function by comparing
the component SDF, $r_i(x)$,
to the global SDF, $r_\domain{}(x)$.
However, due to certain Boolean operations,
such as difference and XOR,
requiring inverted SDFs ($-r$) during their computation,
we compare the absolute values.
To ensure numerical stability,
we relax the exact equality using a small tolerance $\vartheta>0$ (in our
implementation $\vartheta=10^{-10}$) to define the weights as
\begin{equation}
    \label{eq:unnormweight}
    w_{i,\domain{}}(x)
    =
    \begin{cases}
        1 & \text{if } \big| |r_\domain{}(x)| - |r_i(x)| \big| < \vartheta, \\
        0 & \text{otherwise}.
    \end{cases}
\end{equation}

Following \cref{eq:boundaryPoints},
the zero level sets of the SDFs only intersect on
lower dimensional sets of measure zero. %
Consequently, for almost all $x \in \domain$
(excluding regions affected by the tolerance),
the closest boundary is unique,
meaning that
for a given $i$ when $w_{i, \domain{}} = 1$,
we have $w_{k, \domain{}} = 0$ for all $k \neq i$.
This is a desired property as
the solution in the neighbourhood of
a boundary segment is only influenced by its associated boundary condition.

However, near the intersecting points of different SDFs
(as well as regions from any imperfect Boolean operations),
the \DD{} approach results in the summation of multiple extensions of different boundary conditions.
This can cause a significant error.
We normalise the weights by setting
\begin{equation}
    \omega_i = \frac{w_{i,\domain{}}}{\sum_k w_{k,\domain{}}}.
\end{equation}
Note, since every point matches at least one SDF, a division by zero cannot occur.
Therefore, we have $\sum_i \omega_i = 1$,
and we average the ${\BC{}}_i$ terms near these intersection points.

\begin{figure}[ht]
    \centering
    \begin{subfigure}[t]{.24\textwidth}
        \includegraphics[width=\textwidth]{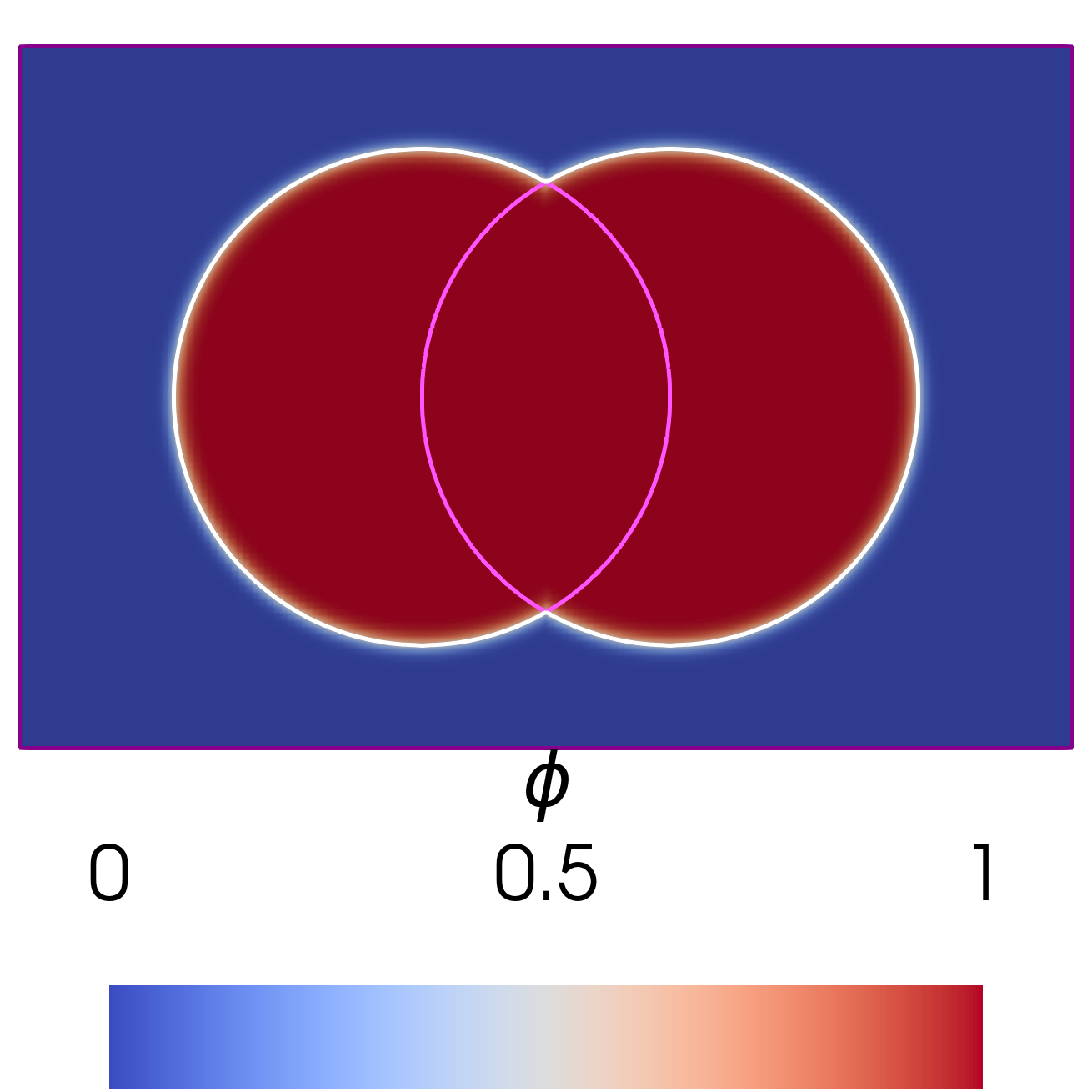}
        \\
        \includegraphics[width=\textwidth]{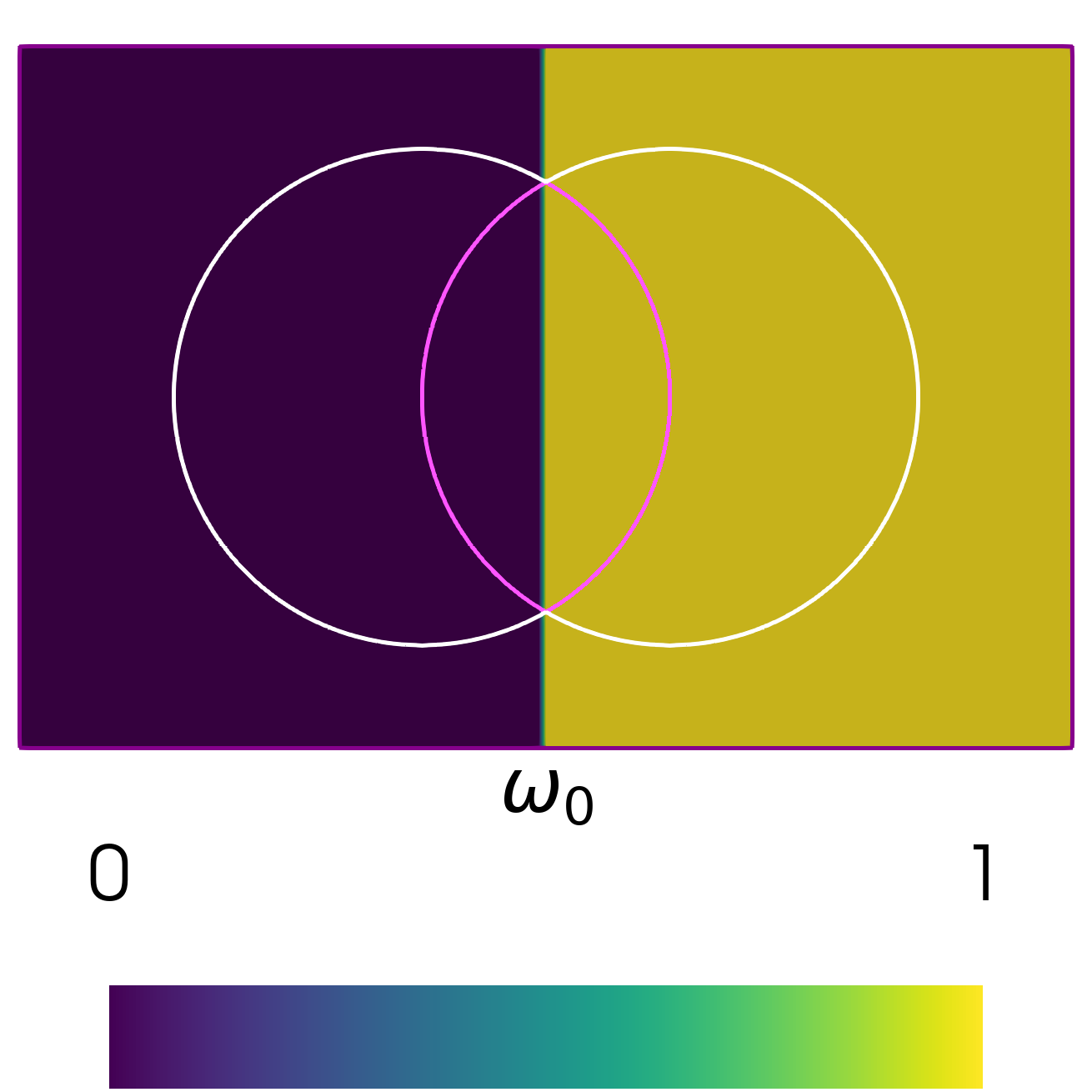}
        \caption{Union}
        \label{fig:sdfunion}
    \end{subfigure}
    \begin{subfigure}[t]{.24\textwidth}
        \includegraphics[width=\textwidth]{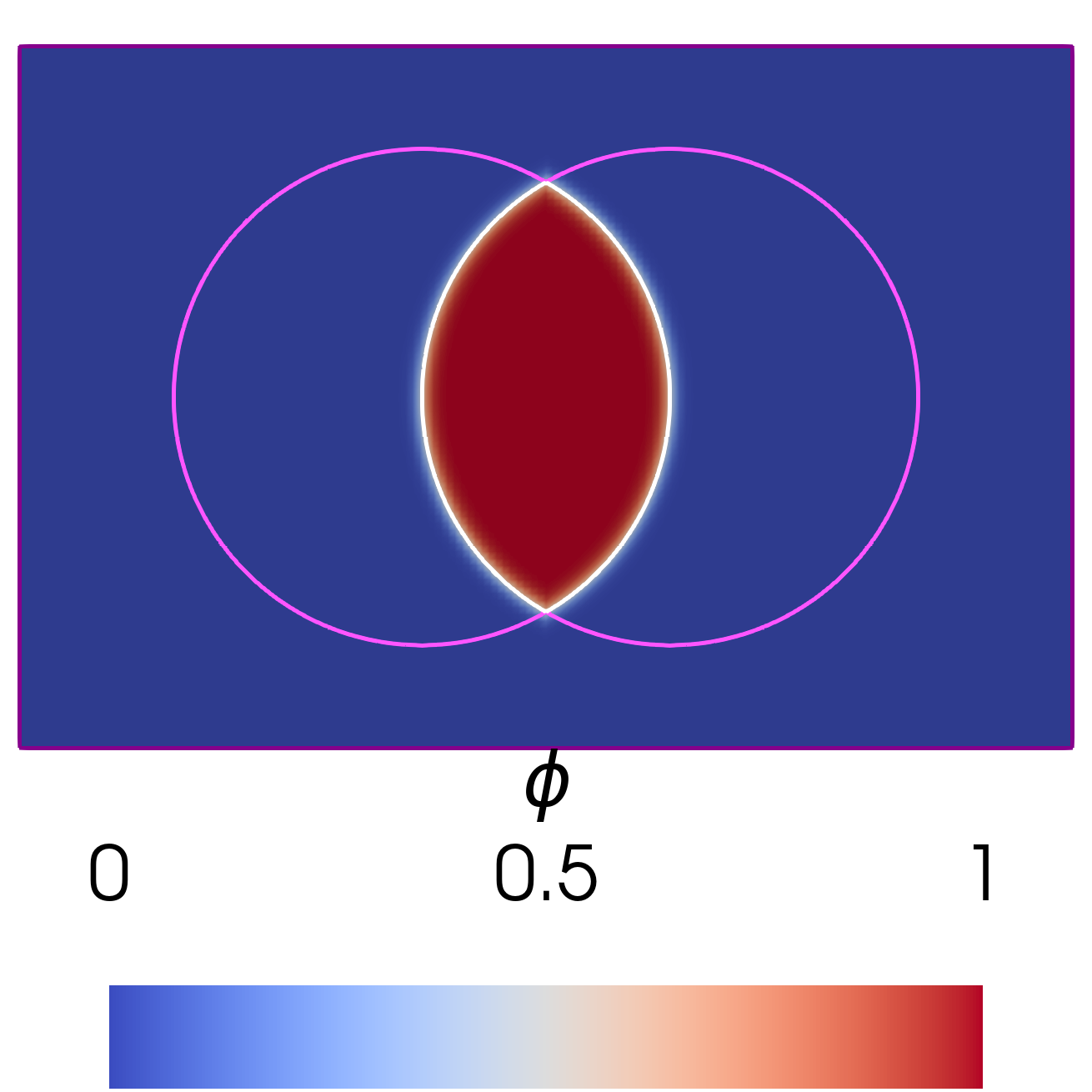}
        \\
        \includegraphics[width=\textwidth]{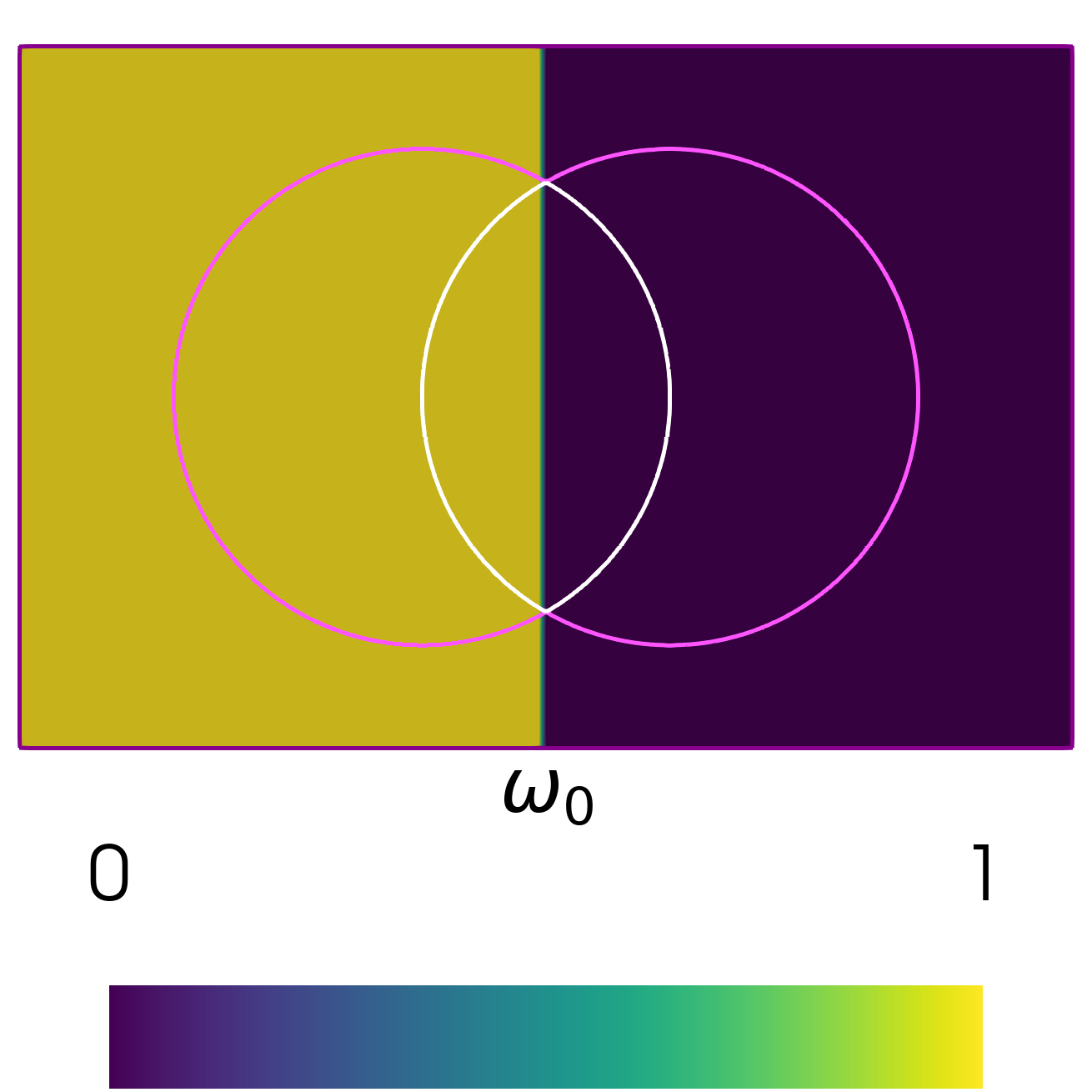}
        \caption{Intersection}
        \label{fig:sdfintersection}
    \end{subfigure}
    \begin{subfigure}[t]{.24\textwidth}
        \includegraphics[width=\textwidth]{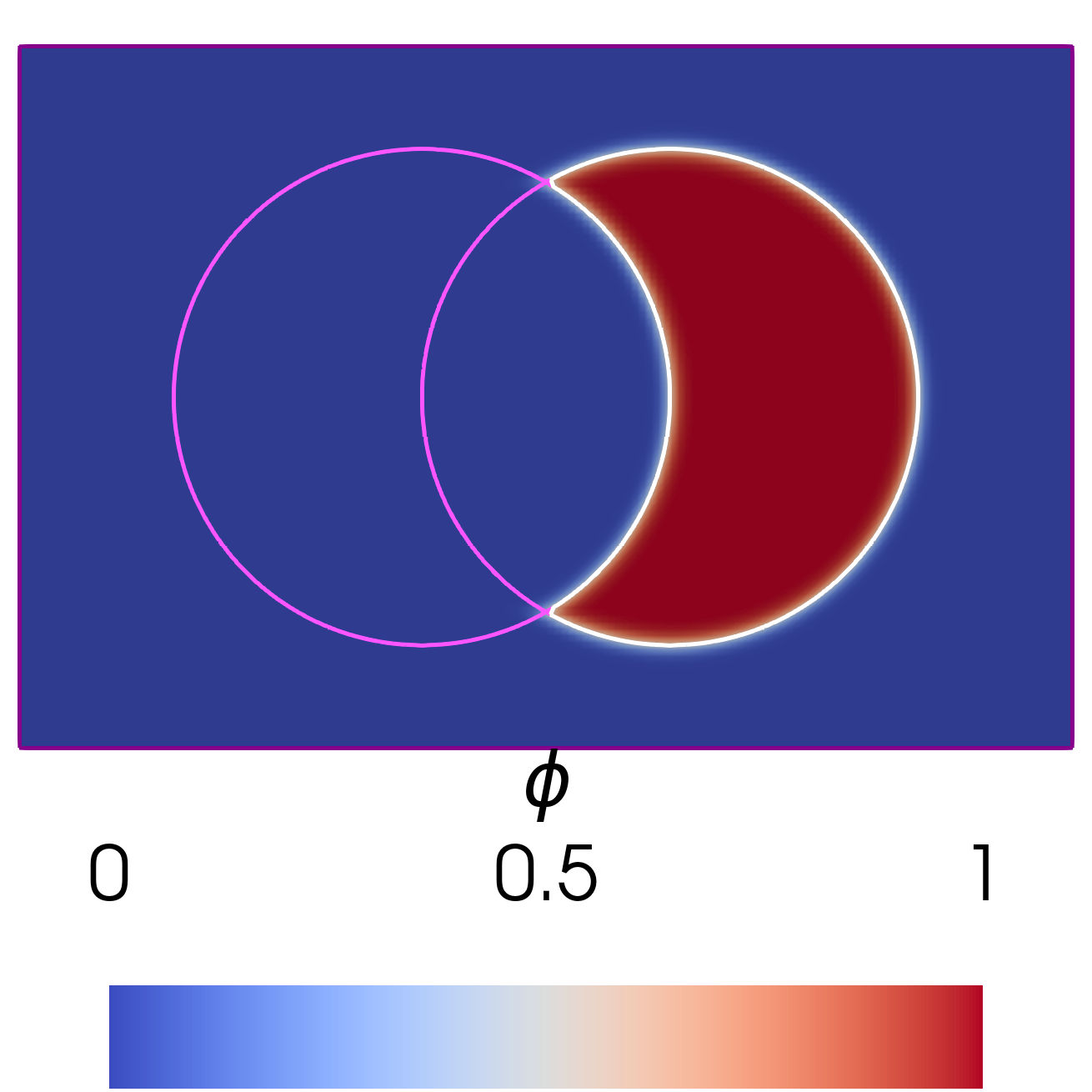}
        \\
        \includegraphics[width=\textwidth]{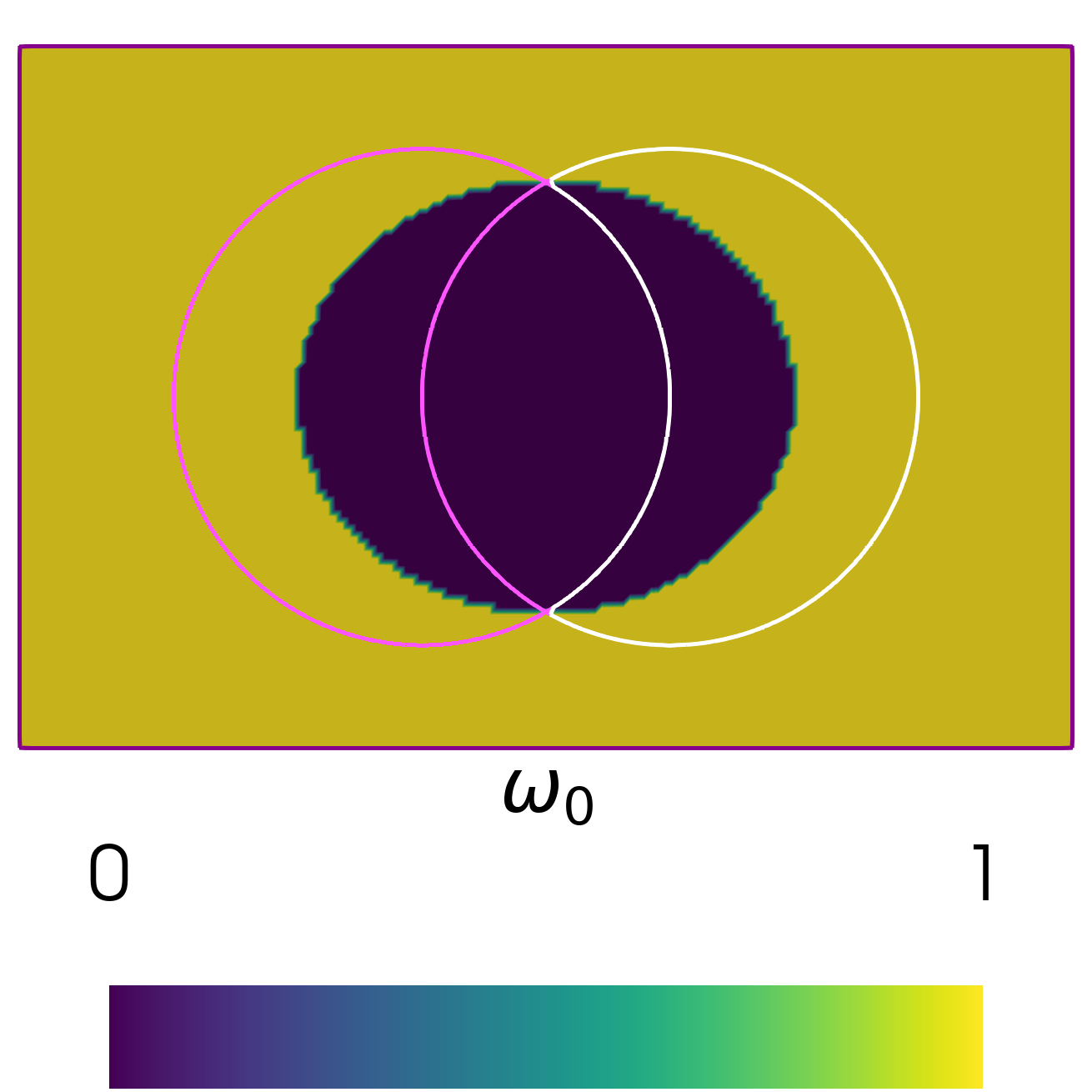}
        \caption{Difference}
        \label{fig:sdfsubtraction}
    \end{subfigure}
    \begin{subfigure}[t]{.24\textwidth}
        \includegraphics[width=\textwidth]{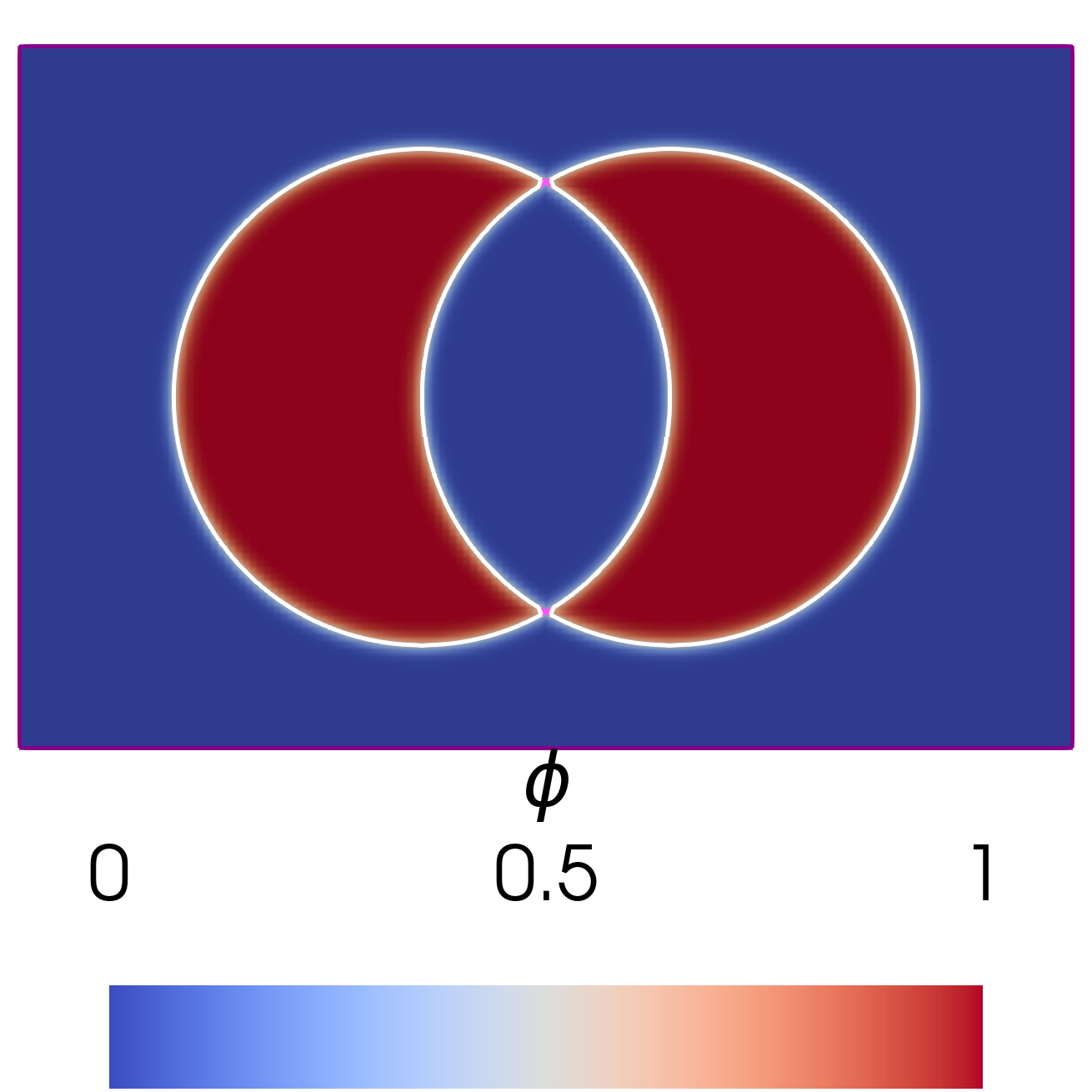}
        \\
        \includegraphics[width=\textwidth]{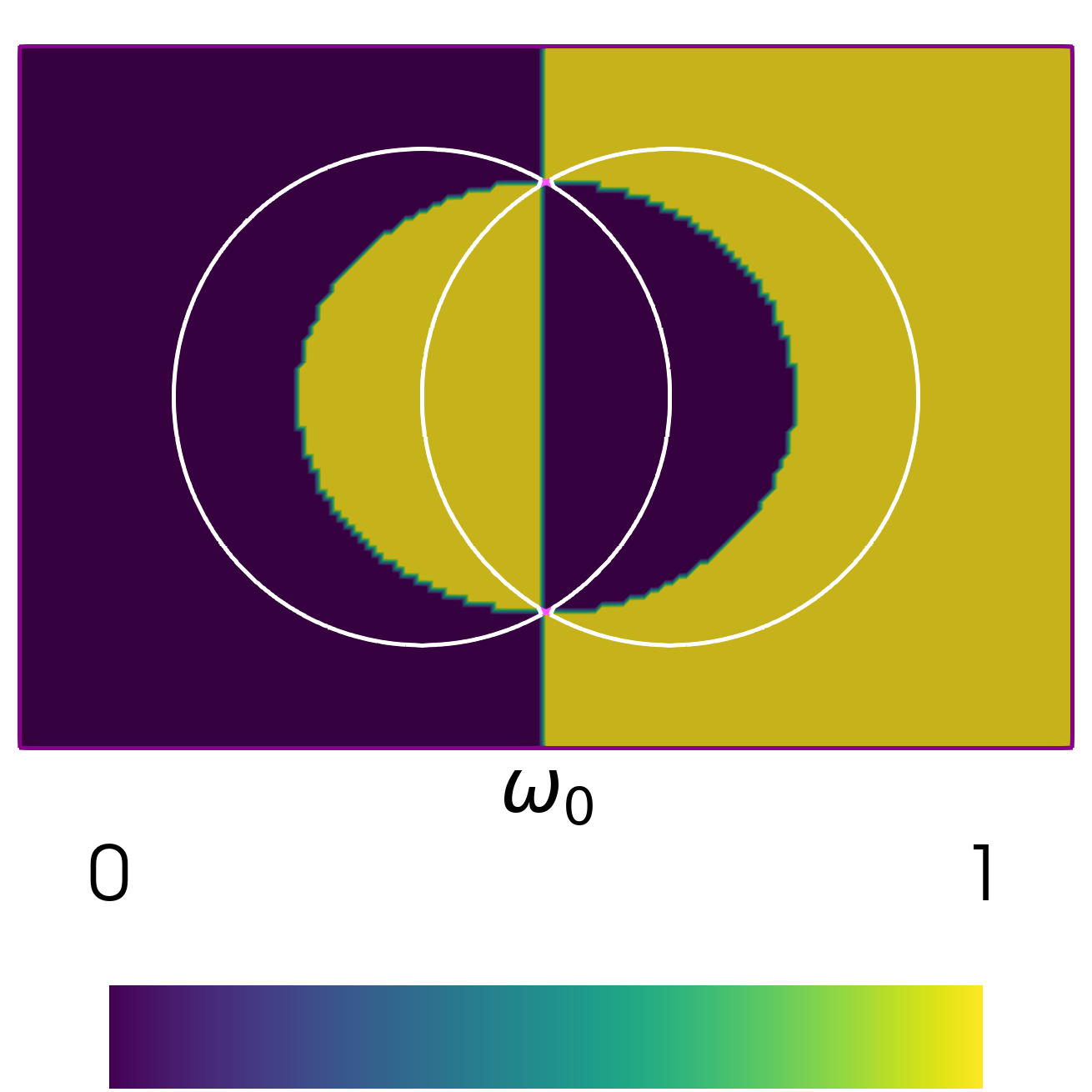}
        \caption{XOR}
        \label{fig:sdfxor}
    \end{subfigure}
    \caption{SDFs for the union, intersection, difference, and XOR operations on two balls. The second shows the boundary weight $\omega_0$ for the ball on the right.}
\end{figure}

Using the four binary operations in \pyth{ddfem.geometry} to combine two SDFs,
we can easily plot the different weightings.
In these examples we will use two balls:
\begin{pythoncode}
c0 = gm.Ball(radius=1, center=(0.5, 0), name="Ball0")
c1 = gm.Ball(radius=1, center=(-0.5, 0), name="Ball1")
\end{pythoncode}
The union, \pyth{c0 | c1},
is shown in \cref{fig:sdfunion}.
The lower plot shows the boundary condition weighting, $\omega_0$,
corresponding to the right ball \pyth{"Ball0"}.
We can see that the different regions are split correctly in a
neighbourhood of each ball.
Similarly,
the intersection, \pyth{c0 & c1},
is shown in \cref{fig:sdfintersection}.
The difference, \pyth{c0 - c1},
is shown in \cref{fig:sdfsubtraction}.
This shows the result with an imperfect SDF away from the boundary
(as seen from $\omega_0=1$ on the left side, outside the domain).
The XOR, \pyth{c0 ^ c1}, is shown in \cref{fig:sdfxor}.

Recall that $g_i$ needs to be defined in the neighbourhood of $\Gamma_i$
for the additional terms from \cref{eq:BCterms},
so we also extended the boundary values in
\cref{eq:boundaryMultiple}
using the same projection
from \cref{eq:gbar}.
In summary, we provide the following functions for the different boundary types:
\begin{align}
    G_{V}   &
    = \sum_{i \in I_{V}} \omega_i \overline{g_{i}},                  \\
    G_{F_v} &
    = \sum_{i \in I_{F}} \omega_i \overline{g_{v, i}} |\nabla \phi|, \\
    G_{F_c} &
    = \sum_{i \in I_{F}} \omega_i \overline{g_{c, i}} |\nabla \phi|.
\end{align}
These boundary equations are collected into the \pyth{ddfem.boundary.BoundaryTerms} class,
which is initialised using a given \pyth{Model} and \pyth{SDF}.
The class automatically constructs an instance of
\pyth{ddfem.geometry.Domain} internally;
however, a pre-modified instance can be provided for advanced applications.
This interface simplifies the implementation of
boundary conditions within a \DD{} model
and extends the boundary data for compatibility with the \DD{} approximation.
It provides the following key methods:
\begin{itemize}
    \item \pyth{BndValueExt}:
          returns $G_{V}$, or \pyth{None} if no Dirichlet boundaries.
    \item \pyth{jumpV}:
          returns a weighted sum with $u$, $\sum_{i \in I_{V}} \omega_i (u - \overline{g_{i}})$.
    \item \pyth{BndFlux_vExt}:
          returns $G_{F_v}$, or \pyth{None} if no flux boundaries.
    \item \pyth{jumpFv}:
          returns a weighted sum with $F_v$, $\sum_{i \in I_{F}} \omega_i (F_v n - \overline{g_{v, i}})|\nabla \phi|$.
    \item \pyth{BndFlux_cExt}:
          returns $G_{F_c}$.
    \item \pyth{jumpFc}:
          returns a weighted sum with $F_c$, $\sum_{i \in I_{F}} \omega_i (F_c n - \overline{g_{c, i}})|\nabla \phi|$.
\end{itemize}

The \pyth{Model.boundary} attribute is a dictionary,
and its keys determine how the boundary conditions are handled.
A key is assumed to be a boundary condition for the diffuse interface when
it is either an instance of the \pyth{SDF} class,
or the string matching the \pyth{SDF.name} attribute of the associated instance.
Other keys are assumed to be boundary conditions for the computational mesh
and will be left unchanged by the \DD{} transformer functions.
The corresponding value in the dictionary is the boundary function,
e.g., $g$, $g_c$, $g_v$.
For example, using the SDF from \cref{fig:Balls,lst:sdfsix} we can set the boundary conditions:
\begin{pythoncode}
Model.boundary = {
    bulk: [BndFlux_c(g_c), BndFlux_v(g_v)],
    "Cutoff": BndValue(g),
}
\end{pythoncode}
Note that the implementation requires each dictionary value to be one of the following types:
\begin{multicols}{2}
    \begin{itemize}
        \item \pyth{BndValue}
        \item \pyth{BndFlux_c}
        \item \pyth{BndFlux_v}
        \item \pyth{[BndFlux_c, BndFlux_v]}
    \end{itemize}
\end{multicols}
The required boundary flux
(\pyth{BndFlux_c}, \pyth{BndFlux_v}, or the list of both)
corresponds to whether the source terms $F_c$ and $F_v$ are defined in the problem.
This wrapper allows easy identification of the different boundary types.
We finally provide the function
\pyth{ddfem.model2ufl}
which takes a model class and its boundary dictionary,
converting it to a full \UFL{} form.

Using the above we arrive at the following updated model for the \DDMO{} approach:
\begin{listing}[H]
    \begin{pythoncode}
class ModelDDM1:
    BT = BoundaryTerms(Model, omega)
    boundary = BT.physical

    bp = BT.domain.boundary_projection
    ep = BT.domain.external_projection

    def S_e(t, x, U, DU):
        return -BT.BndFlux_cExt(t, x, U)

    def S_i(t, x, U, DU):
        source = phi(x) * (f(ep(x)) - c(ep(x)) * U)
        bc_v = -BT.jumpV(t, x, U) * (1 - phi(x)) / (epsilon**3)
        bc_fv = BT.BndFlux_vExt(t, x, U, DU)
        return source + bc_v + bc_fv

    def F_c(t, x, U):
        return phi(x) * ufl.outer(U, b(ep(x)))

    def F_v(t, x, U, DU):
        return phi(x) * D(ep(x)) * DU
    \end{pythoncode}
    \caption{Model of \DDMO{} approximation.}
    \label{lst:DDM1Model}
\end{listing}
Then the completed \UFL{} form is obtained using:
\begin{pythoncode}
ufl_form = model2ufl(ModelDDM1, space)
\end{pythoncode}
\section{Transformers}
\label{section:transformers}

Our main goal is to simplify implementing and using different \DD{} approaches,
so we now introduce the \pyth{transformers} subpackage.
This implements a group of new functions
to transform an existing \pyth{Model} (\cref{lst:ogModel})
and a \pyth{SDF} class,
returning a new model class defined on the domain $\extdomain{}$ based on a wide range of \DD{} methods.

First, any transformer should call the function \pyth{pretransformer} to return a new model class.
This will construct the \pyth{BoundaryTerms} class,
to set up the extended boundaries and additional boundary terms.
To improve performance,
we have found it beneficial to include the extra default boundary condition:
\begin{pythoncode}
boundary[lambda x: domain.omega.chi(x) < 0.5] = BndValue(
    BT.BndValueExt
)
\end{pythoncode}
On any boundaries in the extended domain,
$\extdomainbnd{} \setminus \domainbnd$,
we use the extended value of the Dirichlet boundary conditions.
However, if only flux boundary conditions are given,
we use the extended value of the flux boundary conditions:
\begin{pythoncode}
valFc = BndFlux_c(lambda t, x, U, n: -BT.BndFlux_cExt(t, x, U))
valFv = BndFlux_v(
    lambda t, x, U, DU, n: BT.BndFlux_vExt(t, x, U, DU)
)
boundary[lambda x: domain.omega.chi(x) < 0.5] = [valFc, valFv]
\end{pythoncode}
If the user provides boundary conditions for $\extdomainbnd{}$ in \pyth{Model.boundary},
these take priority and this default extension is not added.
These default conditions assume that $\domain{}$ is fully embedded within $\extdomain{}$.

\begin{rmk}
    It is the user's responsibility to ensure that diffuse interfaces do not intersect $\extdomainbnd{}$.
    This conflict is not handled and will result in both the boundary conditions on $\extdomainbnd{}$
    and the \DD{} penalty terms being added with their respective data.

    For problems combining physical boundaries with diffuse boundaries,
    the internal diffuse interfaces should be specified with an \pyth{SDF},
    whilst external physical boundaries should be prescribed directly with no \pyth{SDF} for those edges.
\end{rmk}

A further important role of
\pyth{ddfem.transformers.pretransformer}
is to extend all the terms and coefficient functions.
This redefines all methods
(e.g., \pyth{S_i}, \pyth{F_v})
using \pyth{external_projection}
to simplify implementing transformers.

A transformer class, \pyth{DDModel},
is implemented using a similar structure to \cref{lst:DDM1Model}.
In \cref{section:mixedboundaries},
we saw that while the viscous $F_v$
and convective $F_c$ flux terms are simply multiplied by $\phi$,
the source terms are augmented by new boundary and potential stability terms.
This requires implementing separate methods
to split the different components of the source terms.
Finally, calling the \pyth{ddfem.transformers.posttransformer} function with \pyth{DDModel}
constructs only the required methods to be compatible with \cref{lst:dunefemdg}
and ensures only existing methods are available in the new model.
The effects of \pyth{posttransformer} on different methods are shown in \cref{tab:ddmodel_methods}.
\begin{table}[htbp]
    \centering
    \begin{tabular}{lll}
        \toprule
        {\pyth{DDModel} Method}             & {Methods Included} & {\pyth{Model} requirement}                 \\ \midrule
        \multirow{3}{*}{\pyth{DDModel.S_e}} & \pyth{S_e_source}  & \pyth{S_e}                                 \\
                                                                 & \pyth{S_e_convection} & \pyth{F_c}         \\
                                                                 & \pyth{S_outside}      & \pyth{outFactor_e} \\ \midrule
        \multirow{3}{*}{\pyth{DDModel.S_i}} & \pyth{S_i_source}  & \pyth{S_i}                                 \\
                                                                 & \pyth{S_i_diffusion}  & \pyth{F_v}         \\
                                                                 & \pyth{S_outside}      & \pyth{outFactor_i} \\ \midrule
        \pyth{DDModel.F_c}                  & \pyth{F_c}         & \pyth{F_c}                                 \\ \midrule
        \pyth{DDModel.F_v}                  & \pyth{F_v}         & \pyth{F_v}                                 \\ \bottomrule
    \end{tabular}
    \caption{Composition of methods in the transformed \pyth{DDModel} class.}
    \label{tab:ddmodel_methods}
\end{table}

The attributes \pyth{outFactor_e} and \pyth{outFactor_i} are used to
define the scaling factors for \pyth{S_outside}.
The key idea is that the penalty term should scale
in proportion to the scale of the problem
(similarly to $\overline{D}$ in \cref{eq:standardDDM1dirichlet}).
These attributes provide separate scaling for the advection and diffusion components.
A default option,
using \UFL{}'s differentiation to evaluate the magnitude of the fluxes,
is implemented to reduce manual configuration.
Note, this automatic penalty scaling is still in testing.
The user can assign these attributes to override the default values.

The decorator \pyth{ddfem.transformers.transformer} is implemented to simplify the use
of \pyth{pretransformer} and \pyth{posttransformer}.
A complete implementation of the \DDMO{} transformer is given in \cref{lst:transformer}.
\begin{listing}[htbp]
    \begin{pythoncode}
@ddfem.transformers.transformer
def DDM1(Model):
    class DDModel(Model):
        def S_e_source(t, x, U, DU):
            return DDModel.phi(x) * Model.S_e(t, x, U, DU)
        def S_e_convection(t, x, U, DU):
            return -DDModel.BT.BndFlux_cExt(t, x, U)
        def S_outside(t, x, U, DU):
            return -( DDModel.BT.jumpV(t, x, U) *
                      (1 - DDModel.phi(x)) / (DDModel.epsilon**3) )
        def S_i_source(t, x, U, DU):
            return DDModel.phi(x) * Model.S_i(t, x, U, DU)
        def S_i_diffusion(t, x, U, DU):
            if DDModel.BT.BndFlux_vExt is not None:
                diffusion = DDModel.BT.BndFlux_vExt(t, x, U, DU)
            else:
                diffusion = zero(U.ufl_shape)
            return diffusion
        def F_c(t, x, U):
            return DDModel.phi(x) * Model.F_c(t, x, U)
        def F_v(t, x, U, DU):
            return DDModel.phi(x) * Model.F_v(t, x, U, DU)
    return DDModel
\end{pythoncode}
    \caption{Usage of the transformer decorator.}
    \label{lst:transformer}
\end{listing}

The transformed model can be used with \pyth{model2ufl}, and
the returned \UFL{} form can then be used with any solver
and mesh defined on $\extdomain{}$,
or further modified.
Currently, the package includes the following transformers:
\begin{multicols}{3}
    \begin{itemize}
        \item \pyth{DDM1 (MDD1)}
        \item \pyth{DDM2 (MDD2)}
        \item \pyth{DDM3 (MDD3)}
        \item \pyth{Mix0DDM}
        \item \pyth{Mix1DDM}
        \item \pyth{SBM}
        \item \pyth{NDDM}
        \item \pyth{NSDDM}
        \item \pyth{Fitted}
    \end{itemize}
\end{multicols}
Furthermore, the transformers are implemented to allow time-dependent coefficients.
Examples of using these transformer functions will be shown in \cref{section:experiments}.
\section{Demonstration}
\label{section:experiments}

We conclude by demonstrating the \DD{} transformers over a range of different problems to illustrate the implementation.
We have already seen \DDMO{}, but we will also make comparisons to \MixedZ{} and \NSDDM{} derived and analysed for Dirichlet boundary conditions in \citep{Benfield2025b}. But we should note that \MixedZ{} was designed with a focus on Dirichlet boundary conditions and improved gradient approximation, and not for full flux boundaries, while our experiments indicate that \DDMO{} works well in this scenario. This motivates the introduction of a new method \MixedZDDM{} which blends \MixedZ{} and \DDMO{} as discussed in \cref{sec:poisson} (see also the code snippet in \cref{app:poisson}). The detailed derivation, analysis and numerical investigation of all these methods is outside the scope of this work,
and can be found in the literature (see for example \citep{Benfield2025b} and the references therein).
For completeness, their formulations are provided in \cref{sec:appMethods}.

\begin{figure}[htbp]
    \centering
    \begin{subfigure}[t]{\textwidth}
        \centering
        \includegraphics{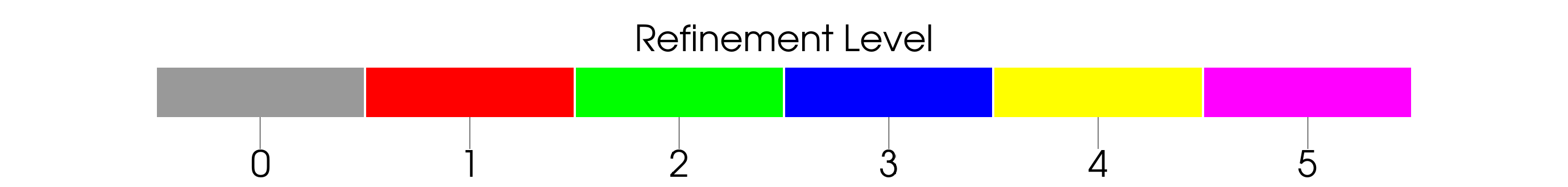}
    \end{subfigure}
    \begin{subfigure}[t]{.48\textwidth}
        \centering
        \includegraphics{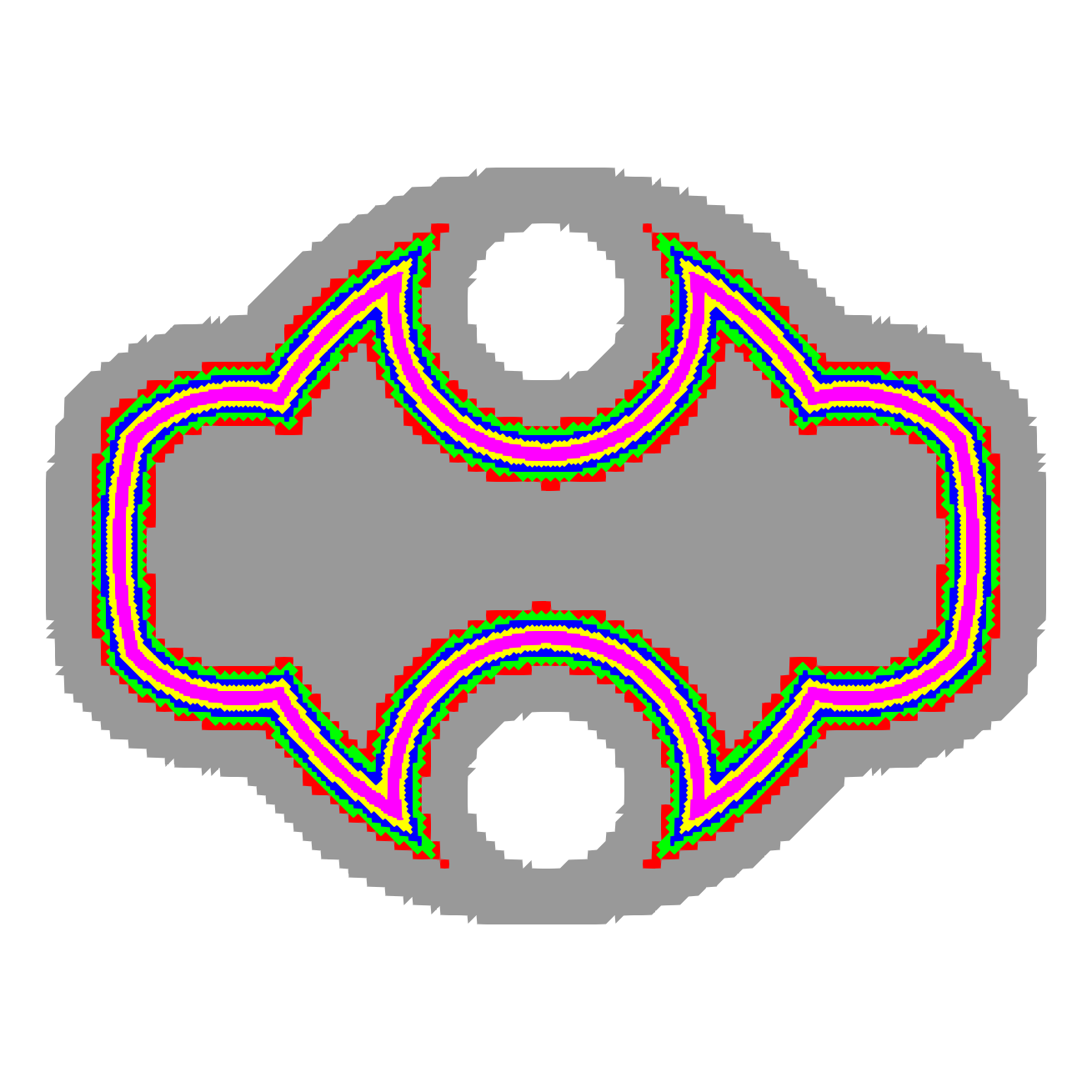}
    \end{subfigure}
    \hfill
    \begin{subfigure}[t]{.48\textwidth}
        \centering
        \includegraphics{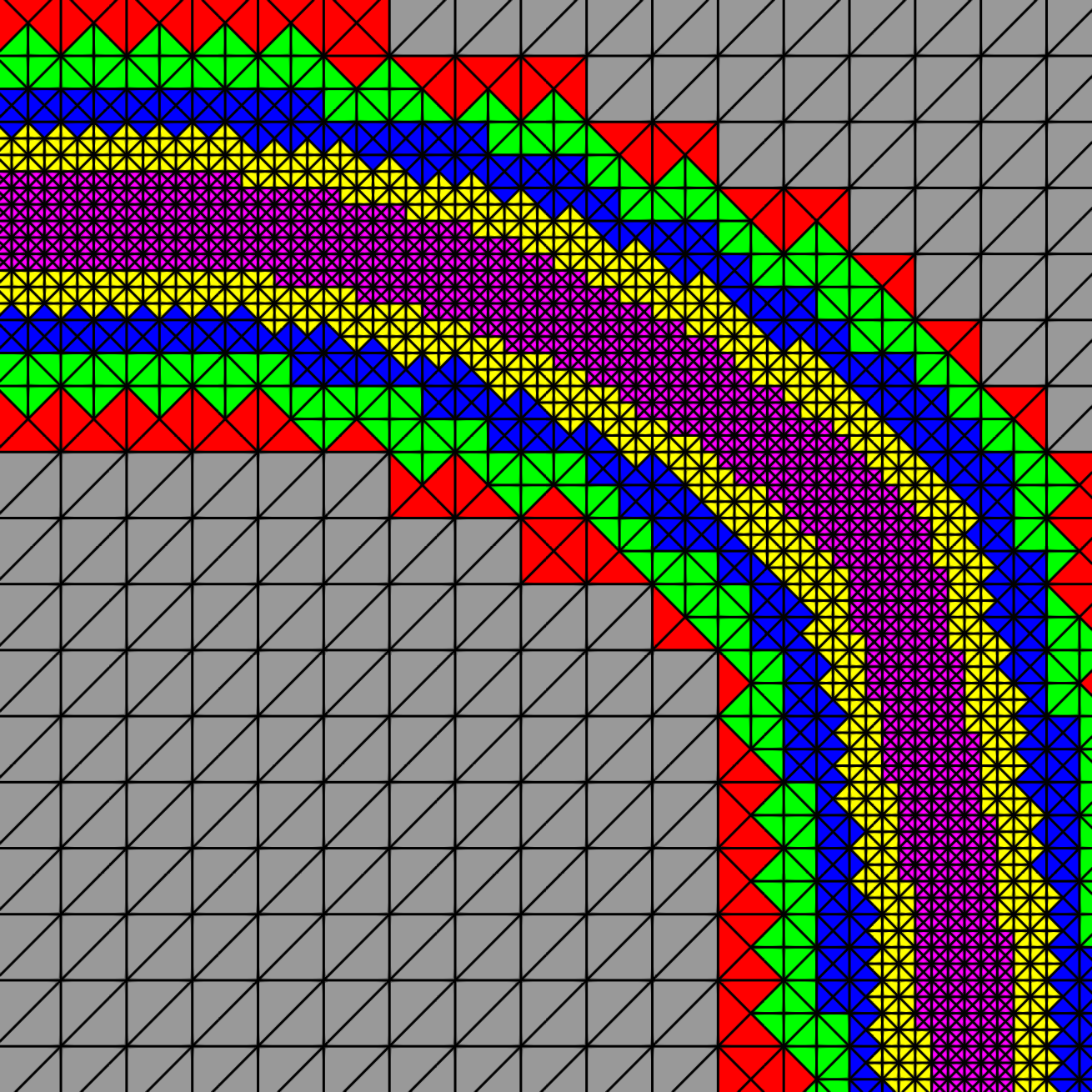}
    \end{subfigure}
    \caption{Adapted and filtered grid view using \dunefem{}.}
    \label{fig:dunemesh}
\end{figure}

The definition of the SDFs and the transformations of models are based solely on \UFL{} as shown and this part of the package can therefore be used with other software based on \UFL{}. For our experiments we use \dunefem{} \citep{Dedner2020}.
We discretise using linear finite elements, and the linear systems are solved using the direct solver UMFPACK \citep{Davis2004UMFPACK} to eliminate additional errors. All fitted grids used for comparison were generated using \gmsh{} \citep{Geuzaine2009}.
A benefit of using \dunefem{} in the \DD{} context is the availability of adaptive and filtered grid views. The filtering allows us to define a simple rectangular mesh, and then remove elements far away from the diffuse boundary, reducing the number of degrees of freedom. Removing mesh elements can be easily done using the SDFs. We decided to remove all mesh elements with $r(x) > 10\epsilon$, following the outer boundary distance recommended in \citet{Yu2020}.
Note that due to the imperfect nature of operating on SDFs, the filtering may include unexpected regions further away than the original boundary, $\domainbnd$.  Using the domain in \cref{lst:sdfsix}, we get an extra loop at the top and bottom of the cut-away regions, which can be seen in \cref{fig:dunemesh}.
When used with \dunefem{}, \ddfem{} also provides a straightforward way to locally adapt the mesh around the boundary $\domainbnd$. We use the surface delta function, which we approximated using $\phi (1-\phi)$ in \cref{eq:unnormweight}, to mark elements for refinement that are in the interfacial region.

For the experiments, we use five levels of refinement across a width of $4\epsilon$.  The levels of refinement for the domain from \cref{lst:sdfsix} are shown in \cref{fig:dunemesh}.  This combination of adaptivity and filtering can produce a mesh with a comparable number of elements to a fitted mesh, although the interface region still requires a fine resolution along straight boundary sections.

\ddfem{} is available through the Python Package Index (PyPI). For details on reproducing the results shown here see \cref{app:code}.

\subsection{Elliptic equation with mixed boundary conditions}
\label{sec:poisson}

\begin{figure}[htbp]
    \centering
    \includegraphics[width=0.475\textwidth]{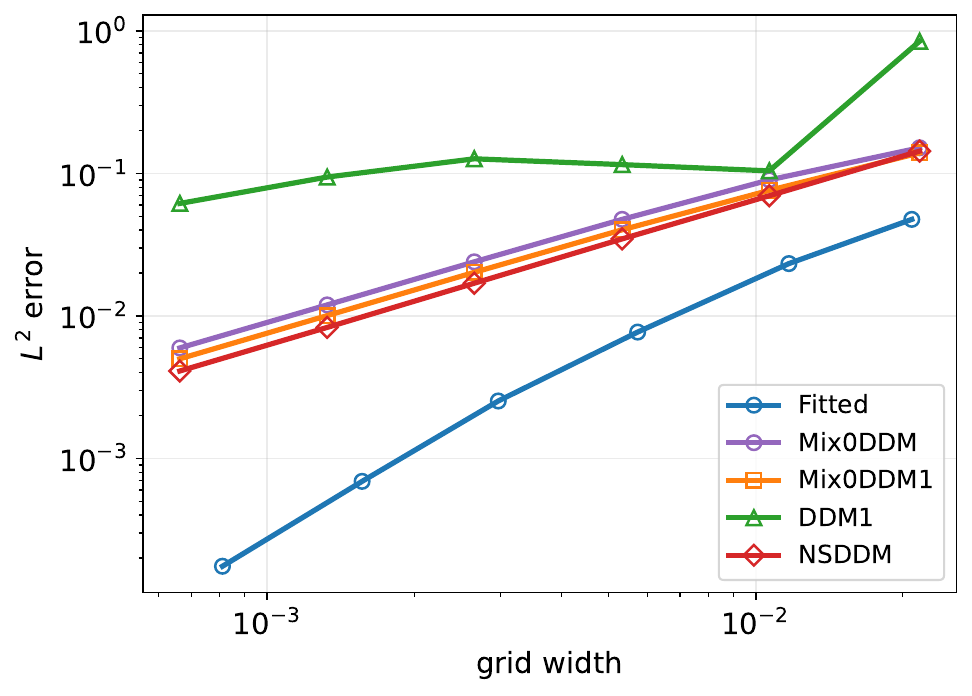}
    \includegraphics[width=0.475\textwidth]{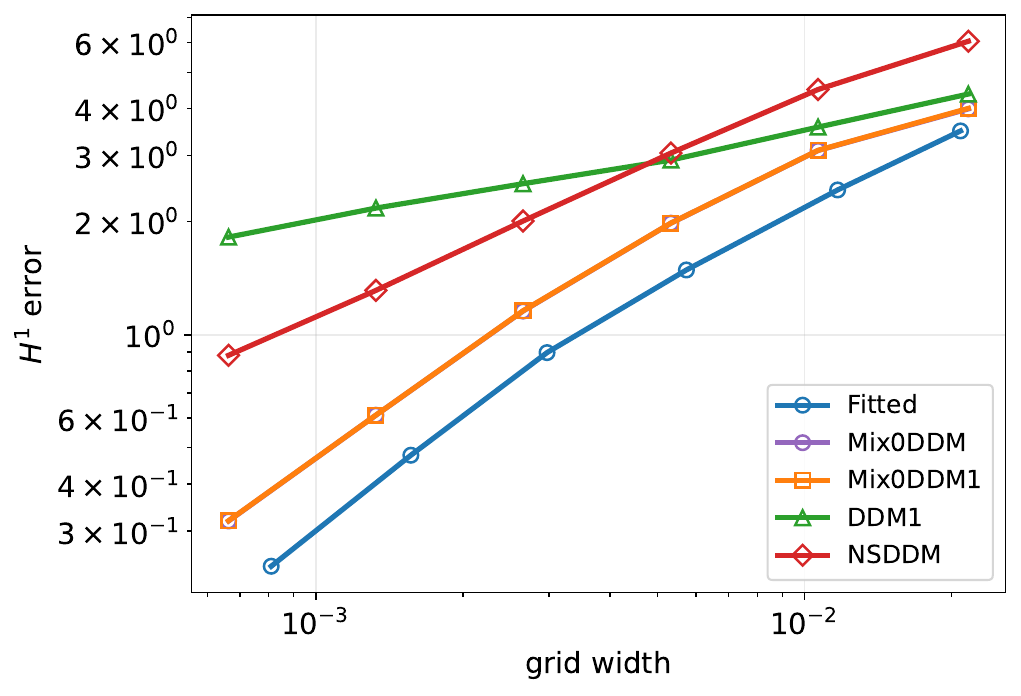}
    \caption{Elliptic problem: comparing errors ($L^2$ on left and $H^1$ on right). Domain is the
        difference of two circles with Dirichlet on outside and Neumann on inside circle.}
    \label{fig:results1}
\end{figure}

In this first experiment we investigate our choice for implementing mixed
boundary conditions. The domain setup is similar to the difference of two
balls shown in \cref{fig:sdfsubtraction}:
$\Omega = B_1((0,0))\setminus B_1((-1,0))$ or using the SDF description
\pyth{omega = Ball(radius=1, center=(0,0)) - Ball(radius=1, center=(-1,0))}.
On this domain we solve the elliptic problem $-\Delta u + 0.01u = f$ with the forcing
and the boundary conditions chosen so that the exact solution is
$u(x,y)=\sin\big(  4\pi (x^2+y^2) \big)$. In \cref{fig:results1} we show $L^2$ errors (left column) and $H^1$ errors (right column) for a number of the methods implemented in our package, and for reference include errors using a fitted grid as well. For these results we use Dirichlet boundary conditions on $\Gamma_{\rm outside}=\partial\Omega\cap \partial B_1(0,0)$ and Neumann conditions on $\Gamma_{\rm inside}=\partial\Omega\cap \partial B_1(-1,0)$. In \cref{app:poisson} we present errors using all Dirichlet conditions \cref{fig:results1DD} and all Neumann in \cref{fig:results1NN}. We already observed in \citep{Benfield2025b} that with Dirichlet boundary conditions \NSDDM{} and \MixedZ{} show very similar $L^2$ norms while \MixedZ{} outperforms the other in the $H^1$ norm. This is due to the natural reconstruction of the gradient of the solution that the mixed method provides, which we do not have for \NSDDM{}. In both norms \DDMO{} is less accurate. With Neumann conditions \MixedZ{} shows some weakness in the $L^2$ error while the other two methods perform identically. In the $H^1$ norm all methods show very similar results at least at smaller values of the grid spacing. These observations (confirmed here) motivates us to consider a new method, which we denote by \MixedZDDM{}. This method blends \MixedZ{} and \DDMO{} so that the former is used for Dirichlet conditions and the latter used for Neumann, thus combining the strengths of either method. To blend the two we use the sum of the weights for all Dirichlet boundaries given by $\omega_D:= \sum_{i \in I_{V}} \omega_i$. In the new model the contributions from \MixedZ{} are weighted with $\omega_D$ and those from \DDMO{} with $1-\omega_D$. The full code defining the blended model is shown in \cref{sec:appMixedZDDM}.
Looking at the results for mixed boundary conditions in \cref{fig:results1} (and \cref{fig:results1ND}) we observe that our approach of handling mixed boundary conditions is effective and that the new blended method \MixedZDDM{} is the best or at least very close to the best \DD{} method we tested in all scenarios considered in this test.

\subsection{Chemical Reaction}
\label{sec:chemR}

We study a reaction-diffusion-advection problem for three species $u_0,u_1,u_2$;
this problem was investigated on a square domain by \citet{Dedner2022}.
Here we will use the domain defined in \cref{lst:sdfthree}.
\begin{listing}[htbp]
    \begin{pythoncode}
sdfs = [
    Ball(radius=0.3, center=(0.15, 0.15)),
    Ball(radius=0.3, center=(-0.15, -0.15)),
    Ball(radius=0.4, center=(0, 0))
]
omega = (sdfs[0] | sdfs[1]) & sdfs[2]
    \end{pythoncode}
    \caption{Defining a three ball SDF.}
    \label{lst:sdfthree}
\end{listing}
We will compare the results across a fitted mesh, which does not use the \DD{} approach,
and an unfitted mesh. Both meshes are locally refined around the boundary of $\Omega$
based on the value of the phase-field function.
However, to avoid any influence of the interior mesh on the amount of chemicals added through the two sources,
we chose to use meshes that match on the original $\domain{}$ for both the fitted and the unfitted simulations.
To this end we used \gmsh{} to produce a conforming mesh with identical interiors.

An interesting aspect of this problem is that we use a discrete
approximation of the velocity field computed also using a \DD{} method.
This problem was investigated more quantitatively in \citep{Benfield2025b} so that we
only show some basic results here to demonstrate the flexibility of our
package for coupling PDEs and solving time-dependent problems.

After discretisation in time, we study the following three
component system for $u=(u_0,u_1,u_2)$
\begin{subequations}
    \label{eq:CRchem}
    \begin{alignat}{2}
        \frac{u^{n} - u^{n-1}}{\tau}
        + \nabla \cdot ( u^{n} \otimes V - D \nabla u^{n})
                                             & = F(t^n,x) - k R(u^{n-1}) & \quad & \text{in } \domain{},    \\
        (D \nabla u^{n} - u^{n} \otimes V )n & = 0                       & \quad & \text{on } \domainbnd{}, \\
        u^{0}                                & = (0, \, 0, \, 0).        &       &
    \end{alignat}
\end{subequations}
We take the diffusion coefficient, $D=10^{-3}$,
a reaction rate, $k=10$, and time step $\tau = 0.05$.

The reaction term is given by two sources around opposite points in the shape
$\Omega$ close to the boundary
\begin{alignat}{2}
    p_0 & = \left(-0.25, \, -0.25\right)~, & \qquad
    p_1 & = \left(0.25, \, 0.25\right) .
\end{alignat}
The resulting source term $F(t,x)$ is active for $t < 10$ in circles around these points,
\begin{align}
    i=0,1 \qquad f_i(x) & =
    \begin{cases}
        5 & \text{if } \|x - p_i\| < 0.04 \\
        0 & \text{otherwise}
    \end{cases}
    \\
    F(t,x)              & =
    \begin{cases}
        \left(f_0(x), \, f_1(x), \, 0\right) & \text{if } t < 10 \\
        \left(0, \, 0, \, 0\right)           & \text{otherwise}
    \end{cases}
\end{align}
The reaction causes $u_0$ and $u_1$ to convert into $u_2$:
\begin{equation}
    R(u)  = \left(u_0 u_1, \, u_0 u_1, \, -2u_0 u_1\right).
\end{equation}
The advective velocity field is given by $V = (\nabla \psi_h)^\perp$,
where we compute a discrete solution $\psi_h$ to
\begin{subequations}
    \label{eq:CRpot}
    \begin{alignat}{2}
        - \Delta  \psi & = - 1 & \quad & \text{in } \domain{},    \\
        \psi           & = 0   &       & \text{on } \domainbnd{}.
    \end{alignat}
\end{subequations}
using the same method as for solving \eqref{eq:CRchem}.

\begin{figure}[htbp]
    \centering
    \begin{subfigure}[t]{0.8\textwidth}
        \centering
        \includegraphics[width=\textwidth]{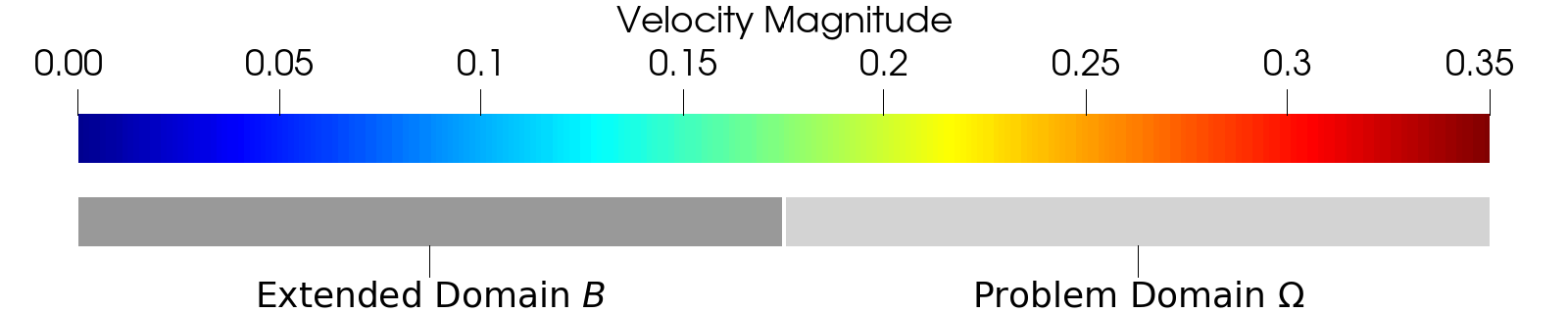}
    \end{subfigure}
    \begin{subfigure}[t]{.32\textwidth}
        \centering
        \includegraphics[width=\textwidth]{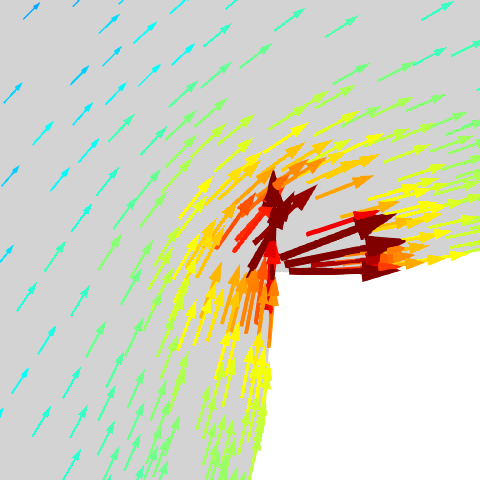}
        \caption{Fitted}
    \end{subfigure}
    \hfill
    \begin{subfigure}[t]{.32\textwidth}
        \centering
        \includegraphics[width=\textwidth]{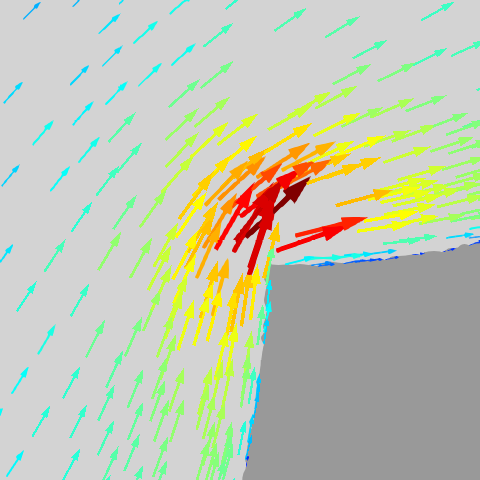}
        \caption{\DDMO{}}
    \end{subfigure}
    \hfill
    \begin{subfigure}[t]{.32\textwidth}
        \centering
        \includegraphics[width=\textwidth]{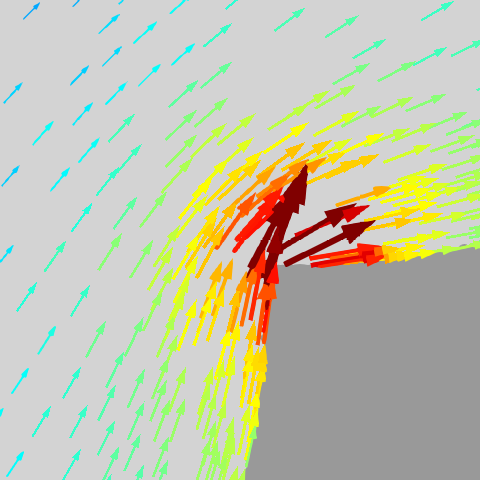}
        \caption{\MixedZ{} and \MixedZDDM{}}
    \end{subfigure}
    \caption{Discrete velocity, zoomed on lower cusp.}
    \label{fig:CRvel}
\end{figure}

Computing this vector field in \cref{fig:CRvel} shows that the velocity fields are closely matched away from the boundary.  However, \DDMO{} shows a small boundary layer with low velocity, and a smaller peak velocity compared to the others. \MixedZ{} has the largest region of velocity at the inner cusps, and smaller magnitude than the fitted mesh.
This is
in line with our investigations in \citep{Benfield2025b} where we designed \MixedZ{} to approximate gradients more accurately.
This demonstrates the importance of choosing the correct \DD{} approximation for a given problem, and the benefit of \ddfem{} collecting multiple transformers for the user.
Note that for the velocity only Dirichlet boundary conditions are present so \MixedZDDM{} is identical to \MixedZ{}.

\begin{figure}[htbp]
    \setlength\tabcolsep{0pt}
    \centering
    \begin{tabular}{@{} >{\centering}m{0.05\linewidth} >{\centering}m{0.22\linewidth} >{\centering}m{0.22\linewidth} >{\centering}m{0.22\linewidth} >{\centering}m{0.22\linewidth} >{\centering\arraybackslash}m{0.05\linewidth} @{}}
        $time$
         & Fitted
         & \DDMO{}
         & \MixedZ{}
         & \MixedZDDM{}
         &                                                      \\
        $16$
         & \includegraphics[width=0.22\textwidth]{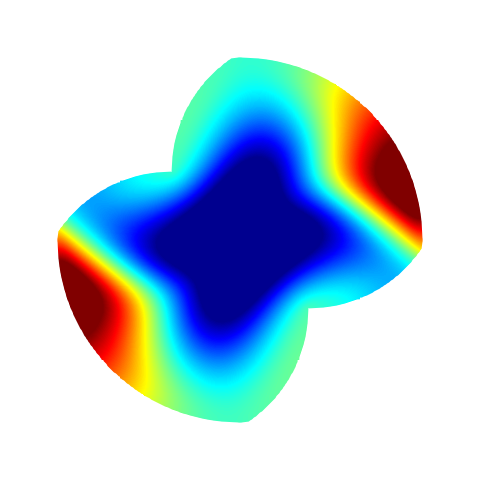}
         & \includegraphics[width=0.22\textwidth]{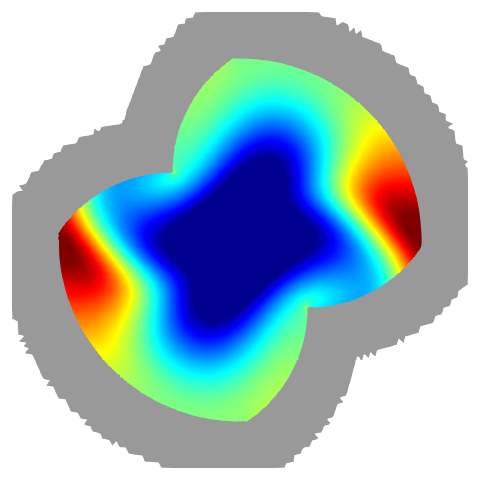}
         & \includegraphics[width=0.22\textwidth]{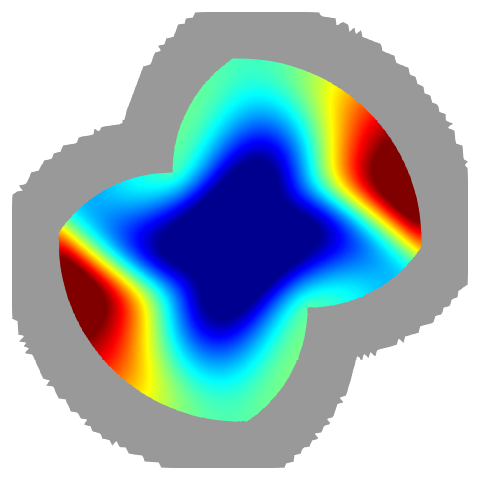}
         & \includegraphics[width=0.22\textwidth]{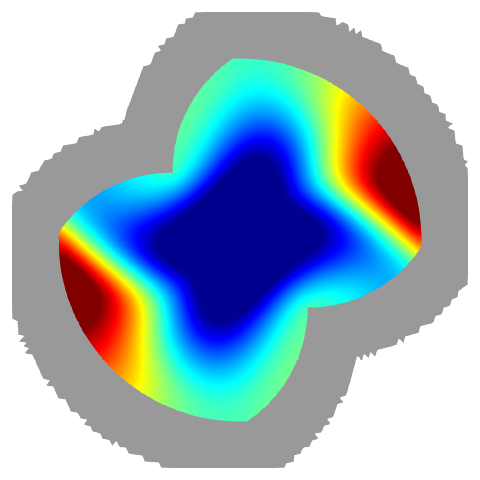}
         &                                                      \\
        $26$
         & \includegraphics[width=0.22\textwidth]{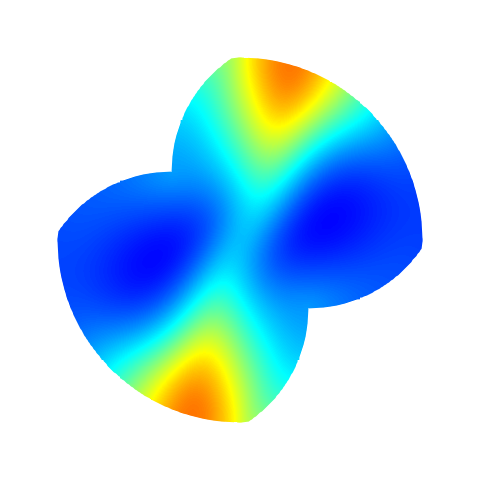}
         & \includegraphics[width=0.22\textwidth]{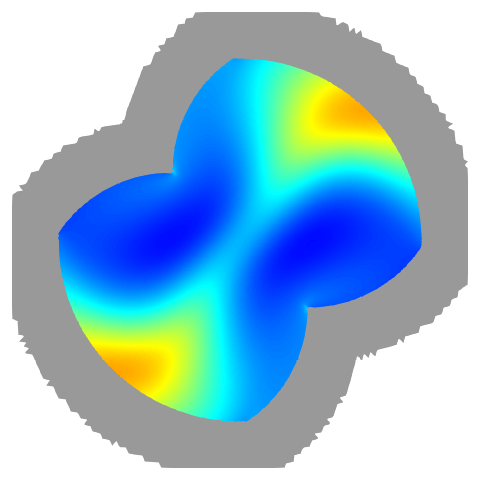}
         & \includegraphics[width=0.22\textwidth]{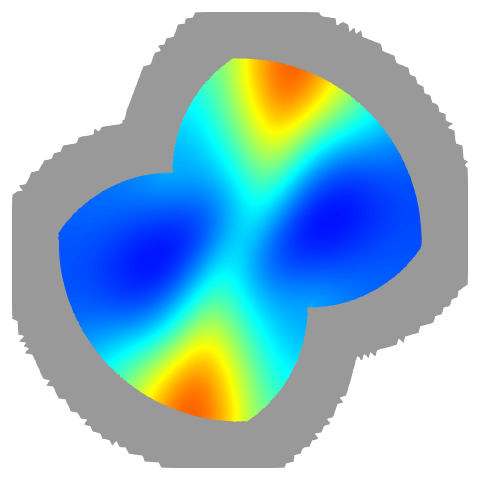}
         & \includegraphics[width=0.22\textwidth]{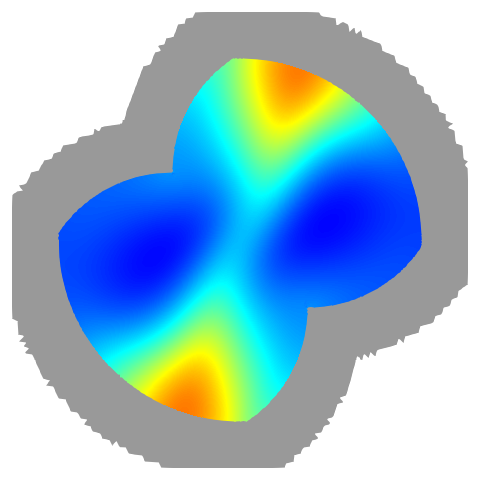}
         &                                                      \\
        \multicolumn{5}{c}{\includegraphics[width=0.8\textwidth]{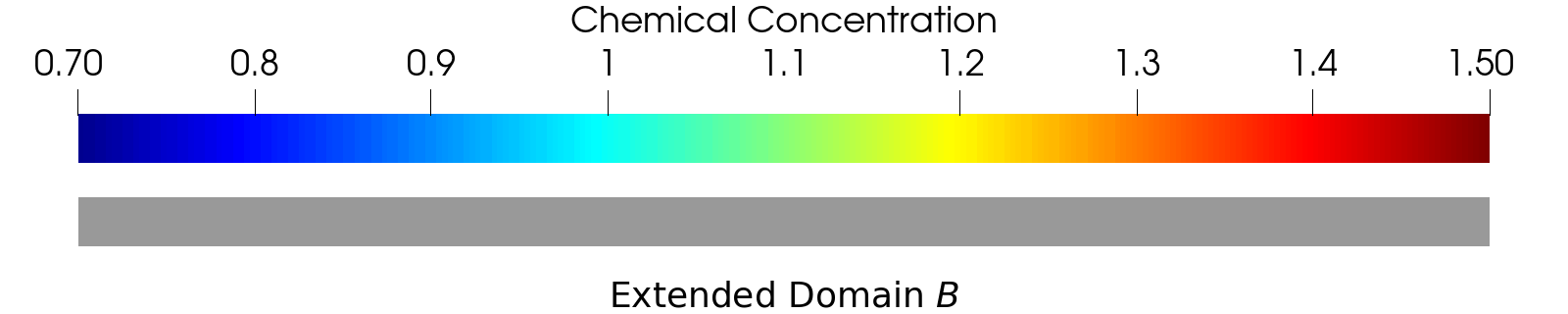}}
    \end{tabular}
    \caption{Approximate solution (concentration magnitude $|u|$) to the chemical reaction problem after the source becomes inactive.}
    \label{fig:CRend}
\end{figure}

\begin{figure}[htbp]
    \centering
    \includegraphics[width=0.9\textwidth]{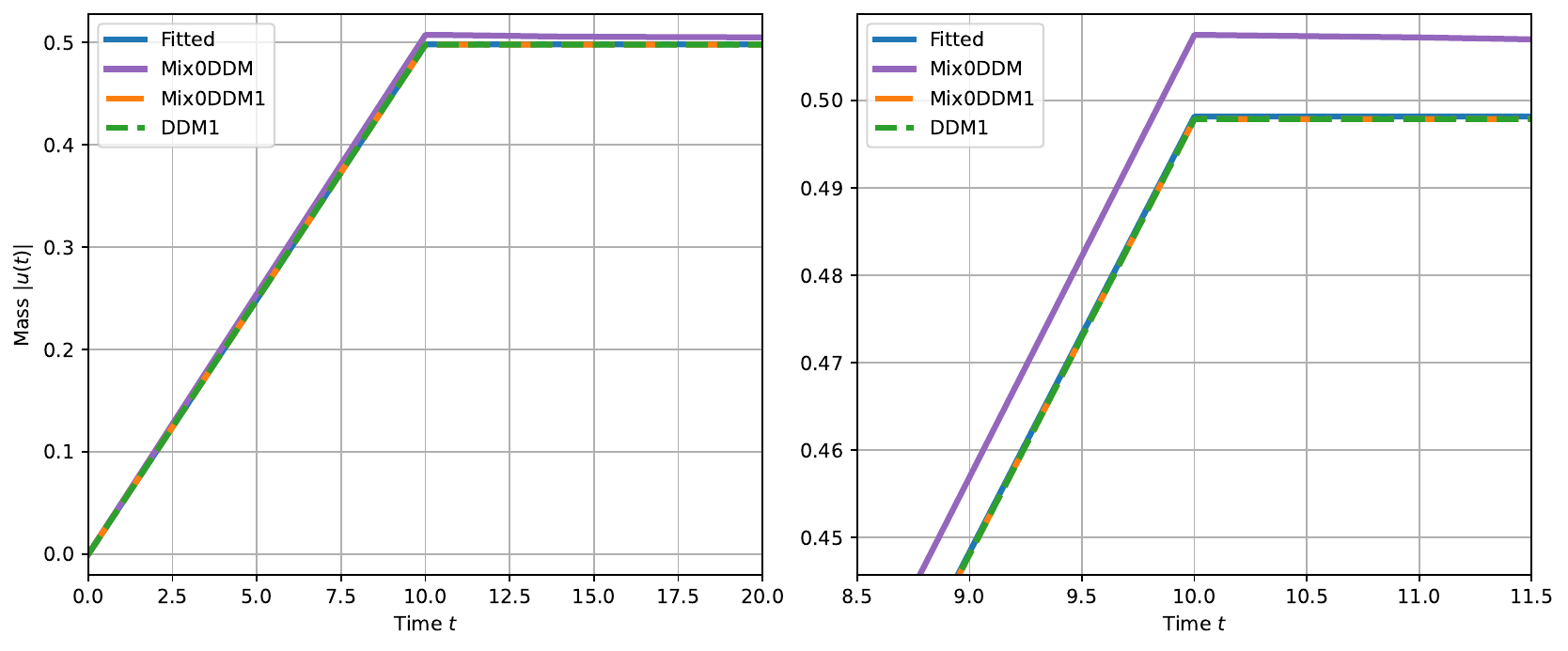}
    \caption{Comparison of computed mass for chemical reaction problem for \DD{} methods and fitted mesh.}
    \label{fig:CRMass}
\end{figure}

We only show results of the reaction process at later times in \cref{fig:CRend}, showing the concentration magnitude $|u(x)| = \sqrt{u_0^2(x) + u_1^2(x) + u_2^2(x)}$, where it can be observed that the transportation using \DDMO{} trails behind the fitted solution, \MixedZ{}, and \MixedZDDM{}, which is a consequence of the lower velocity and boundary layer observed in \cref{fig:CRvel}. In addition, the concentration is lower in \DDMO{}.
For a more quantitative comparison of the methods, we also show the evolution of the overall mass $\int_\Omega u_0+u_1+u_2$, computed on the unfitted domain by using a high-order quadrature to approximate $\int_{\extdomain{}}{ \chi_\domain{} \left(u_0+u_1+u_2\right) }$.
On the left of \cref{fig:CRMass} we show the evolution over the full time interval and on the right around $t=10$ where no extra chemicals are added to the system, so that the mass should be conserved. We can see that, with the exception of \MixedO{}, the unfitted methods conserve mass quite well, if at a slightly lower level compared to the fitted simulations. \MixedZDDM{} is almost indistinguishable from \DDMO{}, while \MixedZ{} not only ends up with too much mass at $t=10$ but then noticeably loses mass after that point in time.
.


\subsection{Cahn-Hilliard equation}
\label{sec:cahnh}

The Cahn-Hilliard equation has been explored with the \DD{} method before by \citet{Aland2010}.
We will solve the Cahn-Hilliard equation using a mixed formulation,
splitting the fourth-order PDE into a system of second-order equations
for the concentration field $c$ and the chemical potential $\mu$.
With time step $\tau = 5 \times 10^{-6}$ until 0.001
and $\tau = 4 \times 10^{-5}$ until 0.008,
we get,
\begin{subequations}
    \begin{alignat}{2}
        \frac{c^{n} - c^{n-1}}{\tau} - 0.5 \Delta(\mu^{n} + \mu^{n-1}) & = 0 & \quad & \text{in } \domain{}, \\
        \mu^{n} -  \frac{d f(c^{n})}{d c} + \lambda \Delta c^{n}       & = 0 & \quad & \text{in } \domain{}.
    \end{alignat}
\end{subequations}
This uses the double-well potential,
\begin{equation}
    f(c) = 100 c^2 (1-c)^2,
\end{equation}
with
\begin{equation}
    \frac{d f}{d c} = 200 c (2c^2 - 3c + 1),
\end{equation}
and gradient energy $\lambda=0.01$.

The \ddfem{} package can transform systems of equations.
So, we have the boundary conditions,
\begin{subequations}
    \begin{alignat}{2}
        \nabla \mu^{n} \cdot n       & = 0 & \quad & \text{on } \domainbnd{},
        \\
        \lambda \nabla c^{n} \cdot n & = 0 & \quad & \text{on } \domainbnd{}.
    \end{alignat}
\end{subequations}

\begin{figure}[htbp]
    \setlength\tabcolsep{0pt}
    \centering
    \begin{tabular}{@{} >{\centering}m{0.09\linewidth} >{\centering}m{0.3\linewidth} >{\centering}m{0.3\linewidth} >{\centering}m{0.3\linewidth} >{\centering\arraybackslash}m{0.05\linewidth} @{}}
        $time$
         & Fitted
         & \DDMO{}, \MixedZDDM{}
         & \MixedZ{}
         &                                                      \\
        $1.6e-4$
         & \includegraphics[width=0.25\textwidth]{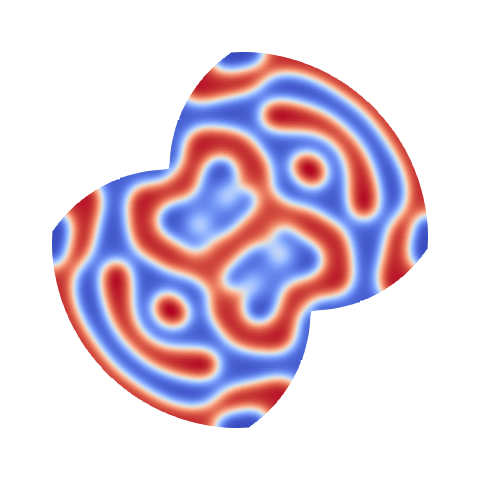}
         & \includegraphics[width=0.25\textwidth]{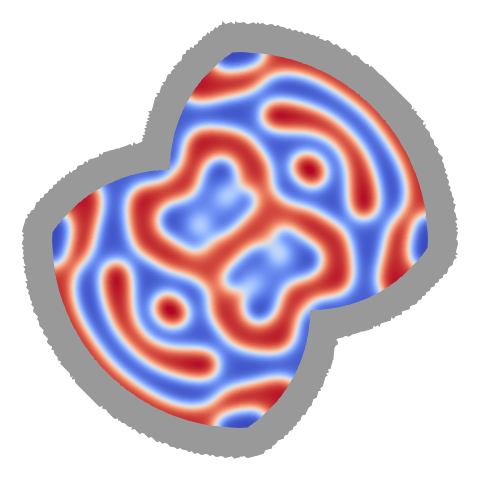}
         & \includegraphics[width=0.25\textwidth]{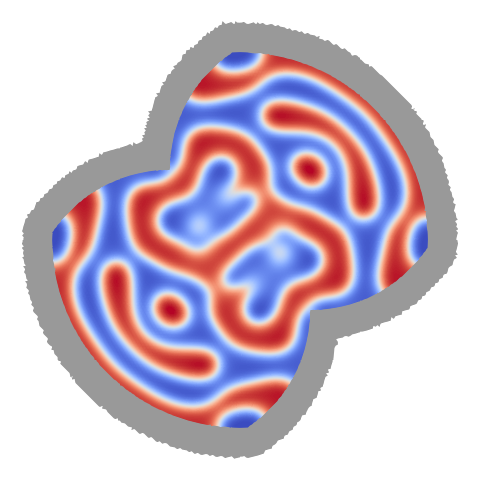}
         &                                                      \\
        $8e-4$
         & \includegraphics[width=0.25\textwidth]{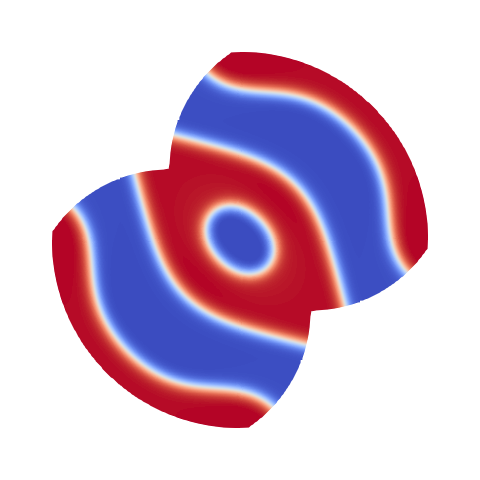}
         & \includegraphics[width=0.25\textwidth]{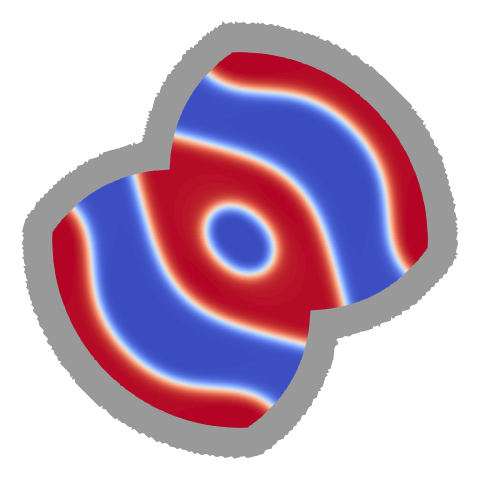}
         & \includegraphics[width=0.25\textwidth]{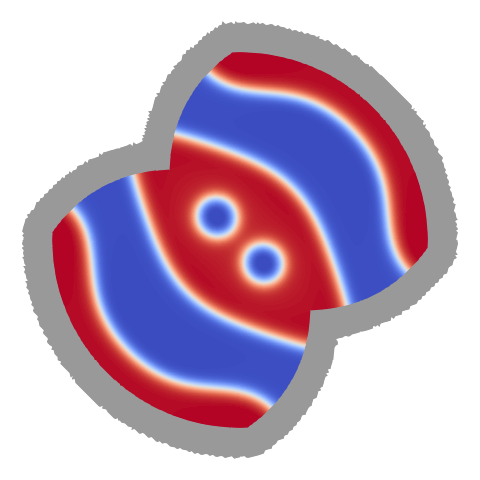}
         &                                                      \\
        $8e-3$
         & \includegraphics[width=0.25\textwidth]{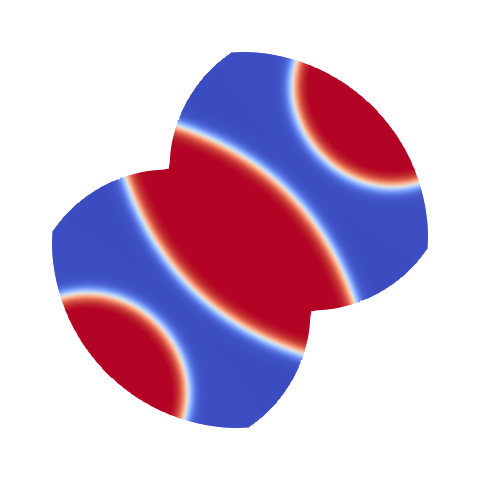}
         & \includegraphics[width=0.25\textwidth]{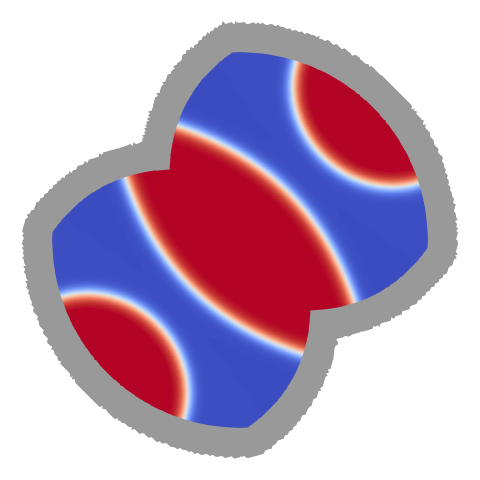}
         & \includegraphics[width=0.25\textwidth]{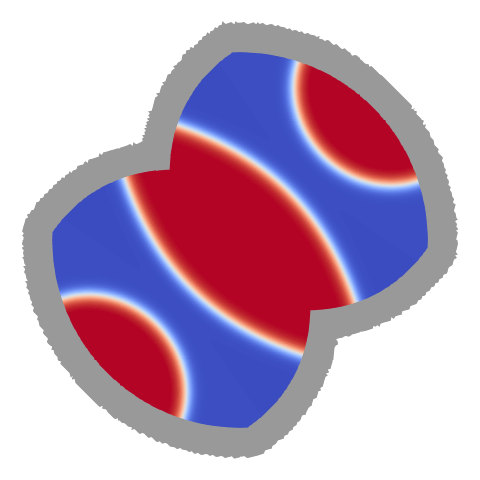}
         &                                                      \\
        \multicolumn{5}{c}{\includegraphics[width=0.75\textwidth]{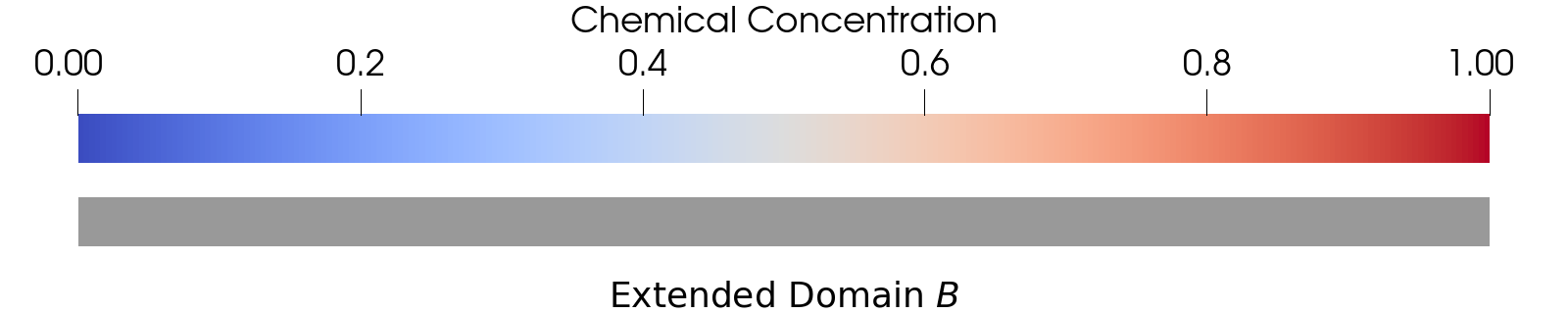}}
    \end{tabular}
    \caption{Comparison of approximate solutions to the Cahn-Hilliard equation using identical interior meshes.}
    \label{fig:CH}
\end{figure}

We initialise the concentration with a small oscillating region in the centre of the domain:
\begin{subequations}
    \begin{align}
        c^{0}(x) & =
        \begin{cases}
            0.5 + 0.1 \sin(20\pi x_1) \sin(20\pi x_2) & \text{if } |x_1|, |x_2| < 0.1 \\
            0.5                                       & \text{otherwise}
        \end{cases}
        \\
        \mu^{0}  & = \frac{d f(c^{0})}{d c} - \lambda \Delta c^{0}.
    \end{align}
\end{subequations}
We take the domain from the previous example with the SDF given in \cref{lst:sdfthree} and again construct a direct extension of the fitted mesh to avoid issues with the projection of the initial data differing between fitted and \DD{}, guaranteeing identical initial mass. The grids are refined once compared to that used in the previous example.
This gives the results in \cref{fig:CH}.
We observe no significant difference between the fitted and \DDMO{}. Although the problem only uses flux boundary conditions, \MixedZ{} only shows very minor deviations. Refining the grid one level does bring this method in line with the others. Note that \MixedZDDM{} is identical to \DDMO{} for this problem since no Dirichlet boundaries are present and \NSDDM{} produced results in line with \DDMO{}.

\begin{figure}[htbp]
    \centering
    \includegraphics[width=0.9\textwidth]{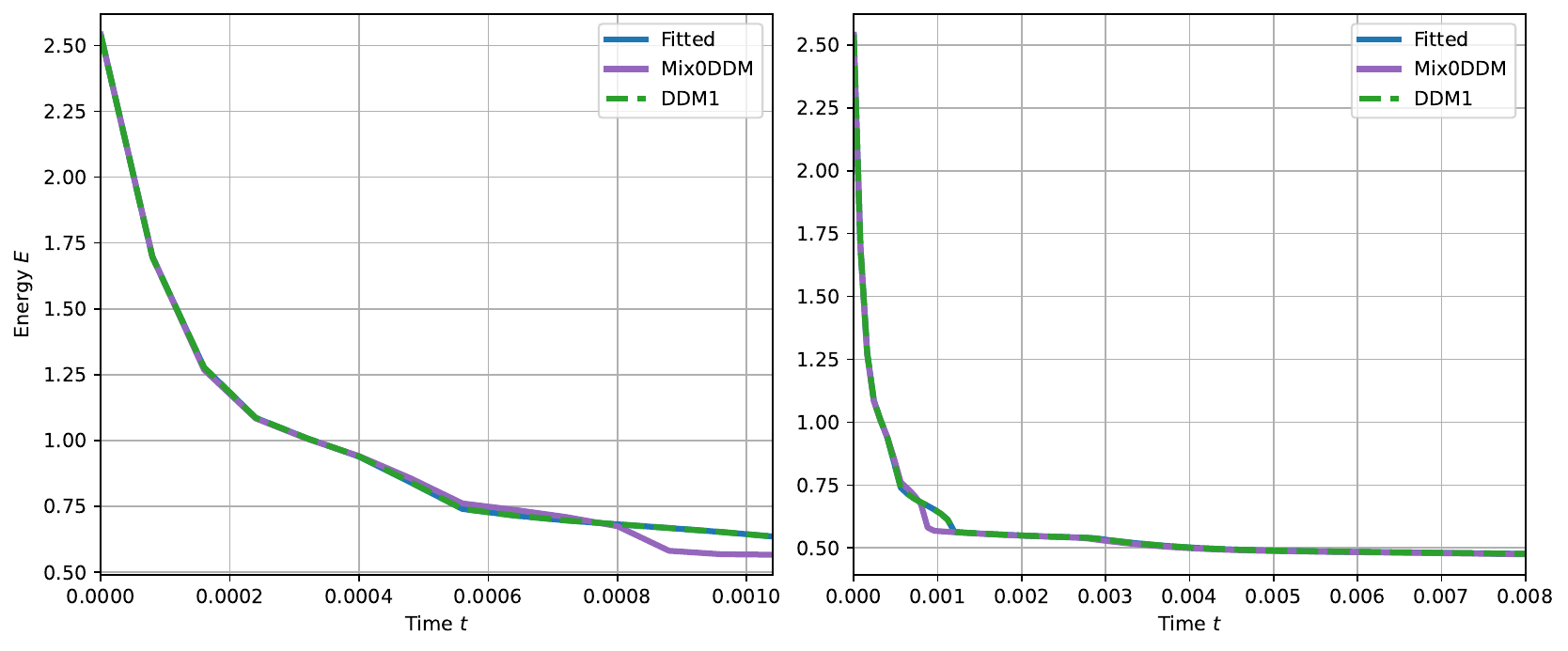}
    \caption{Comparison of computed energy for Cahn-Hilliard problem for \DD{} methods and fitted mesh.}
    \label{fig:chEnergy}
\end{figure}

For further insight we can also compare the energy defined by
\begin{equation}
    E(c)
    = \int_{\domain{}}{ \frac{\lambda}{2} |\nabla c|^2  + f(c)\; dx}
    = \int_{\extdomain{}}{ \chi_\domain{} \left( \frac{\lambda}{2} |\nabla c|^2  + f(c) \right) \; dx}.
\end{equation}
For the \DD{} methods, energy is computed over $\extdomain{}$ using the characteristic function $\chi_\domain{}$ with a suitably high-order quadrature.
Results are shown in \cref{fig:chEnergy}.
All \DD{} methods produce close to identical results, remaining close to the fitted energy curve, although \MixedZ{} does show a different behaviour towards the end of the simulation, which again vanishes after grid refinement. We also compare the mass in \cref{app:ch} which again demonstrates that \MixedZ{} does not conserve mass as already shown in the previous test case, again confirming the accuracy issue when using Neumann boundary conditions.

\subsection{Linear Elasticity}
\label{sec:linelast}

With this example, we continue testing our treatment of mixed boundaries.
The plane-stress model of the small elasticity deformation of a cantilever beam under
gravity can be described by
\begin{subequations}
    \begin{alignat}{2}
        - \nabla \cdot \sigma(u) & = \left( 0, \, -9.81\rho\right) & \quad & \text{in } \domain{},                              \\
        u                        & = 0                             & \quad & \text{on } \domainbnd{}_D,                         \\
        \sigma \cdot n           & = 0                             & \quad & \text{on } \domainbnd{} \setminus  \domainbnd{}_D.
    \end{alignat}
\end{subequations}
For linear elasticity \citep{Aland2012} the stress tensor is given by
\begin{equation}
    \sigma(u) = \mu (\nabla u + \nabla u^T) + \lambda (\nabla \cdot u) I~.
\end{equation}
We use $\lambda=E \nu / (1 - \nu^2)$ and $\mu=E/2 (1 + \nu)$ where we choose representative material constants $E=2\cdot 10^{11},\nu=0.3$, and $\rho=7800$ for generic structural steel \citep{en1993}. All quantities are expressed in SI units, i.e., seconds, kilograms, and metres.
The beam is defined by the domain described by the SDF from \cref{lst:sdfsix}. The beam is fixed on the left boundary $\domainbnd{}_D$ but is free to move under its own weight. Consequently, the deformation given by $u$ satisfies Dirichlet boundary conditions $u=0$ along the boundary of the intersection of the balls with the names \pyth{left} and \pyth{bounding}; traction free boundary conditions are used on the remaining boundary. The full model description used is shown in \cref{app:LEM}.

\begin{figure}[htbp]
    \centering
    \begin{subfigure}[c]{0.45\textwidth}
        \hspace*{\fill}
        \includegraphics[align=c,width=\textwidth]{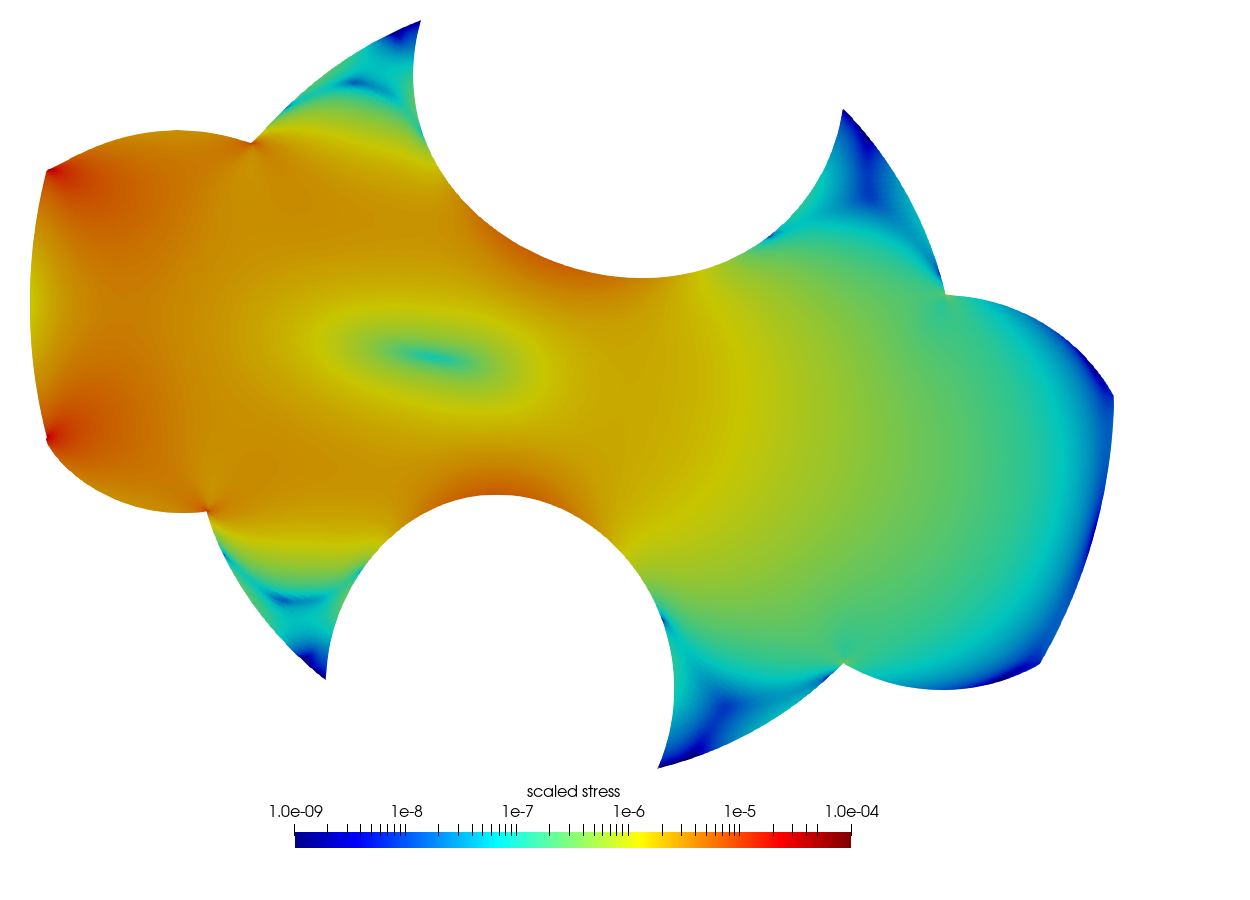}
        \caption{Normalised von Mises stress of the fine grid reference solution on scaled deformed domain. Colour bar is logarithmic.}
        \label{fig:LErefStress}
    \end{subfigure}
    \begin{subfigure}[c]{0.49\textwidth}
        \includegraphics[align=c,width=\textwidth]{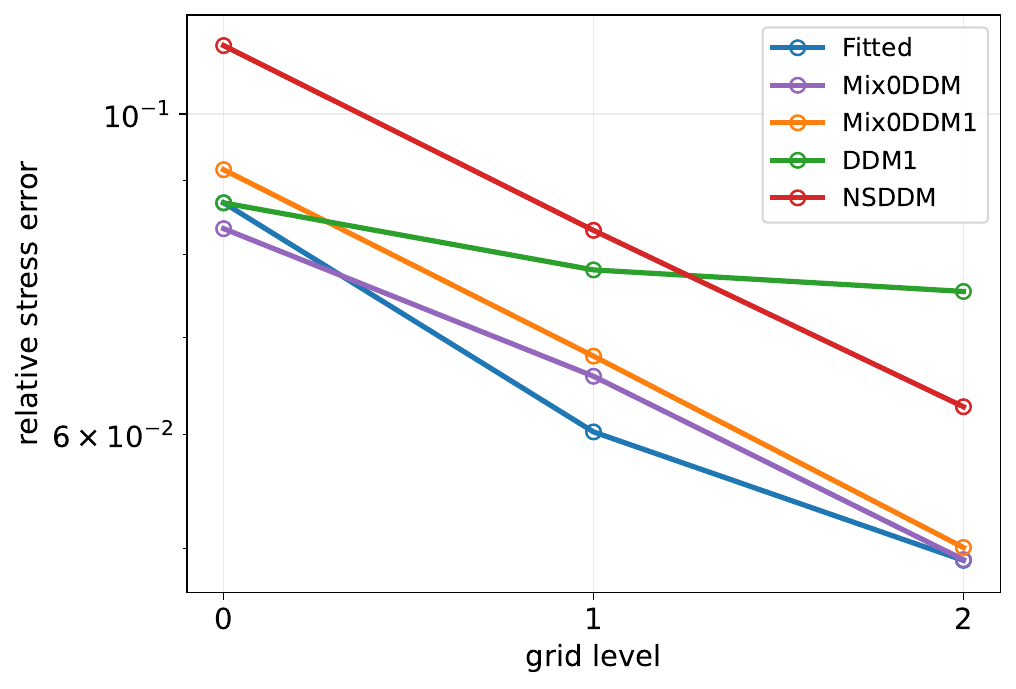}
        \caption{Relative $L^2$ stress errors between the fine grid reference solution and approximation computed on three coarser grids.}
        \label{fig:LEL2ErrorStress}
    \end{subfigure}
    \caption{Linear elasticity example: comparison of the von Mises stress
        normalised by $\lambda+2\mu$.
        The deformed domain is shown with the displacement scaled by $10 000$.}
    \label{fig:LE}
\end{figure}

Although a thorough quantitative comparison of the \DD{} methods is outside the scope of this paper, we include a few results for this problem due to the mixed boundary conditions used. To this end we have computed a reference solution on a fitted grid refined towards the boundary of the domain with slightly more than $136\ 000$ elements and interpolated the coarser grid \DD{} and fitted solutions into a piecewise linear Lagrange space over the reference grid. Errors are then computed on this grid. We show results for three grid levels with $27\;938,59\;766,131\;220$ elements, respectively. In addition, we compute on three coarser fitted grids with $14\;738,34\;065,83\;621$ elements.
In \cref{fig:LE} we study the resulting von Mises stress given by the norm of
\begin{equation}
    s = \sqrt{\sigma_{xx}^2-\sigma_{xx}\sigma_{yy}+\sigma_{yy}^2+3\sigma_{xy}^2}~.
\end{equation}
The reference von Mises stress is shown in \cref{fig:LErefStress} and the $L^2$ error between the fine fitted grid reference solution and the coarser grid solutions is shown in \cref{fig:LEL2ErrorStress}. Results are quite similar to those obtained in our first test case with \DDMO{} showing the highest error and either \MixedZ{} or \MixedZDDM{} closest to the fitted results.

\begin{figure}
    \centering
    \includegraphics[width=0.8\textwidth]{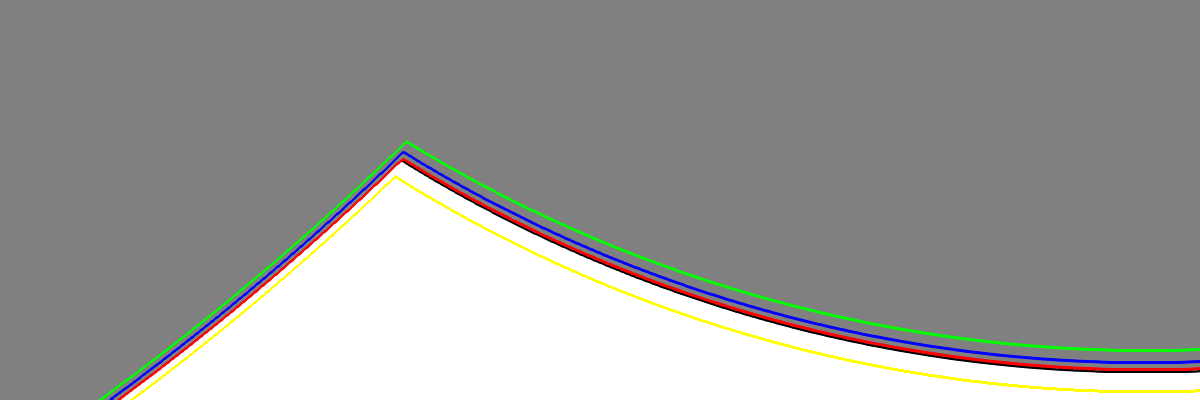}
    \caption{Deformed domain (displacement scaled by $10 000$) around the lower right reentrant corner. Reference domain given by black curve and blue showing the boundary from the coarser fitted simulation. Yellow, green, red are computed with \MixedZ{}, \DDMO{}, \MixedZDDM{}, respectively. All from the level one grids.}
    \label{fig:LEdeformed}
\end{figure}

In \cref{fig:LEdeformed} we plot the boundary of the lower right of the deformed domain. We can see that all methods produce results close to the reference domain with the \MixedZDDM{} being closer to the reference then the fitted coarser grid results. The \NSDDM{} method is so similar to fitted that we did not include it here, \MixedZ{} and \DDMO{} produce similar deformed domains.

In \cref{app:linelast} we provide results giving more insight into the performance of the methods along the boundary of the reference grid.

\section{Conclusion}
\label{section:conclusion}

In this paper we have presented the \ddfem{} Python module,
which allows straightforward application of the \DD{} method for solving a wide
range of PDEs on complex domains.
By formulating constructive solid geometry and boundary approximations directly in \UFL{},
the package remains independent of mesh generation and the underlying solver.
A key feature of the framework is its modularity, and therefore easily extendable.
Different \DD{} formulations can be interchanged via transformer classes,
defined using decorators to compose the components of processing the user's model.
The weighting approach allows mixed boundary condition types to be
combined without requiring additional phase-field variables.
Numerical experiments across
Poisson,
advection-diffusion-reaction,
Cahn-Hilliard,
and elasticity problems
to demonstrate that the package effectively handles complex domains.
These experiments reproduced the behaviour of traditional fitted meshes,
including useful boundary flux and stress approximations,
whilst avoiding the computational overhead of fitted mesh generation.

\begin{figure}[htbp]
    \setlength\tabcolsep{0pt}
    \centering
    \begin{tabular}{@{} >{\centering}m{0.05\linewidth} >{\centering}m{0.45\linewidth} >{\centering}m{0.45\linewidth} >{\centering\arraybackslash}m{0.05\linewidth} @{}}
        $\theta$
         & \DDMO{}
         & \NSDDM{}
         &                                \\
        $\frac{\pi}{2}$
         & \includegraphics{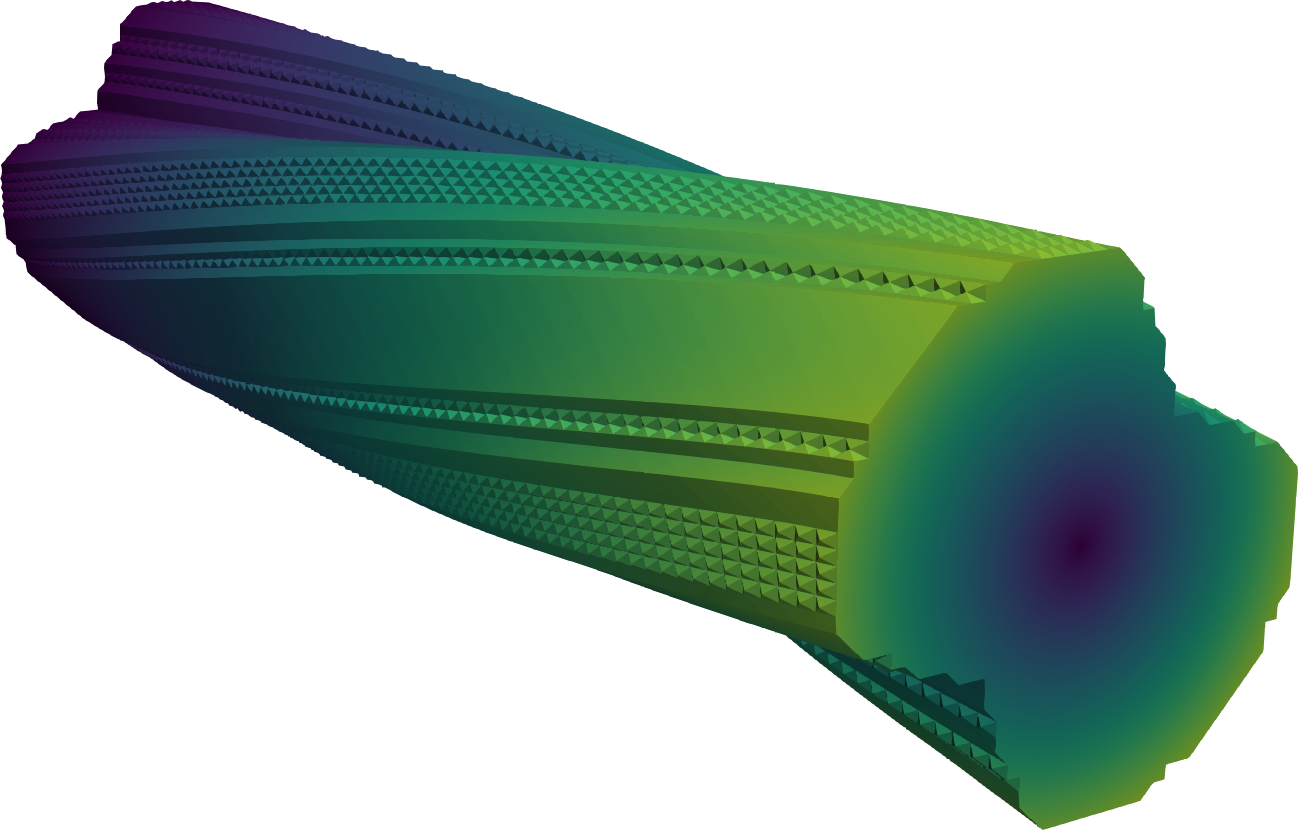}
         & \includegraphics{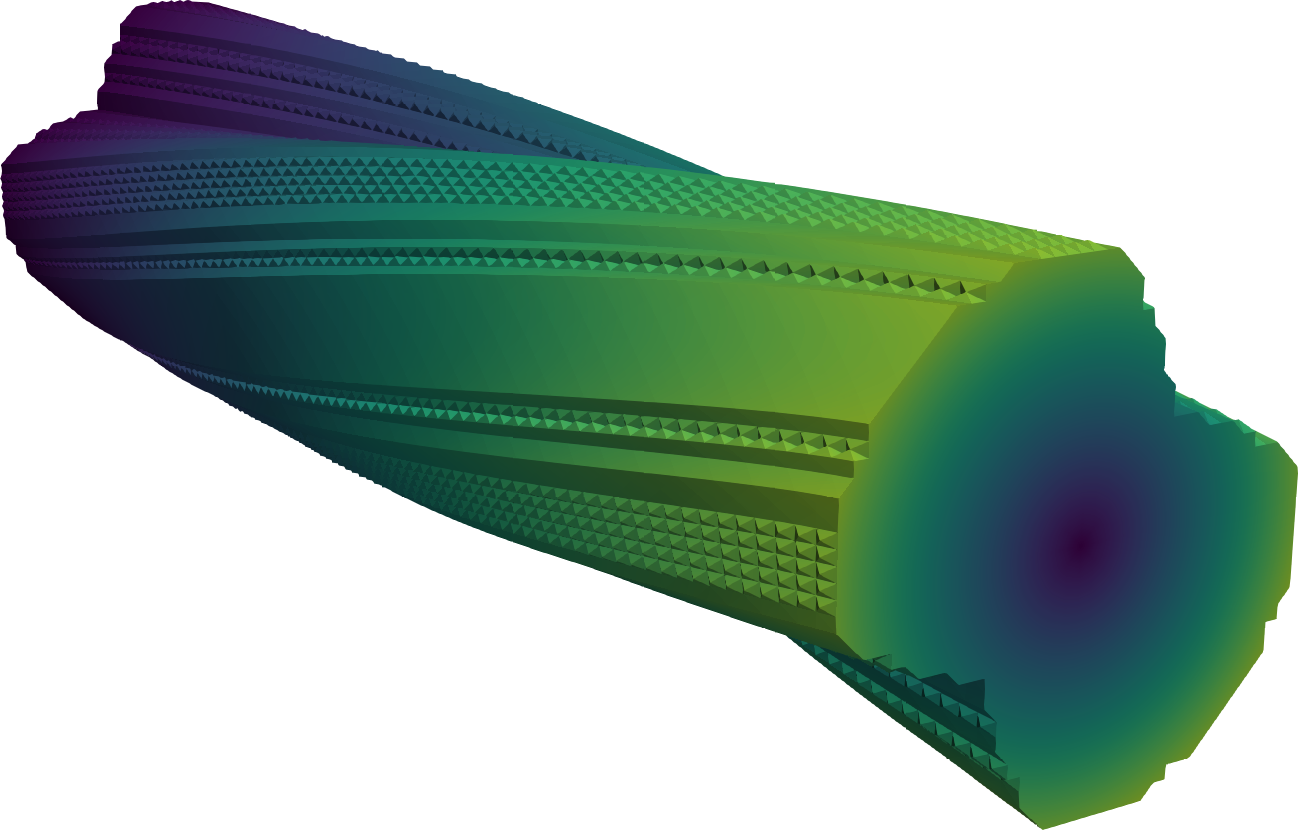}
         &                                \\
        $\pi$
         & \includegraphics{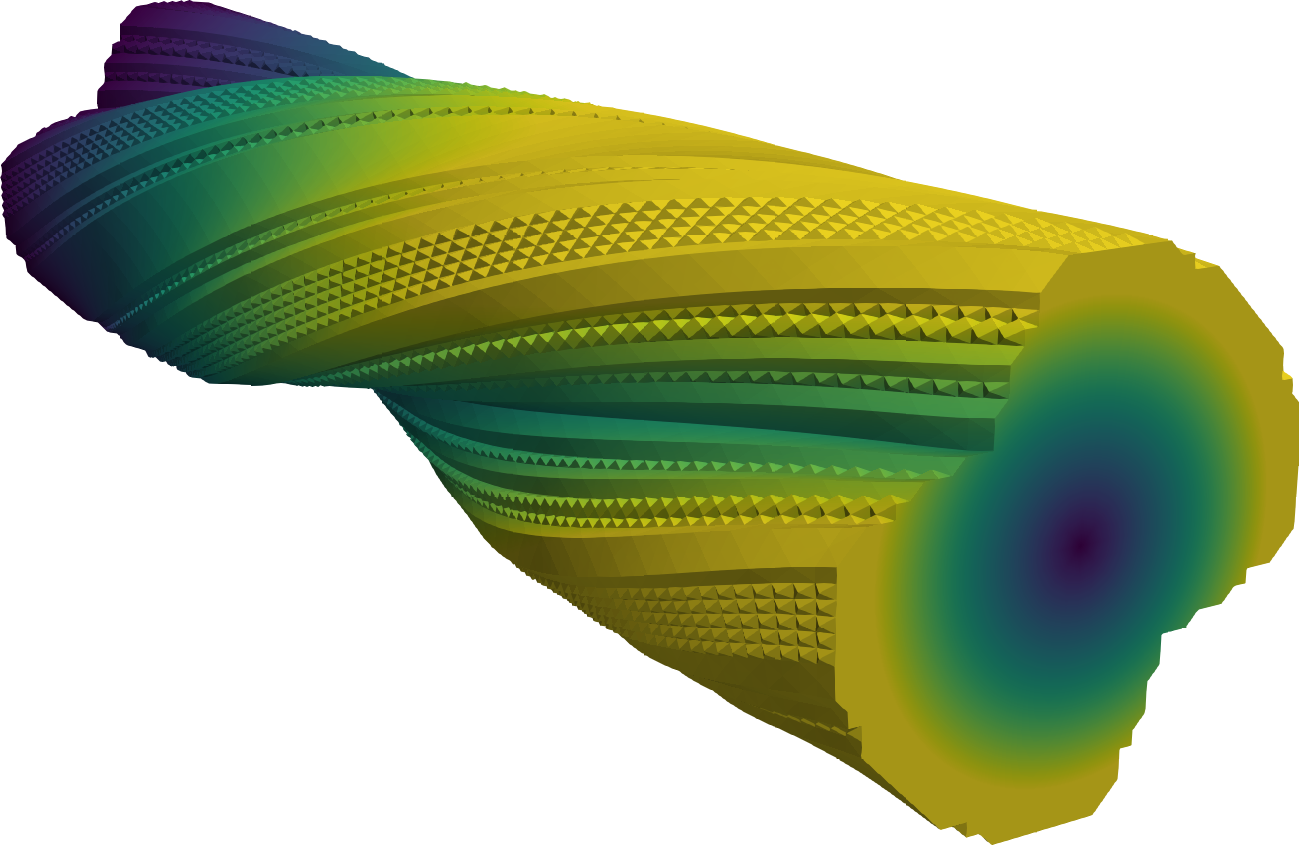}
         & \includegraphics{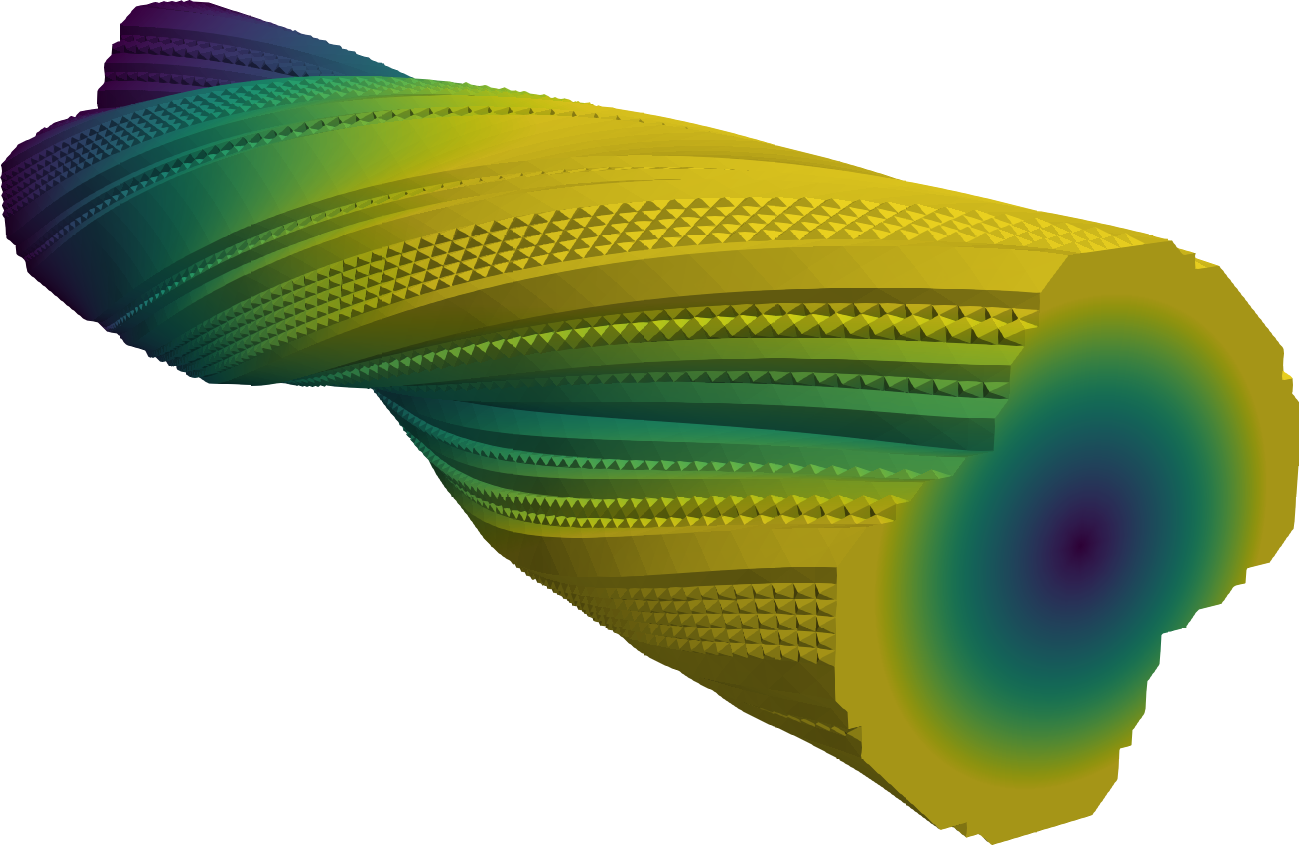}
         &                                \\
        \multicolumn{4}{c}{\includegraphics{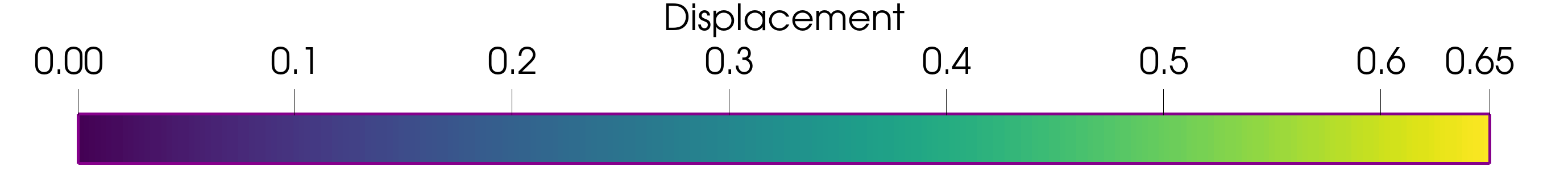}}
    \end{tabular}
    \caption{
        The problem domain has been warped to show the approximate
        solutions for displacement to the 3D hyperelasticity problem.
    }
    \label{fig:HE}
\end{figure}

All the above examples were computed on 2D domains. The implementation is set up to work in 3D and we have included a number of 3D SDF functions and two specific modifiers for 2D SDFs, extrusion \pyth{ddfem.geometry.Extrusion} and revolution around an axis \pyth{ddfem.geometry.Revolution}.
So far we have not thoroughly tested the package on 3D examples. Also, testing to solve non-linear problems has only been started. As an example, we have studied a \emph{3D Hyperelasticity} model from \citep{bleyer2024}, simulating a 3D beam twisting with a compressible neo-Hookean model.
We start with the domain from \cref{lst:sdfthree}, which we extrude using the Python code below to create a domain of a complex beam.
\begin{pythoncode}
omega3D = omega.extrude(3.2)
\end{pythoncode}
With this, the domain will only be extended around the complex boundaries, so the flat faces on the end of the beam will match the extended 2D domain. On one of these flat boundaries the beam is fixed, on the other it is rotated by an angle $\theta$.
On the remaining boundaries of the beam, zero traction is prescribed, which allows the beam to twist. The problem is solved incrementally by stepping the rotation angle, $\theta$, across $60$ values from $0$ to $\pi$.
An example of the results is shown in \cref{fig:HE}, both \DD{} methods produce similar results.
Note that so far we have only studied the \DD{} method qualitatively.

At present, the package is limited in handling vector-valued systems
where each component is restricted to
the same type of boundary condition when defined within a single model class.
We are actively testing and validating the automated selection of
penalty scaling factors, introduced in \cref{section:transformers},
to further reduce the need for problem-specific tuning.
We are also working on extending the \pyth{ddfem.geometry} subpackage to
generate SDFs from a closed loop of line segments, arcs, and Bézier segments,
and using a marching cubes algorithm to extract SDFs from image data.
The current implementation can also handle discrete SDFs, including the SDF operators,
data extension, and boundary weight computation.
To test this we have added a class \pyth{ddfem.GridSDF} that can handle a SDF computed from
polygonal boundaries and triangulated surfaces. A more detailed investigation of these possibilities,
including facilitating the solving of surface-bulk coupling problems are ongoing.
Finally, we plan to carry out further quantitative investigations of 3D and non-linear problems,
as well as problems involving incompressible and compressible fluid flow,
and explore time-dependent problems on moving domains.

\bibliography{references}

\cleardoublepage
\appendix

\section{Reproducing results}
\label{app:code}
\ddfem{} is available on PyPI and can be installed using \pyth{pip install ddfem}. The results shown here were produced using \dunefem{} as the backend. To test the package we recommend setting up a Python virtual environment and using \pyth{pip install ddfem[dune]} to install \dunefem{} together with \ddfem{}. More details on installation and usage can be found on the GitLab page for \ddfem{} \url{https://gitlab.dune-project.org/dune-fem/ddfem} and the tutorial available under \url{https://ddfem.readthedocs.io/en/latest/index.html}.
The software and experiment scripts used for this paper will be tagged as \pyth{v1.1.0} in the repository.

To reproduce the results shown here, bash and Python scripts are included in the git repository and the Python package. On a Unix type operating system, running for example
\begin{pythoncode}
python3 -m venv ddfem-env
source ddfem-env/bin/activate
pip install ddfem[dune]
python -m ddfem getExperiments
cd ddfem-experiments
\end{pythoncode}
will install the package and required dependencies to run all experiments. The required scripts are now located in the \pyth{ddfem-experiments} folder.
\begin{note}
    Depending on your system it is likely that you will have to also install \pyth{mpi4py} using \pyth{pip install mpi4py} before you can run the \pyth{python -m ddfem getExperiments} command. A corresponding error message will be shown if \pyth{mpi4py} is not available on a system with \pyth{MPI} compilers installed.
\end{note}

To run the experiments with reduced resolution compared to those used in this paper, and consequently lower runtime, run
\begin{pythoncode}
bash ./all.sh
\end{pythoncode}
and for full resolution run
\begin{pythoncode}
bash ./all.sh paper
\end{pythoncode}
This will run all experiments and produce all figures from \cref{section:experiments} and the appendices.
To run an individual test using the \pyth{paper.sh} script see the usage instructions:
{\footnotesize
\begin{pythoncode}
usage: paper.sh [-h] problem method [refine] [startRefine]

Run an experiment from the ddfem paper. The first argument gives the
subsection of the experiment, e.g., 3 is the Cahn-Hilliard example from Section 6.3.
For the full set of different boundary conditions for the elliptic problem (section 6.1 and appendix D)
set problem to
       11=Dirichlet, 12=Neumann, 13=Dirichlet outside circle only, 14=Neumann outside circle only
For the twisted beam shown in the conclusion use problem=5 (uses reduced resolution).

positional arguments:
  problem       Experiment from section 6.x with problem=x one of 1, 2, 3, 4
                Use problem=5 from twisted beam and problem=11,12,14 for the results in Appendix D
  method        Method: Fitted, DDM1, Mix0DDM, Mix0DDM1, NSDDM, or reference (only problem=4)
  refine        Refinement level (default: 0)
  startRefine   Start level for error calculations (problem 11,12,13,14, and 4) (default: 0)

options:
  -h, --help    show this help message and exit
\end{pythoncode}
}
After each run, results are located in the folder \pyth{results_x} where \pyth{x} is the number provided as the \pyth{problem} argument. Each folder contains subfolders with the \pyth{vtu} files of each run, these can be opened e.g., with \pyth{paraview}.

\section{Diffuse Domain Method Formulations}
\label{sec:appMethods}

We present mixed boundary versions of the Dirichlet formulations in \citet{Benfield2025b}
for the scalar problem \cref{eq:standardAD}.
We assume a scalar diffusion coefficient $D>0$ and a divergence-free extended velocity $b$.
These assumptions do not restrict the general model \cref{eqn:general}.

Using the boundary partition and extended data $G_V$, $G_{F_v}$, and $G_{F_c}$
from \cref{section:mixedboundaries},
and the boundary jump functions are
\[
    \begin{aligned}
        J_V(u)     & = \sum_{i \in I_V}\omega_i(u-\overline{g_i}),                        \\
        J_{F_v}(F) & = \sum_{i \in I_F}\omega_i(F\cdot n-\overline{g_{v,i}})|\nabla\phi|, \\
        J_{F_c}(F) & = \sum_{i \in I_F}\omega_i(F\cdot n-\overline{g_{c,i}})|\nabla\phi|.
    \end{aligned}
\]
The penalty coefficient $\alpha$ is the sum of the active scaling factors
\pyth{outFactor_e} and \pyth{outFactor_i}
from \cref{section:transformers}.
For prescribed Dirichlet data $g_B$ on $\extdomainbnd{}$, extended to $\extdomain{}$,
take $V = H^1_0(\extdomain)$ and $U = V + g_B$.
In each method below, find $u \in U$ such that the stated equation holds for all $v \in V$.

\begin{meth}[DDM1]
    \label{meth:ddm1adv}
    The weak formulation without stabilisation is
    \begin{equation}
        \begin{aligned}
             & \int_{\extdomain{}}
            \phi D \nabla u \cdot \nabla v
            + \phi(b\cdot\nabla u+cu)v
            \\
             & \qquad
            + \frac{\alpha}{\epsilon^{3}}(1-\phi)J_V(u)v
            - J_{F_c}(bu)v
            \;dx
            \\
             & \qquad =
            \int_{\extdomain{}}(\phi f+G_{F_v})v\;dx.
        \end{aligned}
    \end{equation}
\end{meth}

\begin{meth}[\MixedZ{}]
    \label{meth:mixed0adv}
    Define the gradient by
    \[
        \phi\sigma(u) = \phi\nabla u+J_V(u)\nabla\phi.
    \]
    This gives the weak form,
    \begin{equation}
        \begin{aligned}
             & \int_{\extdomain{}}
            \phi D\sigma(u)\cdot\left[\nabla(\phi v)+v\nabla\phi\right]
            + \phi^2\left[b\cdot\sigma(u)+cu\right]v
            \\
             & \qquad
            + \left[\frac{\alpha}{\epsilon^{2}}(1-\phi)^2
                + \frac{D}{4}|\nabla\phi|^2\right]J_V(u)v
            \\
             & \qquad
            + \phi\left[2J_{F_v}(D\sigma(u))-J_{F_c}(bu)\right]v
            \;dx
            \\
             & \qquad =
            \int_{\extdomain{}}\phi^2 f v\;dx.
        \end{aligned}
    \end{equation}
\end{meth}

\begin{meth}[\NSDDM{}]
    \label{meth:nsddm}
    The Nitsche stabilisation coefficient and the positive inflow term are defined as,
    \[
        \beta = \frac{3D}{2\epsilon}\left[4(1-\phi)+(1-\phi)^2\right],
        \qquad
        [w]^+ =
        \begin{cases}
            w & w > 0,    \\
            0 & w \leq 0.
        \end{cases}
    \]
    This gives the weak form,
    \begin{equation}
        \begin{aligned}
             & \int_{\extdomain{}}
            D\nabla u\cdot\nabla(\phi v)
            + D J_V(u)\nabla\phi\cdot\nabla v
            \\
             & \qquad
            + \phi(b\cdot\nabla u+cu)v
            + \left[J_{F_v}(D\nabla u)-J_{F_c}(bu)\right]v
            \\
             & \qquad
            + \left[\frac{\alpha}{\epsilon^{2}}(1-\phi)
                + \beta|\nabla\phi|+{[b\cdot\nabla\phi]}^+\right]J_V(u)v
            \;dx
            \\
             & \qquad =
            \int_{\extdomain{}}\phi f v\;dx.
        \end{aligned}
    \end{equation}
\end{meth}

\section{Linear Elasticity Model}
\label{app:LEM}
\begin{listing}[H] \footnotesize
    \begin{pythoncode}
class LEModel:
    dimRange = 2
    E, nu = 2e11, 0.3  # generic steel with 2D plane-stress
    rho, g = Constant(7800, name="density"), Constant(9.81, name="gravity")
    # Plane-stress Lame parameters
    lamb = Constant( (E * nu) / (1 - nu**2), name="Lame_1") # = 69 230 769 230.76923 Pa
    mu = Constant( E / (2 * (1 + nu)), name="Lame_2") # = 76 923 076 923.07692 Pa
    outFactor_i = Constant(lamb+2*mu, "outFactor")    # Dirichlet penalty factor
    vd = Constant([0, 0], "stationary")
    vn = Constant([0, 0], "traction")
    I = Identity(dimRange)
    def sigma(U, DU):
        return LEModel.lamb * tr(sym(DU)) * LEModel.I + 2 * LEModel.mu * sym(DU)
    def F_v(t, x, U, DU):
        return LEModel.sigma(U, DU)
    def S_i(t, x, U, DU):
        return as_vector([0.0, -LEModel.rho * LEModel.g])
    def von_Mises(U, DU):
        sigma = LEModel.sigma(U, DU)
        sxx, syy = sigma[0, 0], sigma[1, 1]
        sxy = 0.5 * (sigma[0, 1] + sigma[1, 0])
        return sqrt( sxx**2 - sxx * syy + syy**2 + 3 * sxy**2 )
    valD = lambda t, x, U: LEModel.vd
    valN = lambda t, x, U, DU, n: LEModel.vn
    boundary = {
        "left": BndValue(valD),      # left boundary
        "right": BndFlux_v(valN),    # right boundary
        "bulk": BndFlux_v(valN),     # remaining boundary
    }
\end{pythoncode}
    \caption{Listing for linear elasticity model class}
    \label{lst:LE}
\end{listing}

\section{Blended Diffuse Domain Method (\MixedZDDM)}
\label{sec:appMixedZDDM}

\begin{listing}[H] \footnotesize
    \begin{pythoncode}
from ddfem.transformer import DDM1, Mix0DDM
from ddfem.boundary import BndFlux_v, BndFlux_c, BndValue, BoundaryTerms

def blended(OriginalModel, DDMDirichlet, DDMNeumann, name, domainDescription, EXTRA_BOUNDARY=True):
    # get original model - the model provided as first argument is pre-transformed
    class DirichletModel(OriginalModel):
        boundary = { k: v for k, v in OriginalModel.boundary.items()
                     if not isinstance(v, (BndFlux_v, BndFlux_c)) }
    class NeumannModel(OriginalModel):
        boundary = { k: v for k, v in OriginalModel.boundary.items()
                     if not isinstance(v, BndValue) }
    # now transform the models by the given transformers
    DDMDirichletModel = DDMDirichlet(DirichletModel, domainDescription, EXTRA_BOUNDARY)
    DDMNeumannModel = DDMNeumann(NeumannModel, domainDescription, EXTRA_BOUNDARY)
    # construct the blended model
    wD = BoundaryTerms(OriginalModel, domainDescription).dbc_weight
    class DDModel(DDMDirichletModel):
        wD_ = wD
        transformerName = name
        domainDesc = domainDescription
        BT = BoundaryTerms(OriginalModel, domainDescription)
        boundary = BT.physical
        boundary_diffuse = BT.diffuse
        domain = BT.domain
        # if a method only exists on DM then take is inherited so we only
        # need to check if it exists on NM and then blend
        # Shown here for S_i,F_v - same code needed for S_e,F_c
        if hasattr(DDMNeumannModel, "S_i"):
            if hasattr(DDMDirichletModel, "S_i"):
                def S_i(t, x, U, DU):
                    return ( wD(x) * DDMDirichletModel.S_i(t,x,U,DU) +
                            (1-wD(x)) * DDMNeumannModel.S_i(t,x,U,DU) )
            else:
                S_i = DDMNeumannModel.S_i
        if hasattr(DDMNeumannModel, "F_v"):
            if hasattr(DDMDirichletModel, "F_v"):
                def F_v(t, x, U, DU):
                    return ( wD(x) * DDMDirichletModel.F_v(t,x,U,DU) +
                            (1-wD(x)) * DDMNeumannModel.F_v(t,x,U,DU) )
            else:
                F_v = DDMNeumannModel.F_v
        def sigma(t, x, U, DU=None):
            return ( wD(x) * DDMDirichletModel.sigma(t,x,U,DU) +
                   (1-wD(x)) * DDMNeumannModel.sigma(t,x,U,DU) )
    return DDModel

def Mix0DDM1(OriginalModel, domainDescription, EXTRA_BOUNDARY=True):
    return blended(OriginalModel, Mix0DDM, DDM1, "Mix0DDM1", domainDescription, EXTRA_BOUNDARY)
    \end{pythoncode}
    \caption{This transformer blends Mix0DDM and DDM1 between Dirichlet and Neumann boundary conditions - for brevity only the implicit source and diffusive flux are shown.}
\end{listing}

\section{Problem from \cref{sec:poisson} with different boundary conditions}
\label{app:poisson}
\begin{figure}[htbp]
    \begin{subfigure}[t]{\textwidth}
        \centering
        \includegraphics[width=0.4\textwidth]{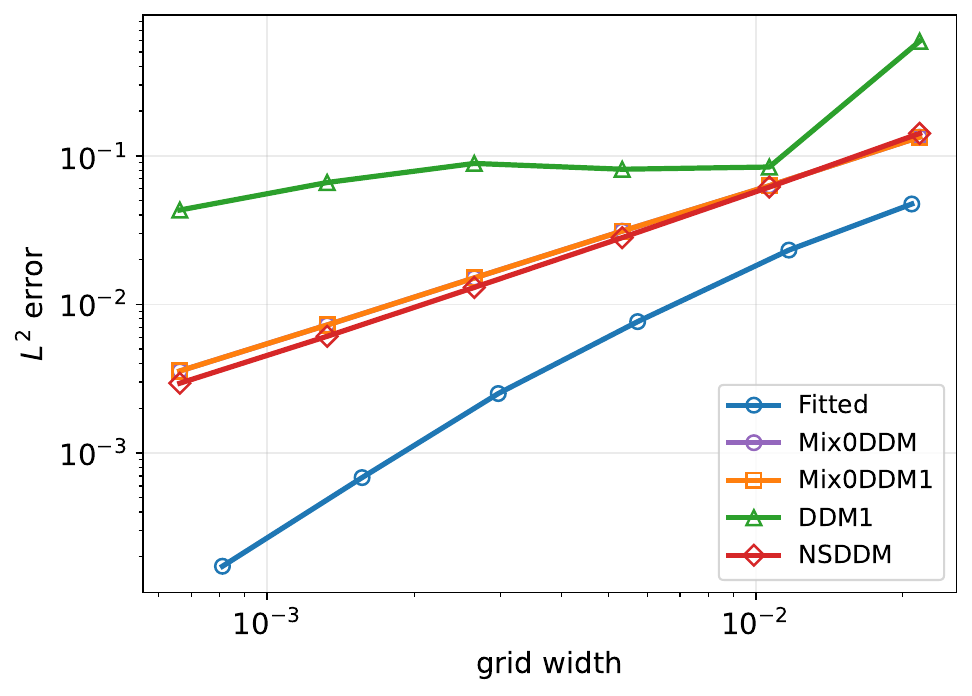}
        \includegraphics[width=0.4\textwidth]{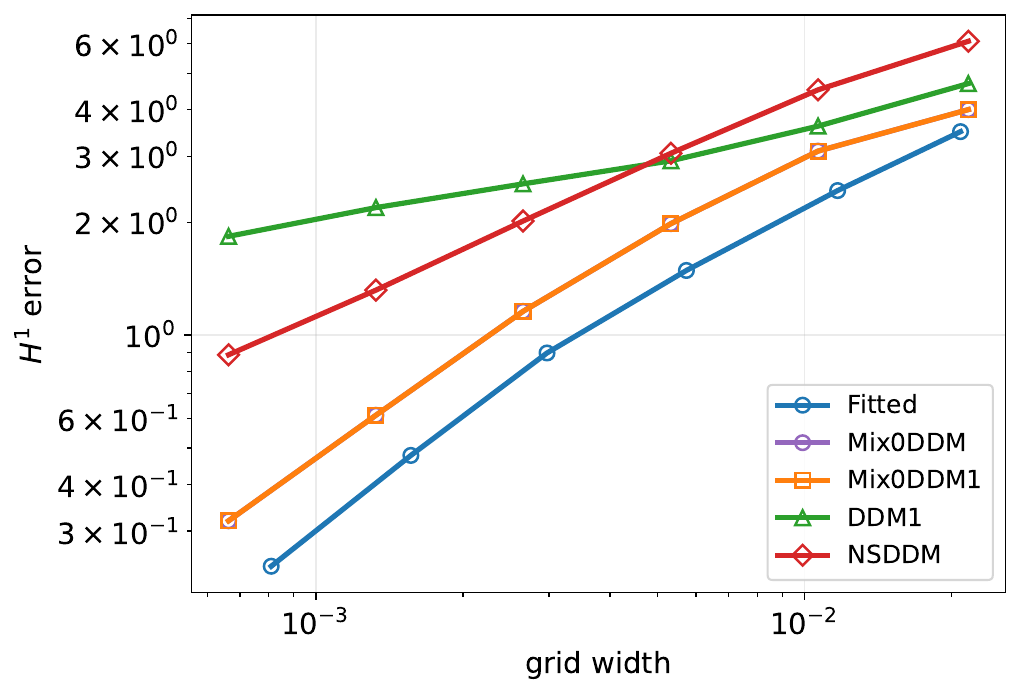}
        \caption{Comparing method with all Dirichlet boundary conditions}
        \label{fig:results1DD}
    \end{subfigure}
    \begin{subfigure}[t]{\textwidth}
        \centering
        \includegraphics[width=0.4\textwidth]{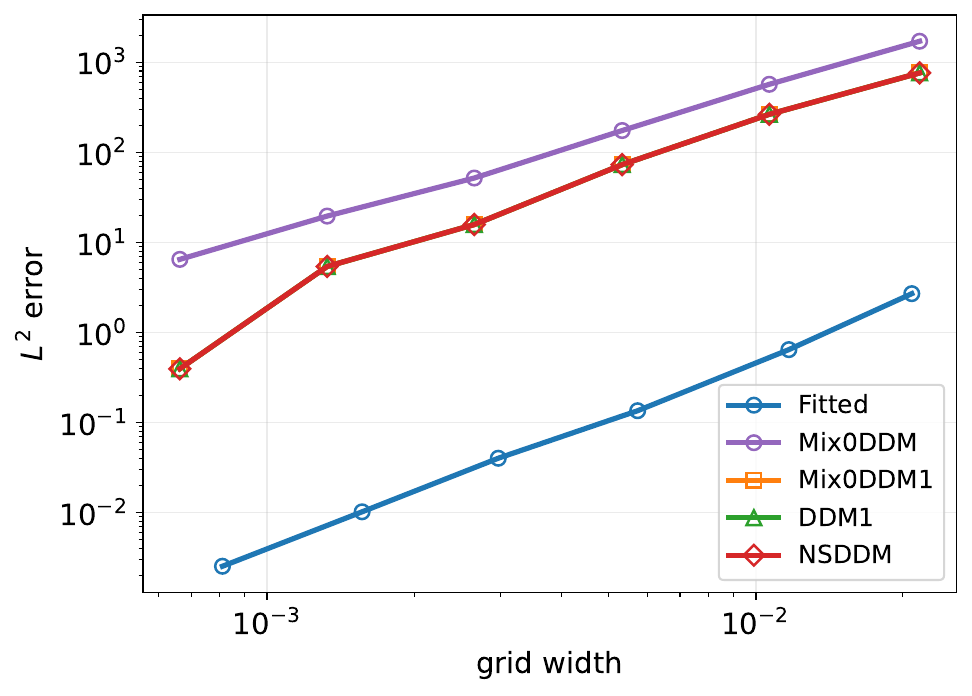}
        \includegraphics[width=0.4\textwidth]{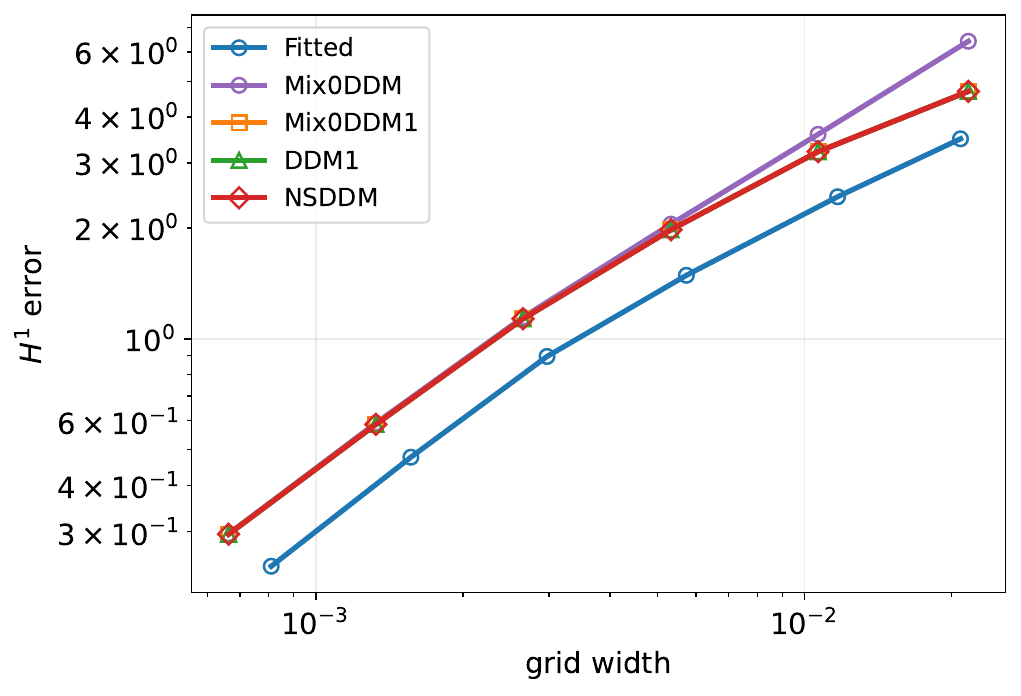}
        \caption{Comparing method with all Neumann boundary conditions}
        \label{fig:results1NN}
    \end{subfigure}
    \begin{subfigure}[t]{\textwidth}
        \centering
        \includegraphics[width=0.4\textwidth]{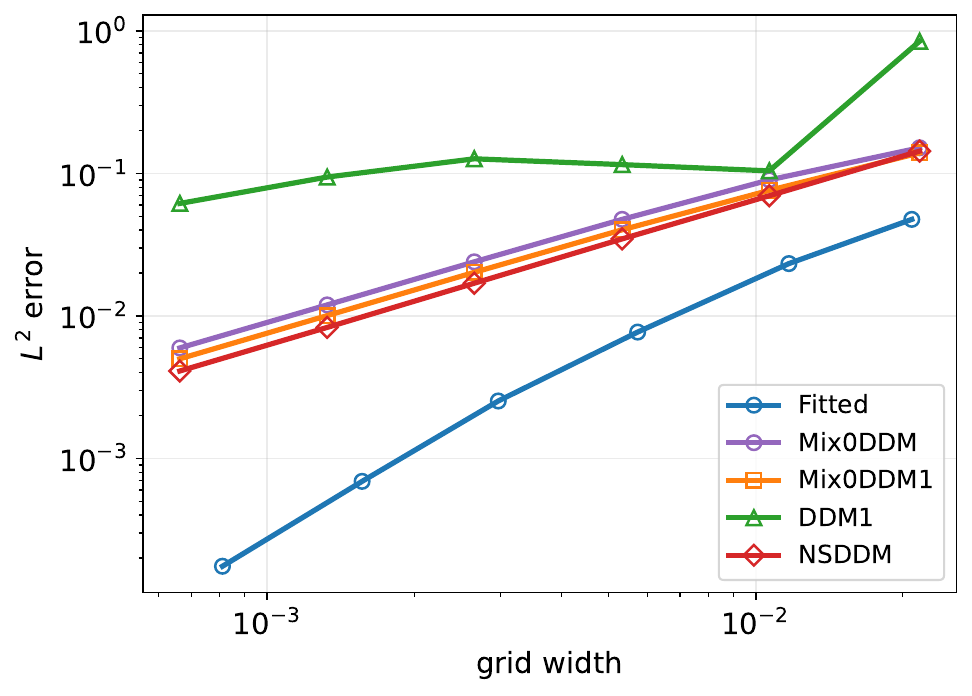}
        \includegraphics[width=0.4\textwidth]{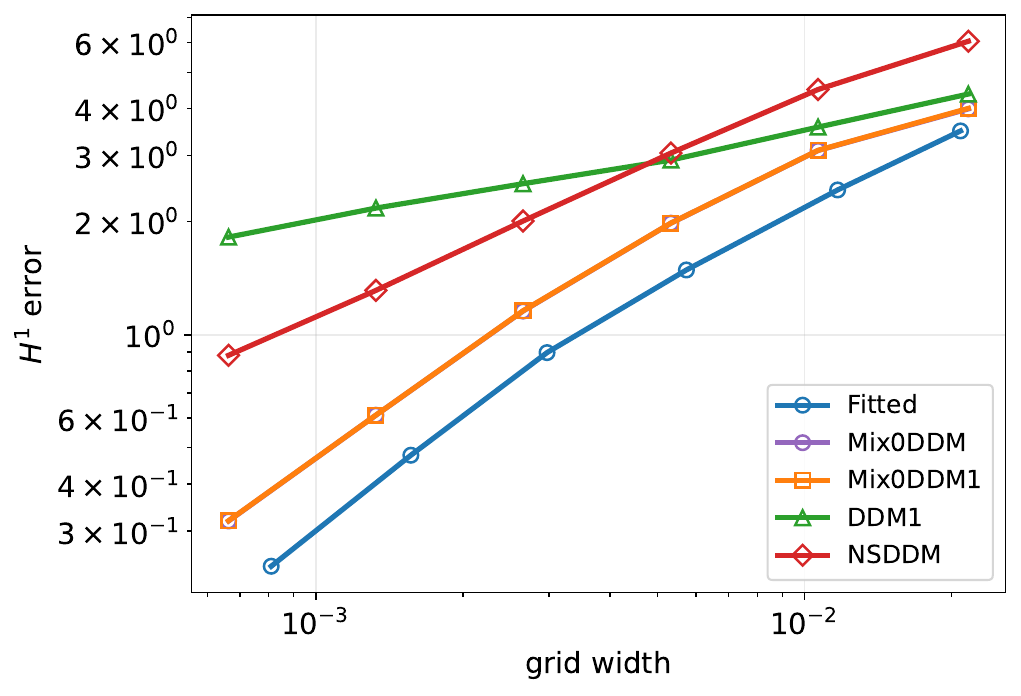}
        \caption{Comparing method with Dirichlet on outside and Neumann on inside circle}
        \label{fig:app:results1DN}
    \end{subfigure}
    \begin{subfigure}[t]{\textwidth}
        \centering
        \includegraphics[width=0.4\textwidth]{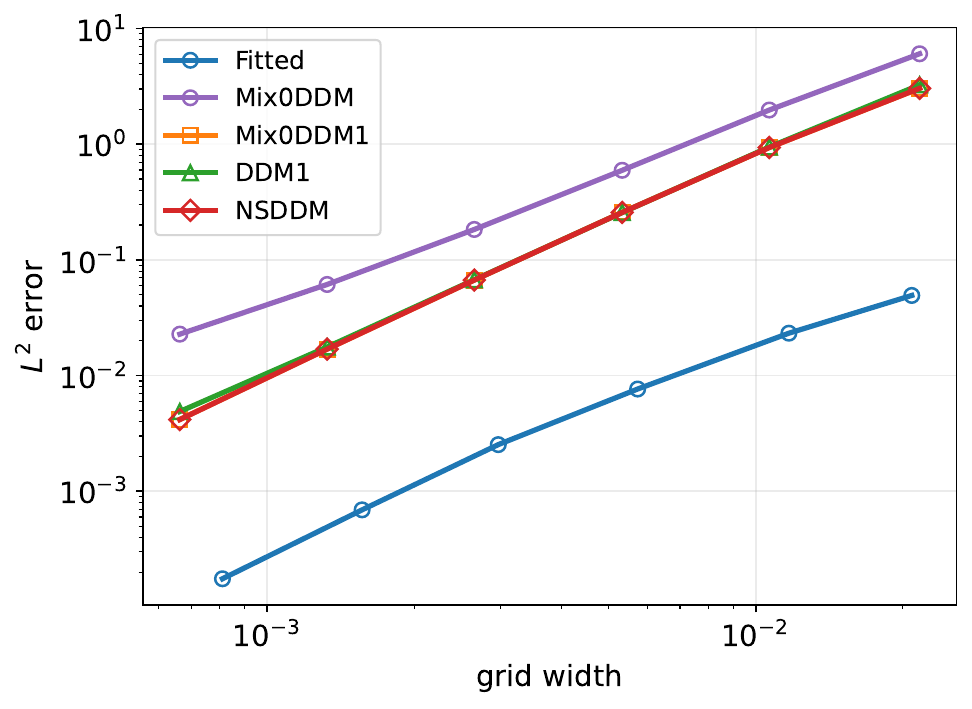}
        \includegraphics[width=0.4\textwidth]{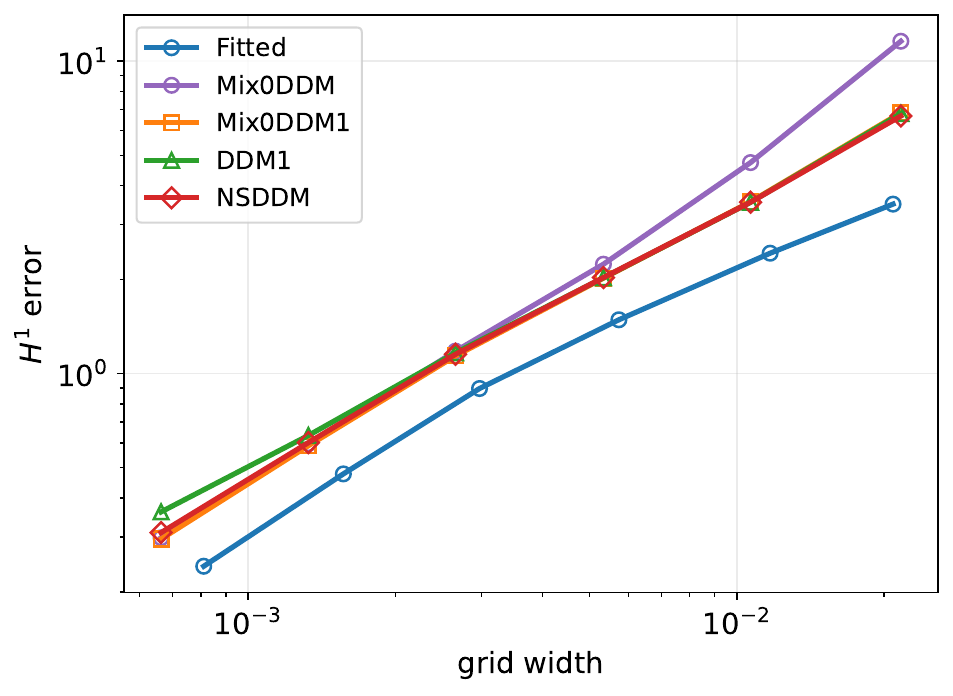}
        \caption{Comparing method with Dirichlet on inside and Neumann on outside circle}
        \label{fig:results1ND}
    \end{subfigure}
    \caption{Problem from \cref{sec:poisson} comparing errors ($L^2$ on left and $H^1$ on right). Domain is the
        difference of two circles.}
\end{figure}

\section{Problem from \cref{sec:cahnh}: mass conservation}
\label{app:ch}
\begin{figure}[htbp]
    \centering
    \includegraphics[width=0.9\textwidth]{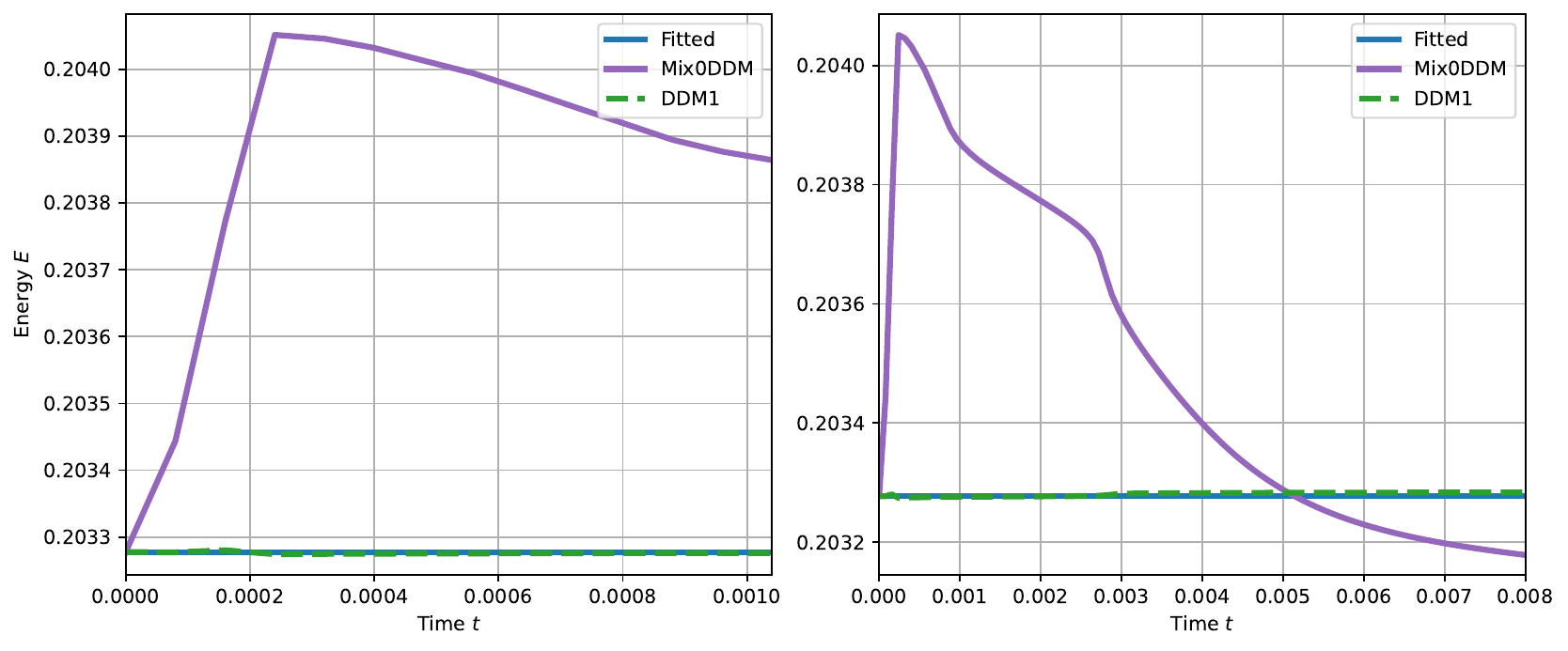}
    \caption{Comparison of computed mass for Cahn-Hilliard problem for \DD{} methods and fitted mesh.}
    \label{fig:chMass}
\end{figure}

\section{Problem from \cref{sec:linelast}: solution along the boundary}
\label{app:LEBnd}
\begin{figure}
    \centering
    \centering
    \includegraphics[width=0.3\textwidth]{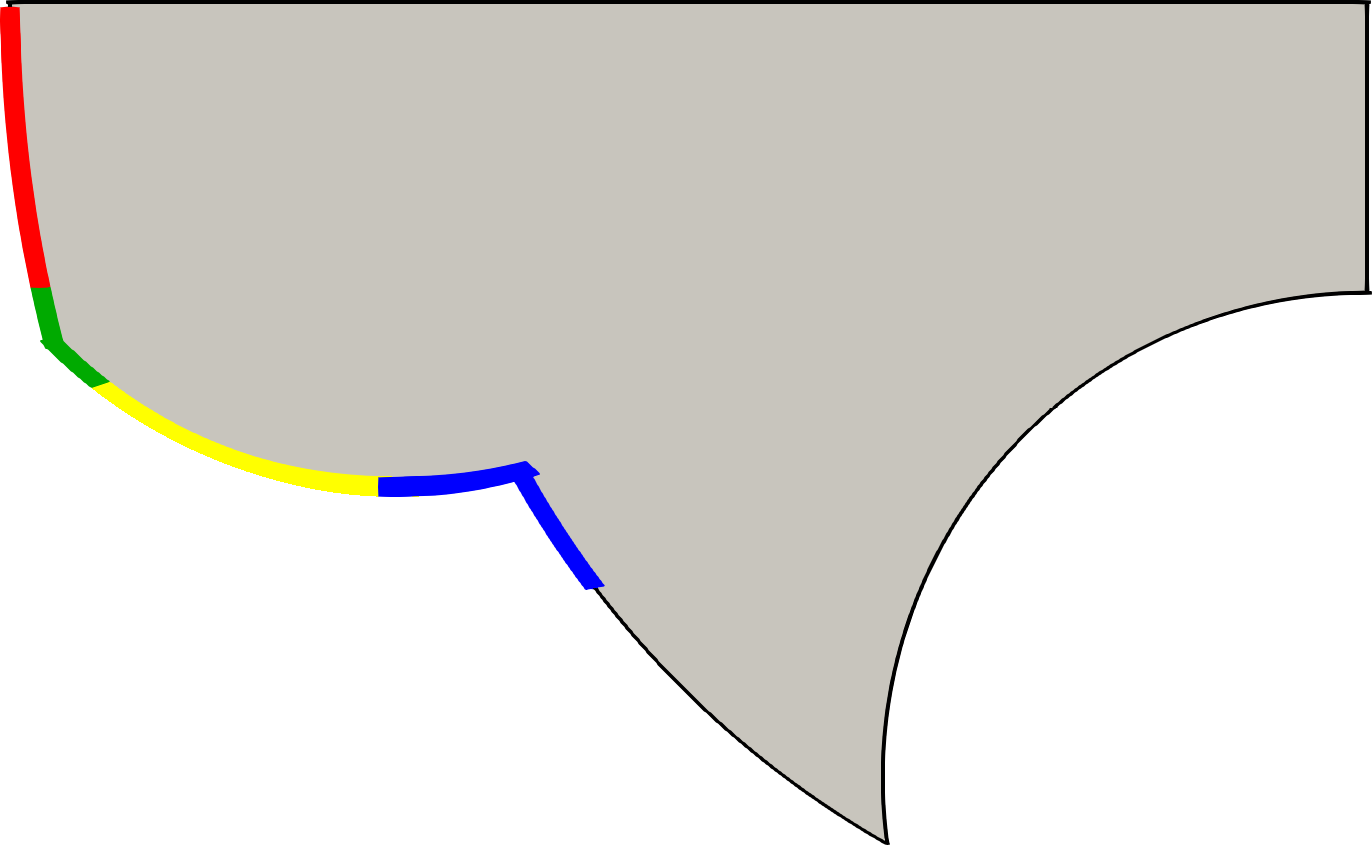}
    \caption{Location of the four boundary segments A (red), B (green), C (yellow), D (blue) along the lower left part of the domain.}
    \label{fig:LEBndMarks}
\end{figure}

\begin{figure}
    \includegraphics[width=0.24\textwidth]{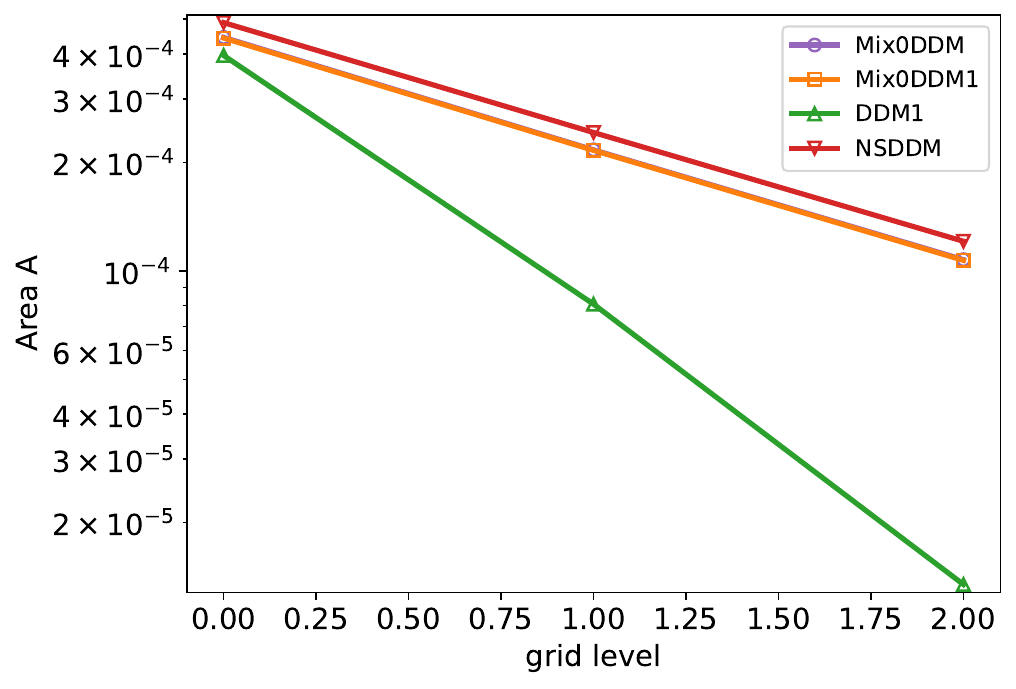}
    \includegraphics[width=0.24\textwidth]{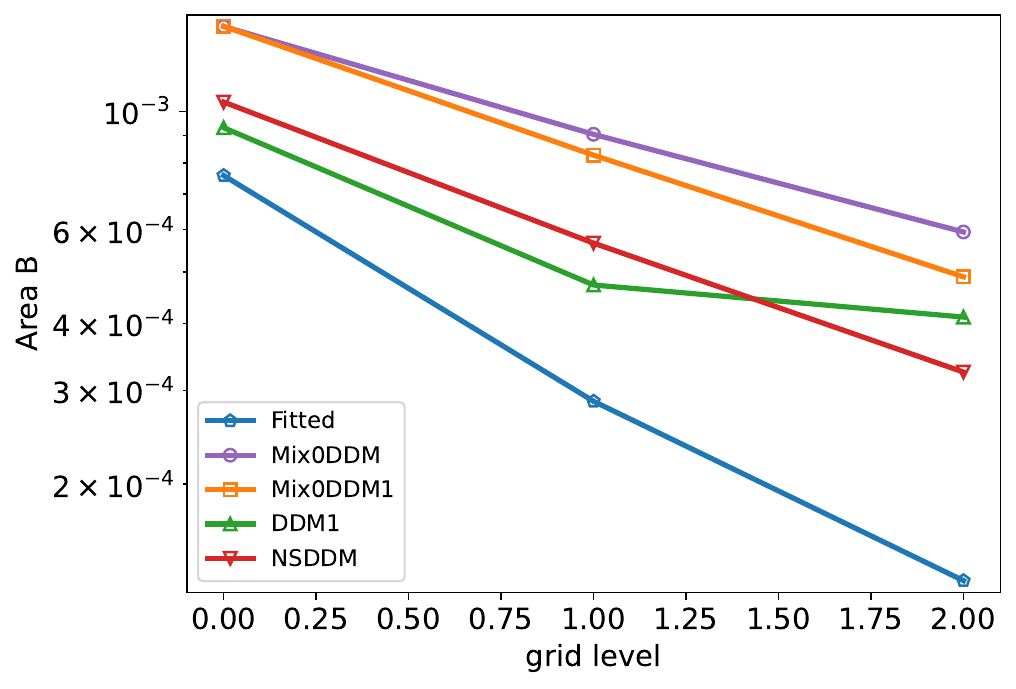}
    \includegraphics[width=0.24\textwidth]{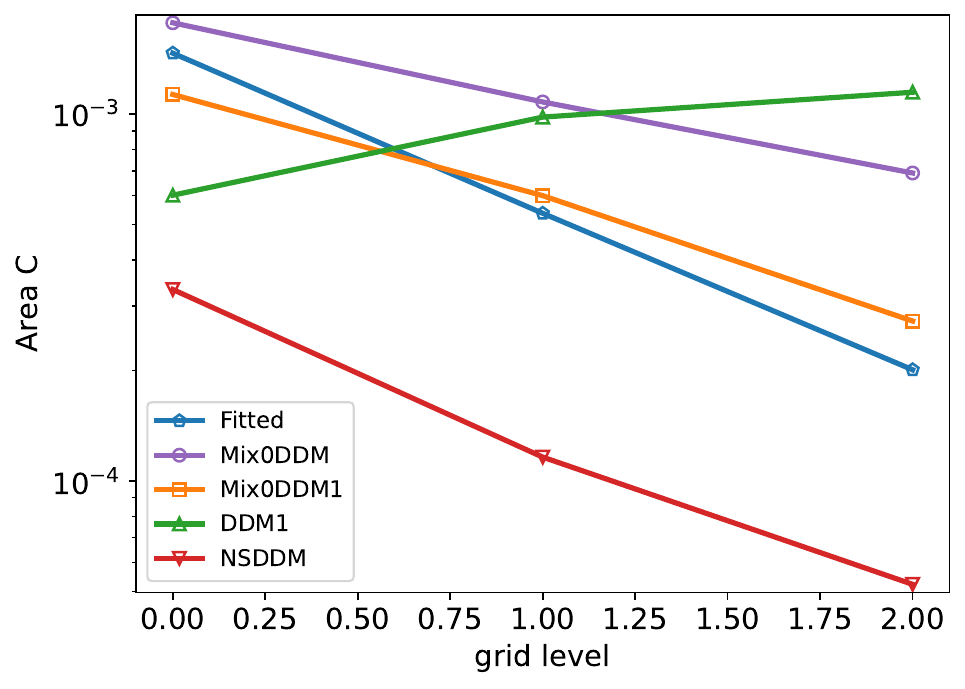}
    \includegraphics[width=0.24\textwidth]{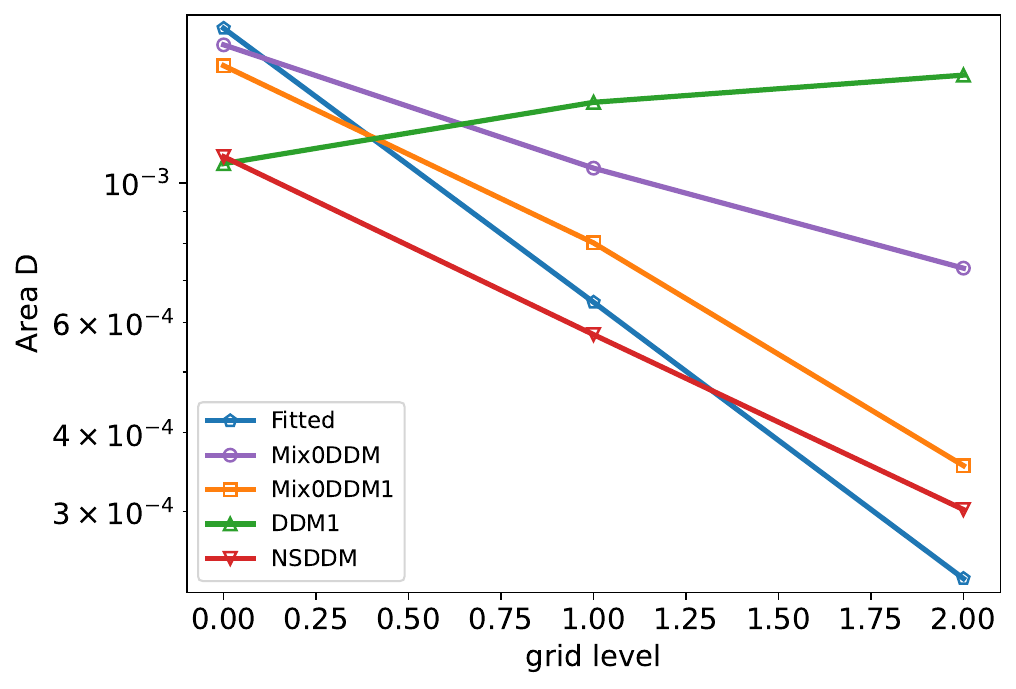}
    \\[1em]
    \includegraphics[width=0.24\textwidth]{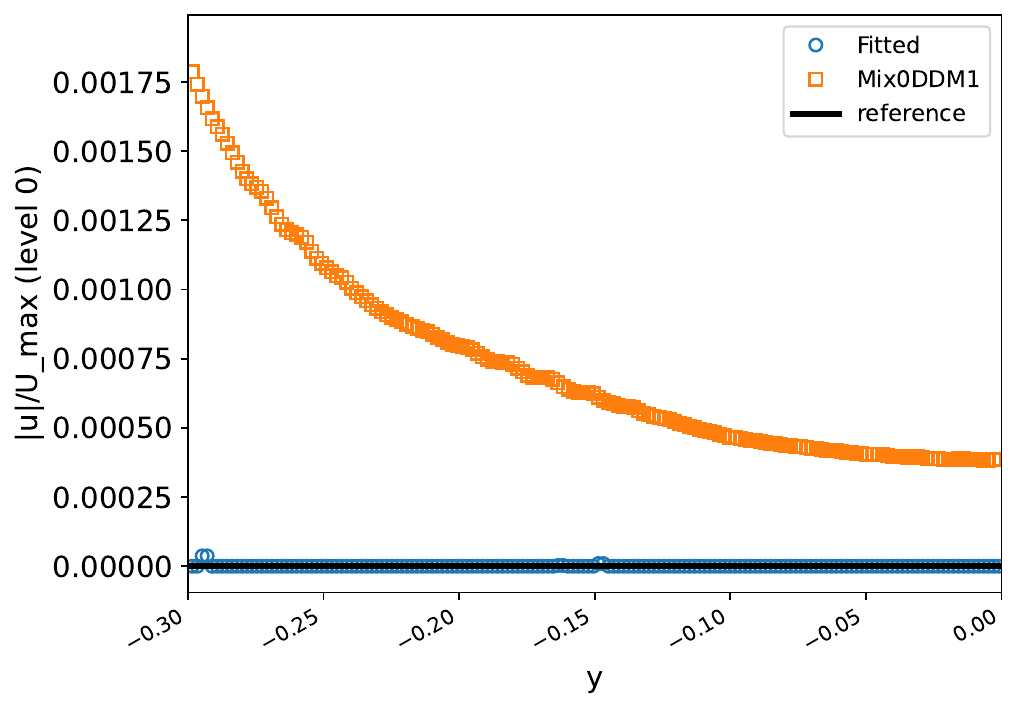}
    \includegraphics[width=0.24\textwidth]{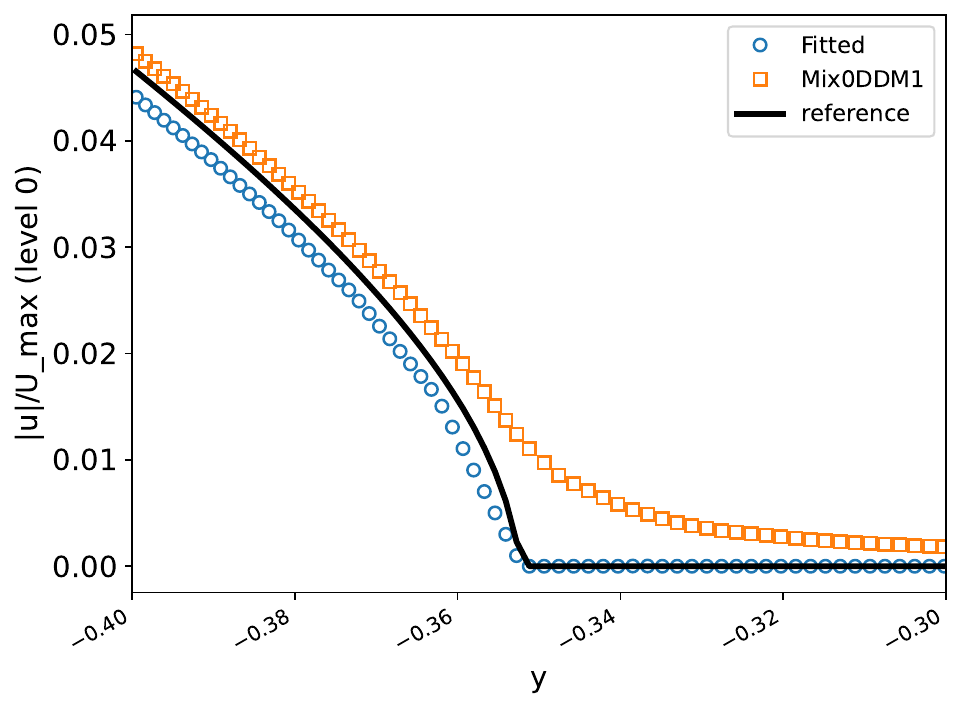}
    \includegraphics[width=0.24\textwidth]{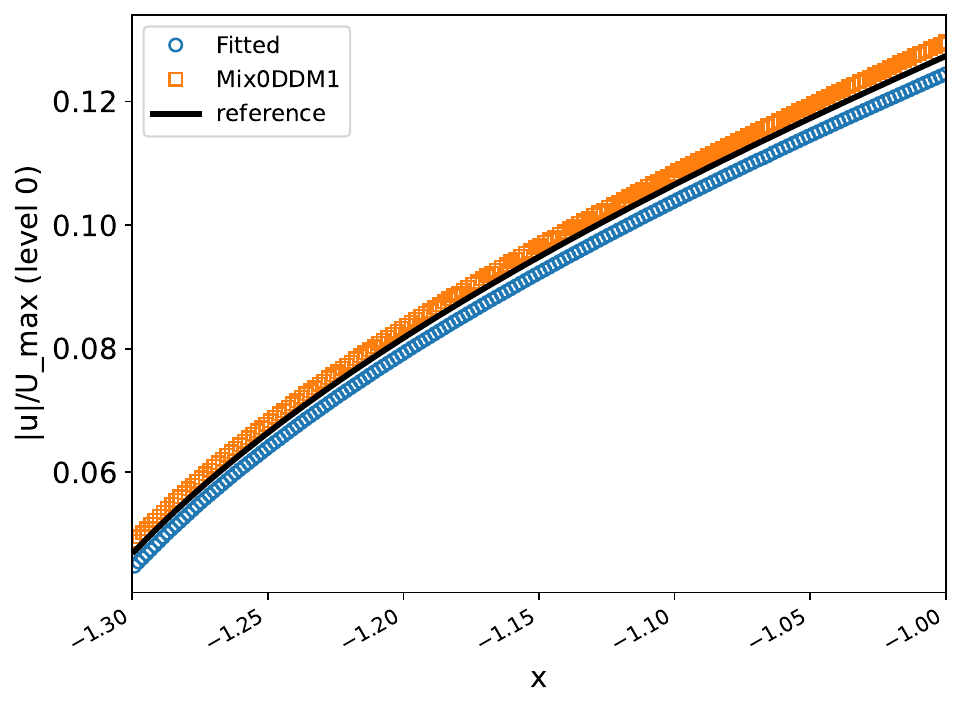}
    \includegraphics[width=0.24\textwidth]{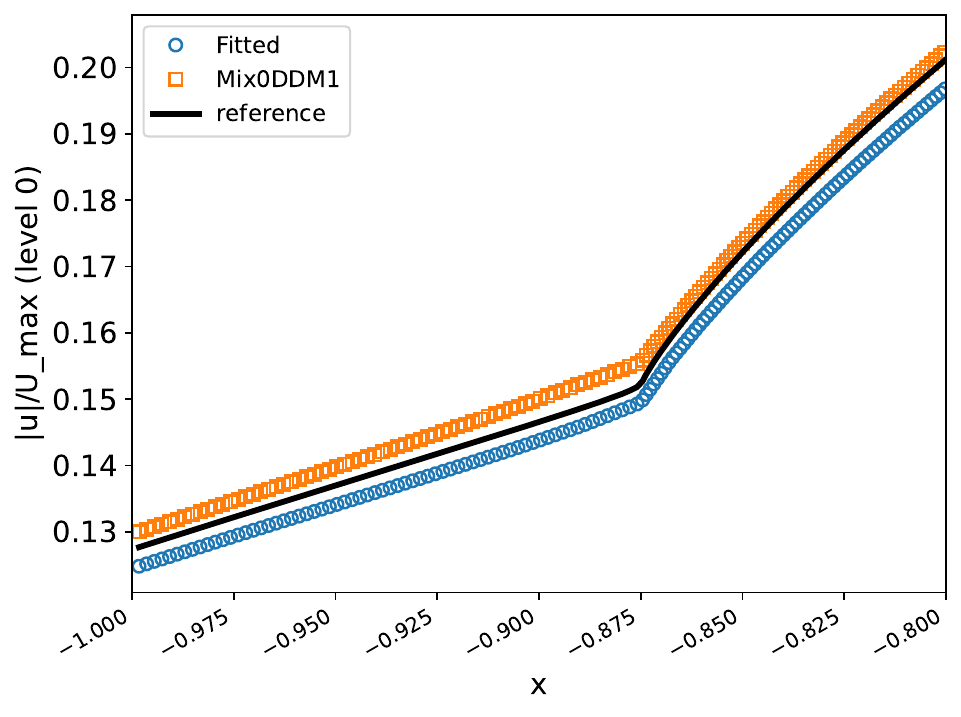} \\
    \includegraphics[width=0.24\textwidth]{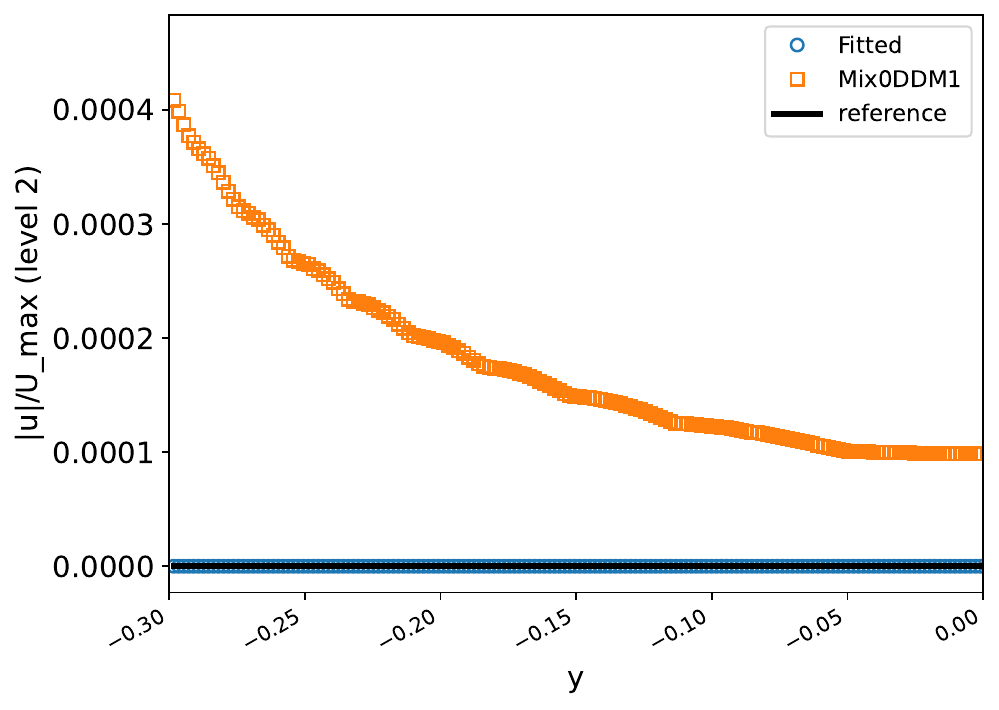}
    \includegraphics[width=0.24\textwidth]{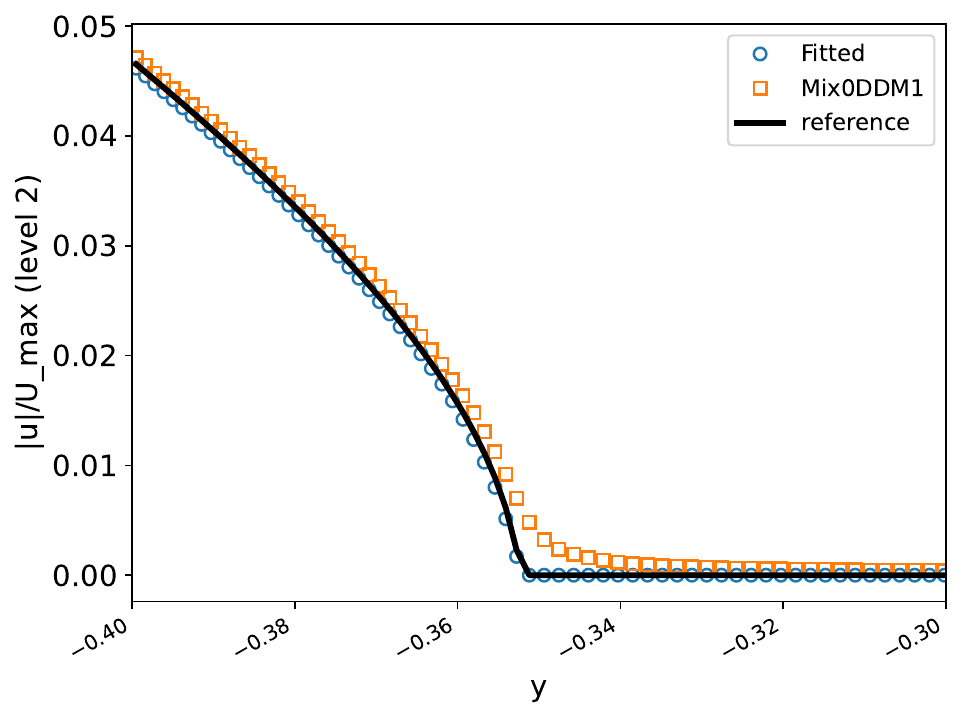}
    \includegraphics[width=0.24\textwidth]{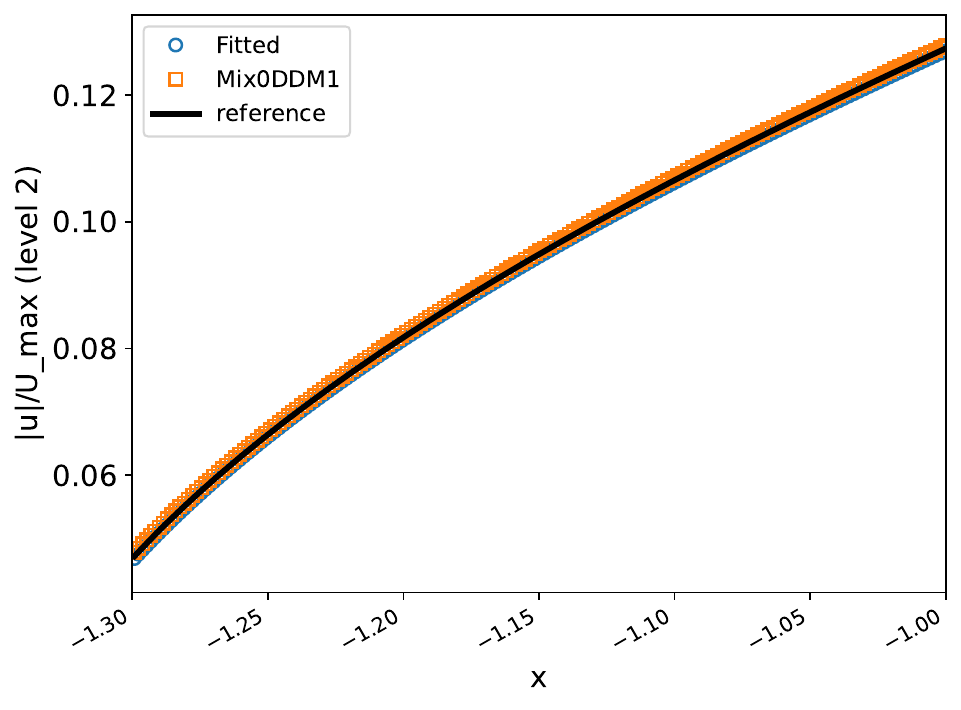}
    \includegraphics[width=0.24\textwidth]{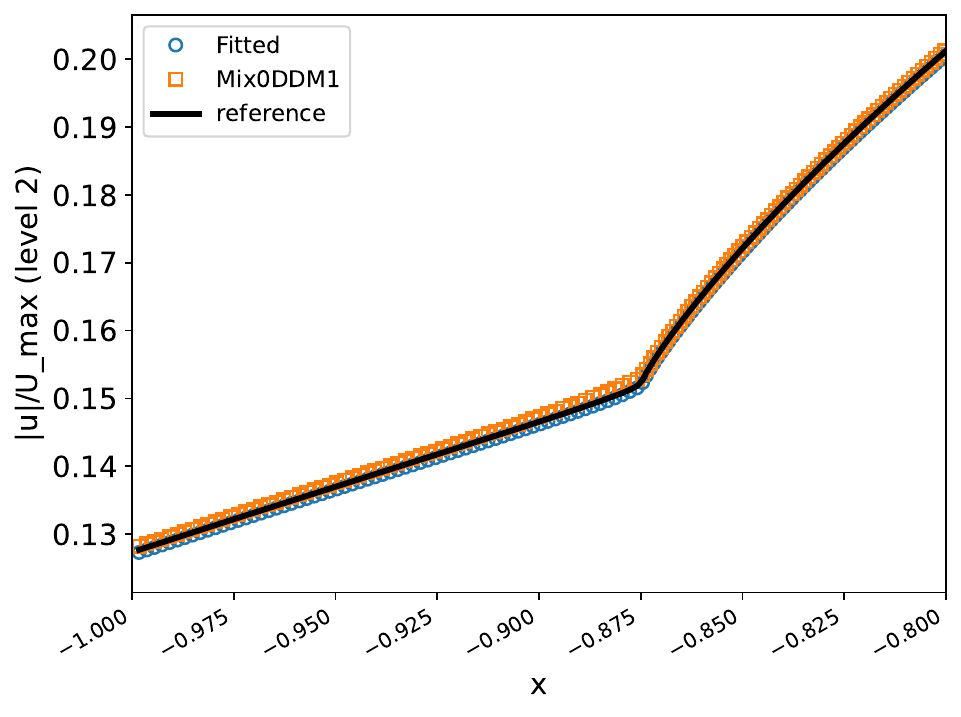}
    \caption{Regions A,B,C,D from left to right: integrated displacement errors (top row) and displacement magnitudes along the boundary segments for reference, fitted, and \MixedZDDM{} on the level zero and two grids (middle and bottom). Displacement normalised with maximum displacement of reference solution $U=6\cdot 10^{-5}$.}
    \label{fig:LEBndD}
\end{figure}

\begin{figure}
    \includegraphics[width=0.24\textwidth]{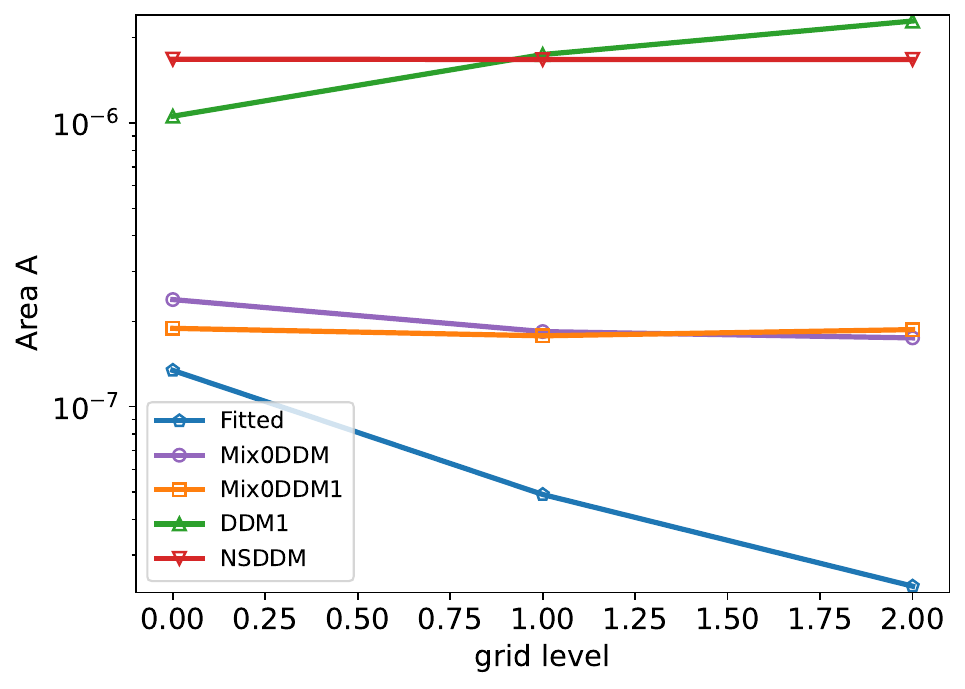}
    \includegraphics[width=0.24\textwidth]{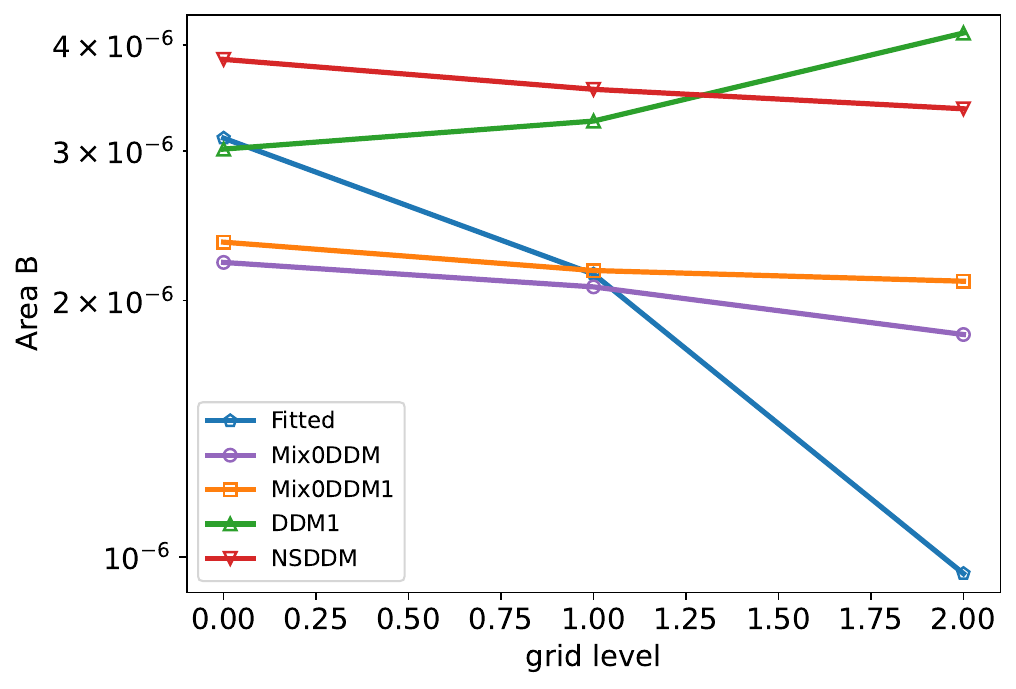}
    \includegraphics[width=0.24\textwidth]{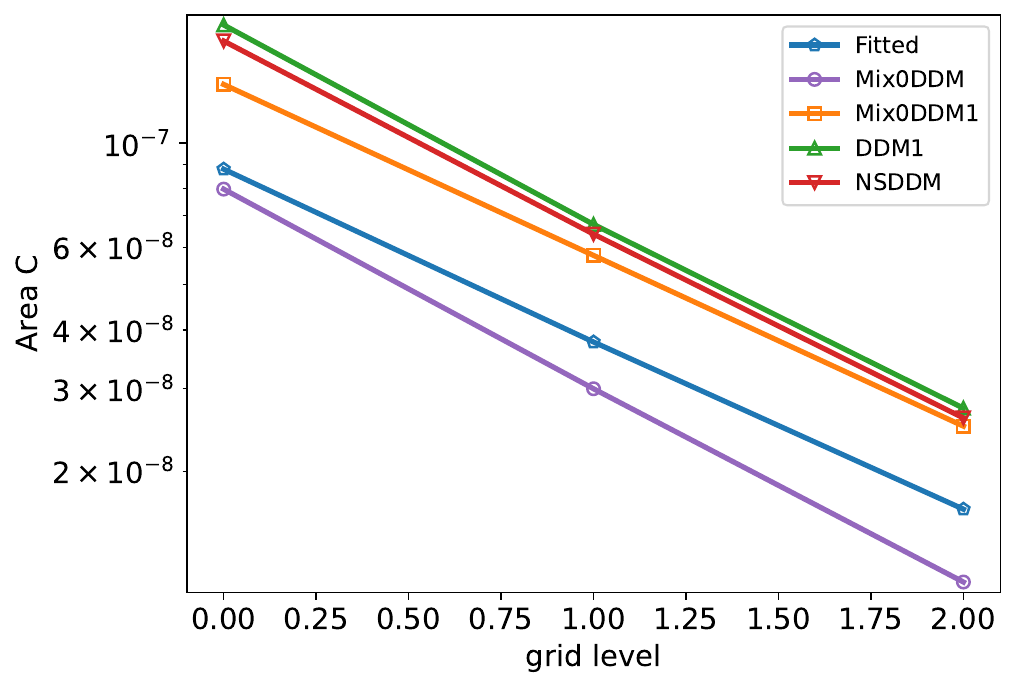}
    \includegraphics[width=0.24\textwidth]{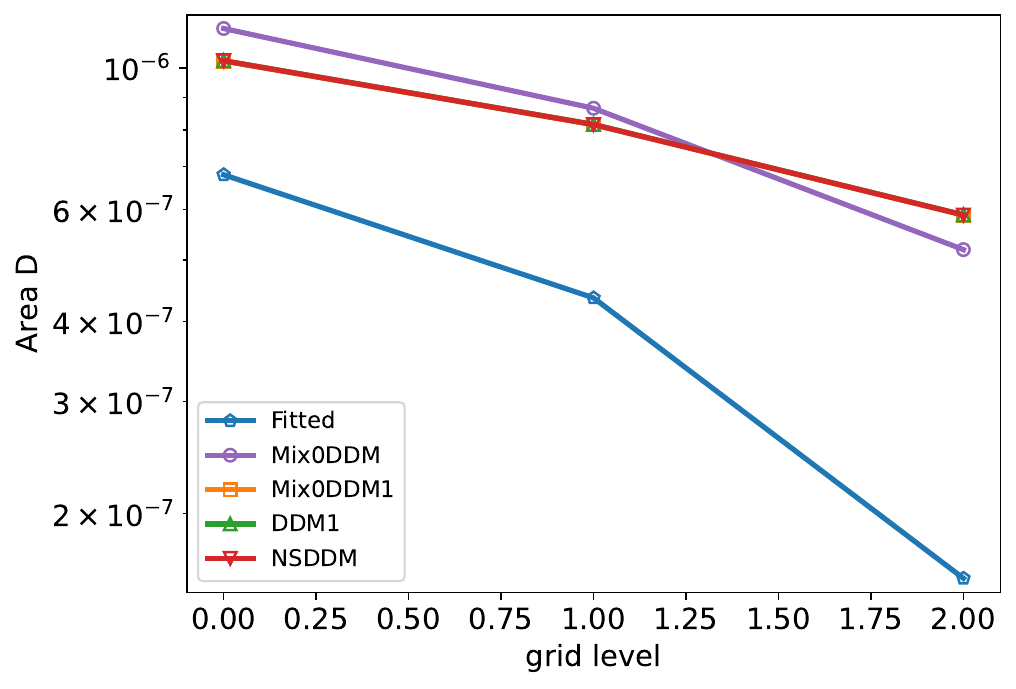}
    \\[1em]
    \includegraphics[width=0.24\textwidth]{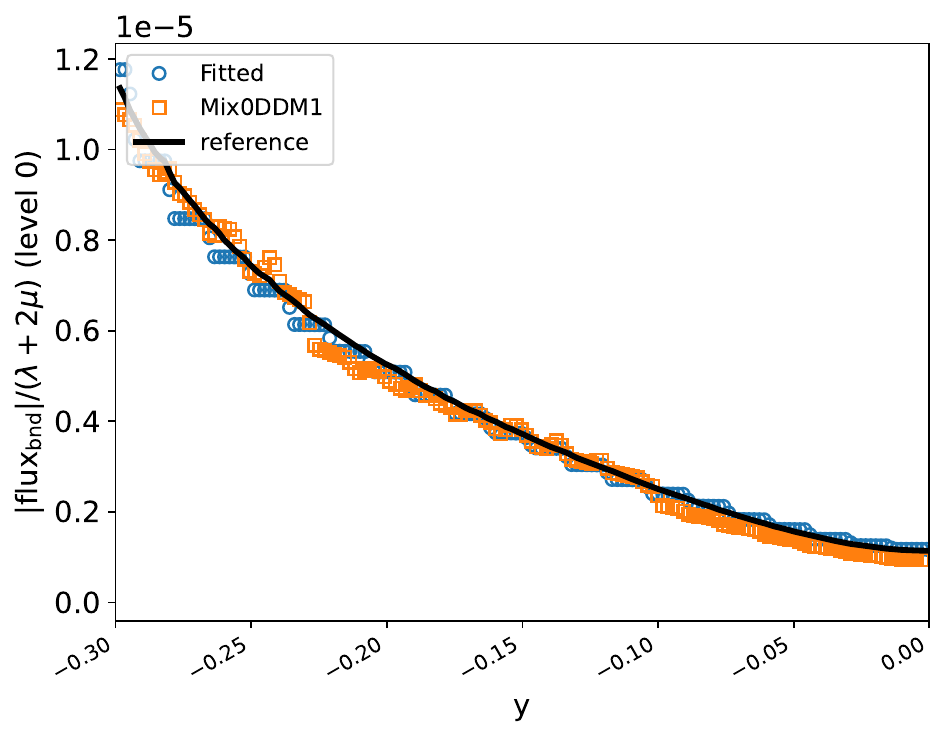}
    \includegraphics[width=0.24\textwidth]{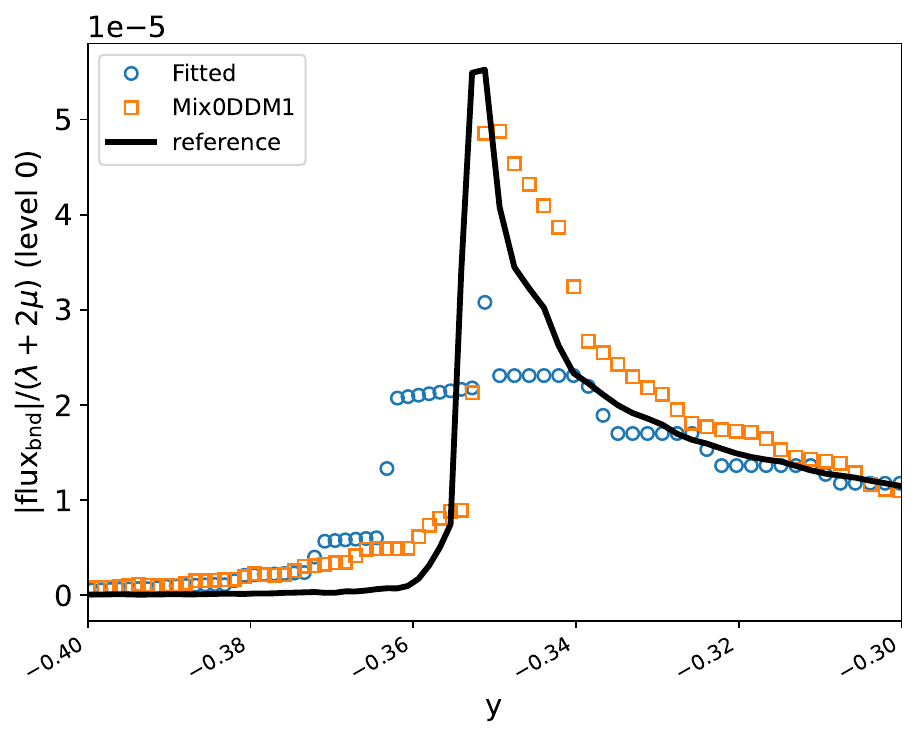}
    \includegraphics[width=0.24\textwidth]{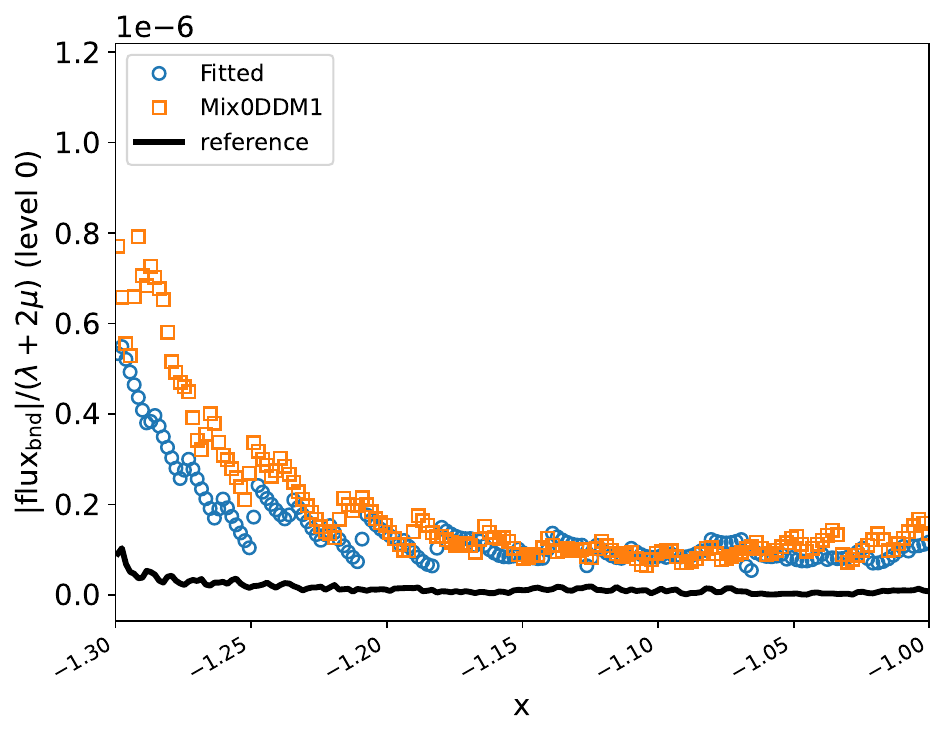}
    \includegraphics[width=0.24\textwidth]{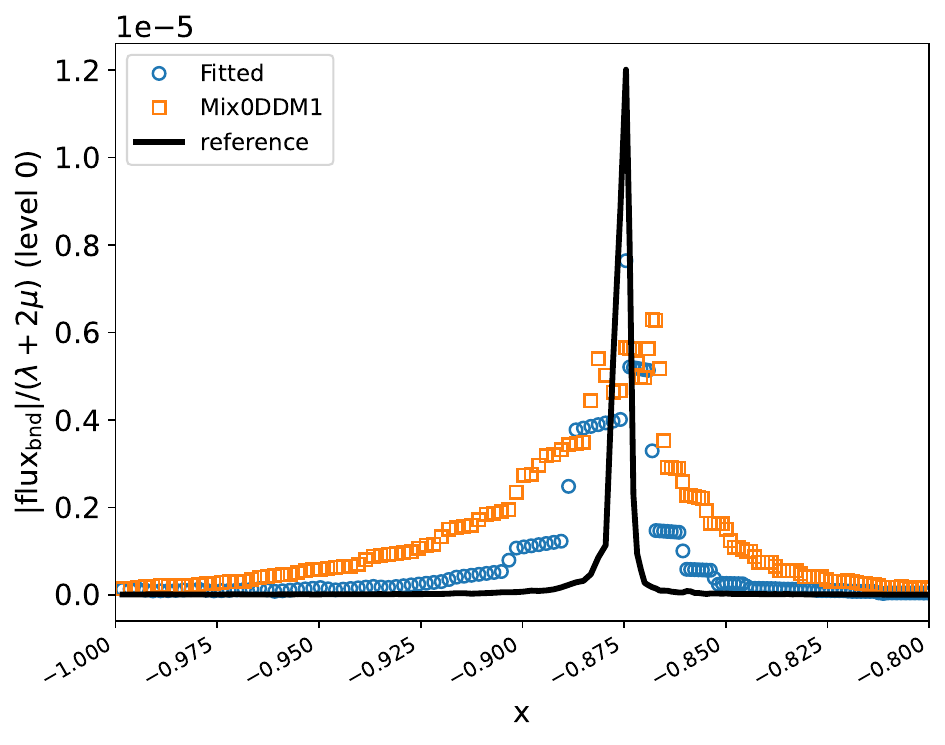} \\
    \includegraphics[width=0.24\textwidth]{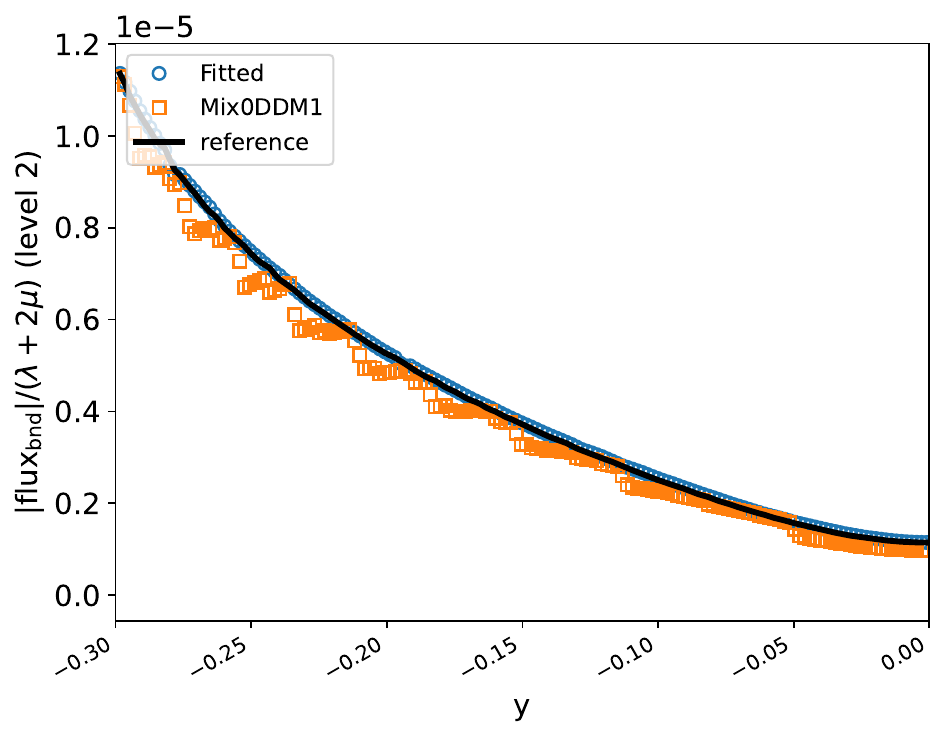}
    \includegraphics[width=0.24\textwidth]{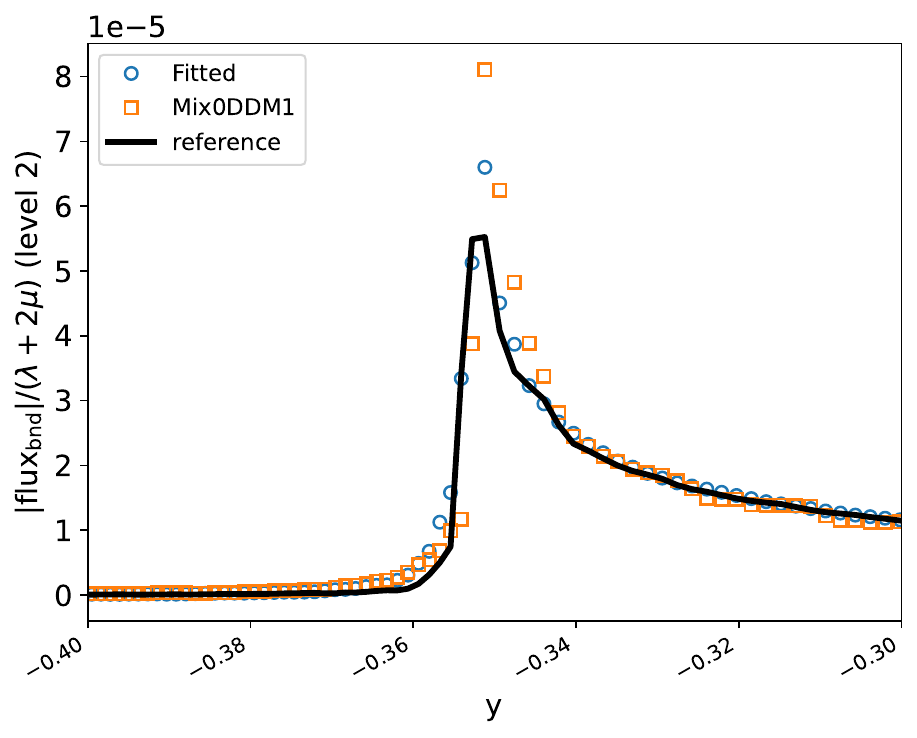}
    \includegraphics[width=0.24\textwidth]{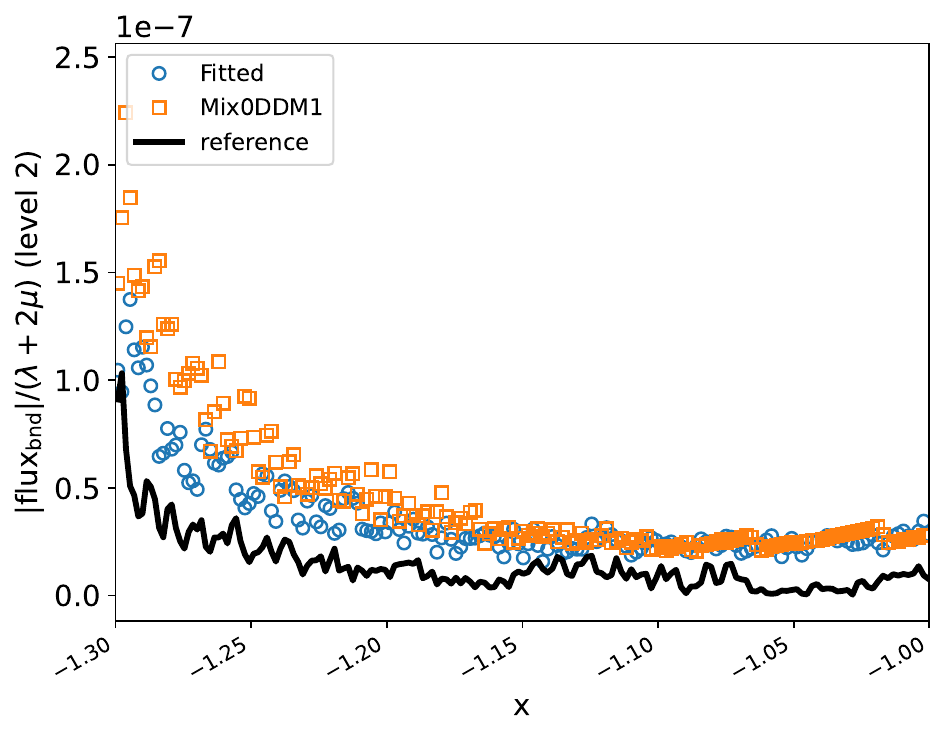}
    \includegraphics[width=0.24\textwidth]{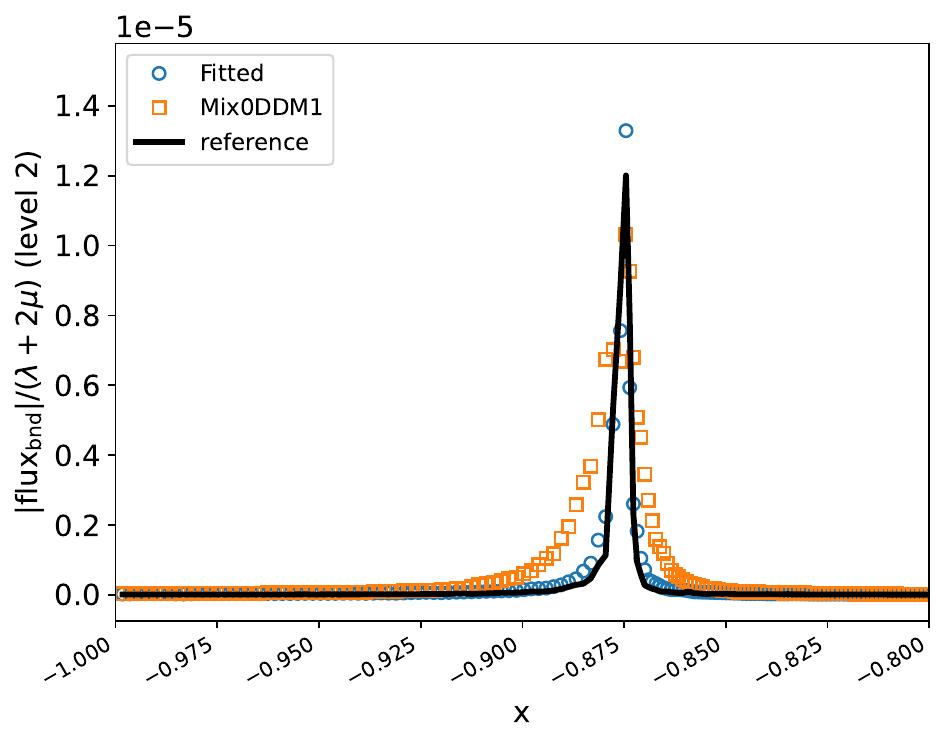}
    \caption{Regions A,B,C,D from left to right: integrated traction errors for all methods (top row) and actual magnitude of the traction in normal direction $|\sigma n|$ along the boundary segments for reference, fitted, and \MixedZDDM{} on the level zero and two grids (middle and bottom) normalised by $\lambda+2\mu$.}
    \label{fig:LEBndF}
\end{figure}

We compare the displacement magnitude and the magnitude of normal traction along different sections of the boundary. All sections are on the lower left as shown in \cref{fig:LEBndMarks}. In \cref{fig:LEBndD} we compare the displacement along these four segments. The top row shows the $L^2$ norm of the difference between the reference displacement and the displacements computed with the \DD{} methods. Displacements are scaled with the maximum displacement in the reference solution $U=6\cdot 10^{-5}$. For segment A, which is along the Dirichlet boundary, the error between fitted and reference is negligible and not shown. Otherwise we see again that convergence of the simplest method \DDMO{} is not always guaranteed - although in region A it outperforms the other \DD{} method by quite some margin. Otherwise results are similar to what we already reported with \MixedZDDM{} and \NSDDM{} being the two methods of choice, outperforming the simulations on the coarser fitted grids in regions C and D (along the zero traction part of the boundary). The second row shows the actual displacement magnitudes along these segments for only fitted and \MixedZDDM{}. Along the Dirichlet boundary we do not see $u=0$ with \MixedZDDM{} but values are small - being only $0.07\%$ of the maximum displacement $U$. On the coarser grid this displacement grows to $0.14\%$ and on the finer grid $0.035\%$ (not shown here) indicating linear convergence to zero.
In \cref{fig:LEBndF} we show the same set of results but now for the magnitude of the normal traction $\sigma\cdot n$ - we again scale results by $\lambda+2\mu\approx 2.6923\cdot 10^{11}$. Along segment A, all \DD{} methods show issues with convergence, although the magnitude of the error is already quite small so the interpolation errors could be influencing the results, however the fitted method seems to be converging. Similar issues can be seen along area B, which like area A, contains part of the Dirichlet boundary. The values plotted along the boundary in the second row look close to the reference values. Along the two regions along the zero traction boundary, all methods seem to converge as expected.
\label{app:linelast}


\end{document}